\documentclass[final,3p]{elsarticle}

\DeclareMathAlphabet{\pazocal}{OMS}{zplm}{m}{n}
\DeclareMathAlphabet{\bpazocal}{OMS}{cmsy}{b}{n}

\RequirePackage[ruled,vlined]{algorithm2e}

\SetCommentSty{xCommentSty}
\RequirePackage{caption,subcaption}

\RequirePackage{lineno}
\RequirePackage{xcolor}
\RequirePackage{amsmath,amsfonts,amsthm,amssymb}
\newtheorem{remark}{Remark}
\RequirePackage[breaklinks]{hyperref}
\RequirePackage{stmaryrd}

\RequirePackage{cleveref}

\RequirePackage{tikz}
\RequirePackage{pgfplots}
\RequirePackage{multicol}\RequirePackage{multirow}

\usetikzlibrary{external}
\usetikzlibrary{arrows.meta}
\usetikzlibrary{shapes.geometric}
\usetikzlibrary{matrix,automata,fit,calc,positioning}

\pgfplotsset{compat=newest}
\usepgfplotslibrary{fillbetween}
\usepgfplotslibrary{groupplots}
\usepgfplotslibrary{external}
\usetikzlibrary{patterns}

\newcommand*{\tikzmk}[1]{\tikz[remember picture,overlay] \node (#1) {};\ignorespaces}
\newcommand{\boxit}[1]{\tikz[remember picture,overlay]{\node[yshift=3pt,fill=#1,opacity=.15,fit={(A)($(B)+(.91\linewidth,.8\baselineskip)$)}] {};}\ignorespaces}

\tikzset{tipA/.tip={Triangle[angle=45:4pt]Bar[width=5.0mm]},
  tipB/.tip={Triangle[angle=45:4pt]} ,
  tipC/.tip={Bar[width=5.0mm]},
  tipD/.tip={Bar[width=2.5mm]},
}
\tikzset{external/only named=true}

\definecolor{blackmy}{RGB}{38, 70, 83}
\definecolor{bluemy}{RGB}{39, 125, 161}
\definecolor{greenmy}{RGB}{42, 167, 143}
\definecolor{yellowmy}{RGB}{233, 196, 106}
\definecolor{brownmy}{RGB}{244, 162, 97}
\definecolor{redmy}{RGB}{249, 65, 68}
\definecolor{cyanmy}{RGB}{5, 255, 255}
\definecolor{purplemy}{RGB}{129, 80, 192}
\tikzset{external/up to date check=md5}\renewcommand{\vec}[1]{\boldsymbol{#1}}

\def \ev{\vec{e}}
\def \pv{\vec{p}}
\def \nv{\vec{n}}
\def \sv{\vec{s}}
\def \tv{\vec{t}}
\def \uv{\vec{u}}
\def \vv{\vec{v}}
\def \wv{\vec{w}}
\def \xv{\vec{x}}
\def \yv{\vec{y}}
\def \zv{\vec{z}}

\def \R{\mathbb{R}}     \def \matC{\vec{C}}
\def \matI{\vec{I}}
\def \matJ{\vec{J}}
\def \matP{\vec{P}}
\def \matU{\vec{U}}

\makeatletter
\newcommand{\nosemic}{\renewcommand{\@endalgocfline}{\relax}}\newcommand{\dosemic}{\renewcommand{\@endalgocfline}{\algocf@endline}}\newcommand{\pushline}{\Indp}\newcommand{\popline}{\Indm\dosemic}\let\oldnl\nl \newcommand{\nonl}{\renewcommand{\nl}{\let\nl\oldnl}}\makeatother

\journal{Computer Methods in Applied Mechanics and Engineering}
\begin{document}

\begin{frontmatter}

  \title{Nonlinear Field-split Preconditioners for Solving Monolithic Phase-field Models of Brittle Fracture}

  \author{Alena Kopani\v{c}\'{a}kov\'{a}\corref{corr_author}}
  \cortext[corr_author]{Corresponding author}
  \ead{alena.kopanicakova@usi.ch}
  \author{Hardik Kothari}
  \ead{hardik.kothari@usi.ch}
  \author{Rolf Krause}
  \ead{rolf.krause@usi.ch}
  \address{Euler Institute, Universit\`{a} della Svizzera italiana, Lugano, Switzerland}

  \begin{abstract}
    One of the state-of-the-art strategies for predicting crack propagation, nucleation, and interaction is the phase-field approach.
Despite its reliability and robustness, the phase-field approach suffers from burdensome computational cost, caused by the non-convexity of the underlying energy functional and a large number of unknowns required to resolve the damage gradients.
In this work, we propose to solve such nonlinear systems in a monolithic manner using the Schwarz preconditioned inexact Newton's (SPIN) method.
The proposed SPIN method leverages the field split approach and minimizes the energy functional separately with respect to displacement and the phase-field, in an additive and multiplicative manner.
In contrast to the standard alternate minimization, the result of this decoupled minimization process is used to construct a preconditioner for a coupled linear system, arising at each Newton's iteration.
The overall performance and the convergence properties of the proposed additive and multiplicative SPIN methods are investigated by means of several numerical examples.
A comparison with widely-used alternate minimization is also performed showing a significant reduction in terms of execution time.
Moreover, we also demonstrate that this reduction grows even further with increasing problem size.

   \end{abstract}

\end{frontmatter}

%\linenumbers
\begin{sloppypar}
  \section{Introduction}
Computational modeling of complex fracture processes, such as crack propagation, branching, or merging, plays an important role in the field of computational mechanics and engineering.
The phase-field approaches to fracture allow to model such complex processes without explicitly representing discontinuities, and have therefore become very popular.
The main idea behind the phase-field approach is to introduce a smooth indicator function that characterizes the state of the material.
This phase-field function is also used to regularize sharp crack interfaces, giving rise to a volumetric approximation of the crack zones.
The regularization is typically based on Ambrosio--Tortorelli elliptic functional~\cite{ambrosio1990approximation}, which models the transition between the broken and intact parts of the domain in a diffused manner.
The size of the transition zone is controlled by a localization limiter, called the length-scale parameter.
Using this volumetric approximation of the sharp crack interfaces, the numerical modeling of the phase-field fractures can be carried out using the standard finite element method (FEM), without any need for computationally tedious tracking/remeshing algorithms.

In this work, we consider a phase-field fracture model based on a variational approach, which has been originally proposed in~\cite{francfort1998revisiting} as an extension to classical Griffith's energy-based principle~\cite{griffith1921phenomena}.
Using such a variational approach, the fracture processes are modeled by minimizing the total potential energy of an underlying system.
The first numerical implementation of the variational phase-field model has been presented by Bourdin et~al.~\cite{bourdin2000numerical}.
Subsequently, a thermodynamically consistent phase-field fracture formulation has been developed by Miehe et~al.~\cite{miehe2010thermodynamically, miehe2010phase}.
Since then, the phase-field fracture models have been extensively used and extended in various directions, for instance dynamics~\cite{Kuhn_dynamic, bourdin2011time}, large-deformations~\cite{del2007variational, hesch2017framework}, or coupled multi-physics applications~\cite{bilgen2017_pneumatic, wick2020multiphysics, heider2021review}.
For a detailed overview of the phase-field fracture models, we refer the reader to~\cite{wu2020phase, zhuang2022phase}.

The biggest computational challenge associated with phase-field fracture modeling is related to its burdensome computational cost, which is caused by three main factors.

Firstly, meshes with high resolution are required to resolve the fracture zones, which gives rise to nonlinear systems with a large number of unknowns.
Adaptive mesh refinement techniques~\cite{patil2018adaptive, ferro2018anisotropic, heister2015primal, mang2020mesh} and multiscale finite element methods~\cite{nguyen2019multiscale} can be employed to alleviate the computational cost by reducing the number of unknowns.
However, the mesh resolution and the number of unknowns increase inevitably with the number of cracks in the structure.
Therefore, in order to obtain an accurate solution of complex fracture problems with feasible computational cost, it is crucial to employ efficient and scalable solution strategies.

Secondly, the non-convexity of the total potential energy, resulting from the coupling between the displacement and the phase-field field, hinders the convergence of standard Newton's method.
Several remedies have been proposed to overcome this difficulty; for example, line-search and trust-region methods~\cite{GERASIMOV2016276, kopanivcakova2020recursive, gerasimov2022second}, arc-length control and path-following approaches~\cite{may2016new, singh2016fracture}, primal-dual active set method~\cite{heister2015primal, noii2021quasi}, augmented Lagrangian approach~\cite{wick2017error}, modified Newton's method~\cite{wick2017modified}, quasi-Newton's methods~\cite{wu2020bfgs, kristensen2020phase}, and Newton's methods with inertia correction~\cite{wambacq2021interior, lampron2021efficient}.
Although all of these methods for solving the arising coupled phase-field fracture problems are promising, they have not been widely adopted in the literature.
This is most likely due to their complexity or lack of robustness.
As a consequence, the majority of the phase-field fracture simulators utilize a staggered solution scheme, also called the alternate minimization (AM) method~\cite{Bourdin2007, miehe2010phase}.
The idea behind the AM method is to minimize the energy functional separately with respect to the displacement and the phase-field, while the other respective variable is kept fixed.
This gives rise to two convex minimization problems, which can be solved efficiently using standard solution strategies.
Despite its robustness, the AM method exhibits slow convergence, which even deteriorates with increasing problem size~\cite{bilgen2019}.
Several approaches to accelerate the AM method have been proposed recently in the literature, e.g., over-relaxation strategies~\cite{farrell2017linear, storvik2021accelerated}, stabilization techniques~\cite{brun2020iterative}, or sub-stepping algorithms~\cite{gupta2020auto}.
However, the applicability of the AM method to large-scale problems remains limited.

Thirdly, we point out that the linear and nonlinear systems arising after FEM discretization are often severely ill-conditioned.
In particular, for second-order partial differential equations (PDEs), the condition number grows as $\pazocal{O}(h^{-2})$, where $h$ denotes the mesh size.
Moreover, in the case of phase-field fracture problems, the ill-conditioning is further aggravated due to the strongly varying stiffness in the structure.
Consequently, the convergence rate of standard iterative solvers deteriorates with the increasing mesh resolution.
As a remedy, preconditioning strategies~\cite{saad2003iterative}, such as multilevel or domain decomposition methods, can be employed.
In the case of phase-field fracture simulations, the efficiency of various linear and nonlinear multilevel methods has been demonstrated by several authors, see for example~\cite{kopanivcakova2020recursive,farrell2017linear,bilgen2017phase,jodlbauer2020matrix,graser2020truncated, zulian2021large, badri2021preconditioning}.
In the context of domain decomposition methods, standard linear (block) Jacobi and additive Schwarz methods have been investigated in~\cite{bilgen2017phase, badri2021preconditioning}.
A variant of the linear additive Schwarz method has been also developed in~\cite{svolos2020updating}, where authors have decomposed the domain into two subdomains, related to localized and healthy parts of the domain.
To the best of our knowledge, there have been no attempts to apply nonlinear domain decomposition methods to solve phase-field fracture problems so far.

Motivated by the observations above, the goal of this work is to propose a solution strategy, which can tackle the non-convexity of the coupled energy functional in a robust manner, and which acts as a scalable preconditioner.
In particular, we develop a variant of the Schwarz Preconditioned Inexact Newton (SPIN) method~\cite{cai2002nonlinearly, liu2015field, liu2018note, gross2021globalization}, specifically tailored to solve phase-field fracture simulations.
The proposed solution strategy is designed to solve the arising nonlinear problems in a monolithic manner by taking advantage of the underlying structure of the coupled problem.
More precisely, we construct a nonlinear preconditioner by partitioning the degrees of freedom into two sets, related to the displacement and the phase-field.
Similar to the AM method, the energy functional is minimized separately with respect to the displacement and the phase-field.
In contrast to the AM method, the result of this minimization process is used to construct the preconditioner for a coupled linear system, arising at each Newton iteration.
Using several representative benchmark problems, we show that the newly developed additive and multiplicative SPIN methods are robust, significantly more efficient than the standard AM method and that their performance does not deteriorate with increasingly larger loading steps or increasing mesh resolution.

This paper is organized as follows.
In \cref{sec:pf}, we review the variational phase-field fracture model, including the finite element discretization of the arising Euler-Lagrange equations.
In \cref{sec:am}, we briefly recall the solution strategies, which are typically used for solving phase-field fracture problems.
\Cref{sec:spin} introduces additive and multiplicative variants of the SPIN method, specifically tailored for the phase-field fracture simulations.
The benchmark problems used to test the proposed SPIN methods are introduced in \cref{sec:benchmark_problems}.
The convergence and performance studies of the SPIN methods, including a comparison with the standard AM method, are presented in \cref{sec:conv_study}.
Eventually, the presented work and possible future extensions are summarized in \cref{sec:conclusion}.

   \section{The phase-field fracture model of brittle fracture}
\label{sec:pf}
In this section, we briefly introduce the phase-field fracture model for brittle fracture.
To this aim, we consider a linear elastic body ${\Omega \subset \R^d}$, ${d \in \{2,3\}}$, with boundary ${\partial \Omega}$.
The boundary ${\partial \Omega}$ is further decomposed into Dirichlet boundary ${\partial \Omega_D}$ and Neumann boundary ${\partial \Omega_N = \partial\Omega \setminus \partial \Omega_D}$, on which Dirichlet and Neumann boundary conditions are applied, respectively.
Due to the application of these boundary conditions or the external body forces, the body $\Omega$ undergoes deformation ${\uv:\Omega\to \R^d}$, which in turn causes crack propagation.
In this work, we consider a quasi-static loading process, where the pseudo time ${t= 1,\ldots, T}$ is used to index the deformation state.
Thus, the Dirichlet and Neumann boundary conditions are prescribed at time $t$ as ${\uv^t_D:\partial\Omega_D\to \R^d}$ and \({\tv_N^t:\partial\Omega_N \to \R^d}\), respectively.

Through this work, we consider linear elastic materials, constitutive law of which is provided by Hooke's law.
The Cauchy stress tensor is then given as
\({
    \boldsymbol{\sigma}(\uv) := 2 \mu \boldsymbol{\varepsilon}(\uv) + \lambda \text{tr}(\boldsymbol{\varepsilon}(\uv))\boldsymbol{I},
  }\)
where $\lambda$ and $\mu$ are the Lam\'e parameters, $\text{tr}(\cdot)$ denotes the trace operator, and $\boldsymbol{I} \in \R^{d\times d}$ stands for the identity tensor.
The second-order strain tensor is given as ${\boldsymbol{\varepsilon}(\uv):= \frac{1}{2}(\nabla \uv +(\nabla \uv)^\top)}$.

\begin{figure}
  \centering
  \begin{subfigure}[b]{0.45\textwidth}
    \centering
    \includegraphics{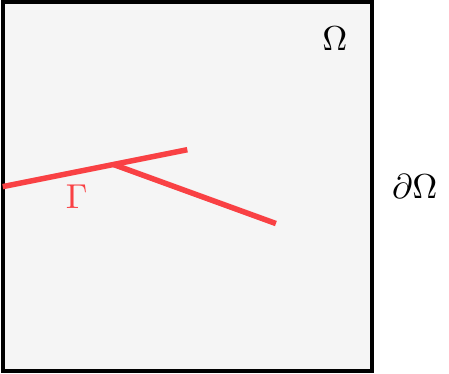}

\end{subfigure}\hspace{1.5cm}
  \begin{subfigure}[b]{0.45\textwidth}
    \centering
    \includegraphics{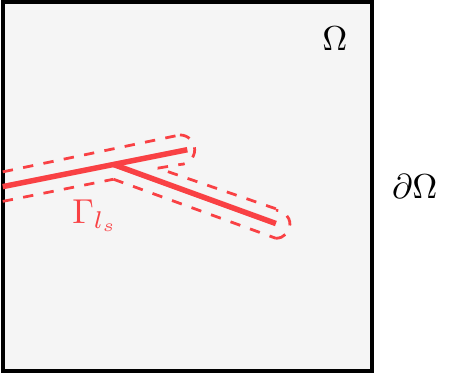}

\label{fig:diffused_crack}
  \end{subfigure}
  \caption{ Left: Computational domain $\Omega$ with the sharp crack $\Gamma$.
    Right: Phase-field regularization of the same crack, denoted by $\Gamma_{l_s}$. The size of the regularized zone is controlled by the length-scale parameter $l_s$.
  }
  \label{fig:notation}
\end{figure}

\subsection{Variational formulation for brittle fracture}
The variational approach to fracture~\cite{francfort1998revisiting} allows us to find the displacement~$\uv^t$ and fracture surface $\Gamma^t$ at time $t$ by minimizing the total potential energy of the body $\Omega$.
Thus, at each pseudo-time-step, the pair $(\uv^t,\Gamma^t)$ can be obtained by solving the following minimization problem:
\begin{linenomath*}
  \begin{equation}
    (\uv^t,\Gamma^t) = \underset{\uv=\uv_D^t \text{ on } \partial\Omega_D \atop \Gamma \supset \Gamma^{t-1}}{\text{arg}\,\min} \pazocal{E}(\uv,\Gamma) := \int_{\Omega \setminus \Gamma} \psi_e(\uv) \ d\Omega + \int_{\Gamma} \pazocal{G}_c \ d\Gamma,
    \label{eq:energy_not_reg}
  \end{equation}
\end{linenomath*}
where the first term denotes the strain energy functional, given as $\psi_e(\uv) = \frac{1}{2}\boldsymbol{\sigma}(\uv):\boldsymbol{\varepsilon}(\uv)$.
The second term describes the fracture energy, defined as the critical energy release $\pazocal{G}_c$ integrated over the sharp fracture surface $\Gamma$.

\begin{remark}
  For simplicity, the energy formulation~\eqref{eq:energy_not_reg} neglects external forces.
\end{remark}

The direct minimization of the energy functional~\eqref{eq:energy_not_reg} is computationally challenging as it requires explicit tracking of the crack paths and the cumbersome remeshing of the domain becomes essential.
The phase-field approach overcomes this difficulty by replacing the fracture energy in~\eqref{eq:energy_not_reg} with its volumetric approximation~\cite{Bourdin2007}.
To this aim, we introduce a phase-field variable ${c:\Omega \to [0,1]}$, which represents the state of the domain, where $c=0$ denotes the intact state, $c=1$ identifies fully broken state and $c\in(0,1)$ denotes a smooth transition between the both states.
Using this phase-field variable~$c$, the fracture is modeled in a diffused manner and the fracture energy from \eqref{eq:energy_not_reg} can be approximated as
\begin{linenomath*}
  \begin{equation}
    \int_\Gamma \pazocal{G}_c \ d \Gamma \approx \frac{\pazocal{G}_c}{c_\omega} \int_\Omega \bigg(\frac{\omega(c)}{l_s} + l_s \left| \nabla c \right|^2 \bigg) d\Omega.
    \label{eq:pf_approx}
  \end{equation}
\end{linenomath*}
Here, the symbol $l_s$ denotes the length-scale parameter that controls the width of the transition zone between the broken and unbroken states.
The function $\omega(\cdot)$ defines the decaying profile of the phase-field and $c_\omega:=4 \int_{0}^1 \sqrt{\omega(c)} \ dc$ denotes an induced normalization constant.
Using \eqref{eq:pf_approx}, we can now reformulate the energy functional \eqref{eq:energy_not_reg} as
\begin{linenomath*}
  \begin{equation}
    \pazocal{E}(\uv, c) := \int_{\Omega} g(c) \ \psi_e(\uv)\ d\Omega + \frac{\pazocal{G}_c}{c_\omega} \int_\Omega \bigg( \frac{\omega(c)}{l_s} + l_s \left| \nabla c \right|^2 \bigg) d\Omega,
    \label{eq:regularized}
  \end{equation}
\end{linenomath*}
where $g(\cdot)$ is a degradation function, which accounts for the loss of stiffness in the fracture.
Particular choices of the functions $g(\cdot)$ and $\omega(\cdot)$ give rise to different phase-field formulations~\cite{kuhn2015degradation,sargado2018high}.
In this work, we adopt the widely used \textit{AT-2} model with ${g(c):=(1-c)^2}$, ${\omega(c):=c^2}$ and ${c_\omega = 2}$, see~\cite{ambrosio1990approximation, bourdin2000numerical,miehe2010thermodynamically} for details.
Given these particular choices of $g(\cdot)$, $\omega(\cdot)$ and $c_\omega$, the minimizer of~\eqref{eq:regularized} asymptotically converges to the minimizer of \eqref{eq:energy_not_reg} as $l_s \to 0$, see $\Gamma$-convergence results in~\cite{giacomini2005ambrosio}.

The main limitation of the phase-field fracture model~\eqref{eq:regularized} is that it allows for crack interpenetration under compressive loading and that the nucleation threshold is symmetric in tension and compression~\cite{de2021nucleation}.
To avoid this unphysical behavior, several modifications of the model~\eqref{eq:regularized} have been proposed.
The most common approach is to decompose the strain energy function $\psi_e$ into a tensile and a compressive part and to allow only the tensile part to be influenced by the damage.
Several types of decomposition have been introduced in the literature, for instance, volumetric-deviatoric split~\cite{amor2009regularized}, or the no-tension model~\cite{freddi2010regularized}.
In this work, we employ the spectral decomposition, originally proposed in~\cite{miehe2010thermodynamically}.
The strain energy functional is decomposed as \( \psi_e(\uv):=\psi^{+}_e(\uv) + \psi^{-}_e(\uv) \), where
\({
    \psi^\pm_e(\uv):= \frac{1}{2} \lambda [\text{tr}^{\pm}(\boldsymbol{\varepsilon})]^2 + \mu \boldsymbol{\varepsilon}^{\pm} \cdot \boldsymbol{\varepsilon}^{\pm}.
  }\)
Here, the positive and negative parts of the strain tensor are obtained as ${\boldsymbol{\varepsilon}^\pm := \sum_i \langle \varepsilon_i \rangle_\pm \nv_i \otimes \nv_i}$, where
$\{\varepsilon_i\}_i$ and $\{\nv_i\}_i$ denote the eigenvalues and the eigenvectors of $\boldsymbol{\varepsilon}$, respectively.
Moreover, the ramp functions are defined as $\langle x \rangle_+ =\max(x,0)$ and $\langle x \rangle_- = \min(x,0)$.
Employing this additive decomposition of the strain energy, and the particular choices of $g(\cdot)$, $\omega(\cdot)$ and $c_\omega$, we can reformulate~\eqref{eq:regularized} as
\begin{linenomath*}
  \begin{equation}
    \pazocal{E}(\uv, c) :=
    \underbrace{ \int_{\Omega} \bigg( (1-c)^2 \ \psi^+_e(\uv) + \psi^-_e(\uv)\bigg)\ d\Omega }_{=:\pazocal{E}_e} + \underbrace{\frac{\pazocal{G}_c}{2} \int_\Omega \bigg( \frac{c^2}{l_s} + l_s \left| \nabla c \right|^2 \bigg) d\Omega}_{=:\pazocal{E}_f}.
    \label{eq:split_regularized}
  \end{equation}
\end{linenomath*}
Please note that the degradation function $g(c)=(1-c)^2$ in~\eqref{eq:split_regularized} affects only the positive part of the elastic energy.

\subsection{Minimization problem}
\label{sec:pf_min_problem}
Using the phase-field fracture energy functional~\eqref{eq:split_regularized}, the quasi-static brittle fracture problem is expressed as
\begin{linenomath*}
  \begin{equation}
    \begin{aligned}
      (\uv^t,c^t) &  & = \underset{\uv=\uv_D^t \text{ on } \partial\Omega_D \atop \partial_t c \geqslant 0}{\text{arg}\,\min} \pazocal{E}(\uv,c).
    \end{aligned}
    \label{eq:energy_reg_constrained}
  \end{equation}
\end{linenomath*}
where the constraint $ \partial_t c \geqslant 0$ imposes the positive evolution of the phase-field, which in turn ensures the fracture irreversibility.
From the numerical point of view, the irreversibility condition can be enforced directly by means of inequality constraints.
Consequently, the constrained optimization algorithms, e.g.,~the active-set methods~\cite{heister2015primal}, or the interior-point methods~\cite{wambacq2021interior} need to be employed.
Alternatively, the irreversibility condition can be modeled indirectly, for example by introducing a history field~\cite{miehe2010phase}, using penalization approaches~\cite{WHEELER201469, gerasimov2019penalization}, or by considering relaxed equality condition~\cite{bourdin2000numerical}.

In this work, we opt for an indirect approach and impose the irreversibility condition using the quadratic penalty method~\cite{gerasimov2019penalization}.
This allows us to reformulate the constrained minimization problem~\eqref{eq:energy_reg_constrained} as the following unconstrained minimization problem:
\begin{linenomath*}
  \begin{equation}
    (\uv^t,c^t) = \underset{\uv=\uv_D^t \text{ on } \partial\Omega_D}{\text{arg}\,\min} \Psi(\uv,c) := \pazocal{E}(\uv,c) + \frac{\gamma}{2} \int_{\Omega} (\langle c-c^{t-1} \rangle_-)^2\ d\Omega,
    \label{eq:energy_reg}
  \end{equation}
\end{linenomath*}
where $\gamma \in \R^+$ denotes the penalty parameter.
The solution of~\eqref{eq:energy_reg} converges to the solution of~\eqref{eq:energy_reg_constrained} as~${\gamma \rightarrow \infty}$.
However, very large values of the penalty parameter~$\gamma$ give rise to ill-conditioning.
Therefore, it is beneficial to employ the smallest possible~$\gamma$, which ensures that the irreversibility condition is satisfied with sufficient accuracy.
Following~\cite{gerasimov2019penalization}, we choose~$\gamma$ as
\begin{linenomath*}
  \begin{equation}
    \label{eq:penalty_parameter}
    \gamma = \frac{\pazocal{G}_c}{l_s} \bigg( \frac{1}{\tau_{\text{irr}}^2} - 1 \bigg),
  \end{equation}
\end{linenomath*}
where $\tau_{\text{irr}}$ denotes a tolerance, with which the irreversibility condition should be satisfied.
This explicit formulation for the parameter $\gamma$ has been derived by exploiting the notion of the optimal phase-field profile as well as the $\Gamma$-convergence in~\cite{gerasimov2019penalization}.

\subsection{Euler-Lagrange equations and finite-element discretization}
We define function spaces ${\bpazocal{V}^t = \{\vv\in [H^1(\Omega)]^d : \vv\vert_{\partial\Omega_D} = \uv_D^t \}}$ and ${\bpazocal{V}_0 = \{\vv\in [H^1(\Omega)]^d : \vv\vert_{\partial\Omega_D} = \boldsymbol{0} \}} $, where $H^1$ denotes the Sobolev space that consists of all square-integrable functions and have square-integrable first weak derivatives.
At each time $t$, the first order necessary optimality conditions of the minimization problem~\eqref{eq:energy_reg} are obtained by differentiating the energy functional~$\Psi$, which gives rise to the following Euler-Lagrange equations:
Find a pair $(\uv, c) \in \bpazocal{V}^t \times H^1(\Omega)$, such that
\begin{linenomath*}
  \begin{equation}
    \begin{aligned}
      \nabla_{\uv} \Psi(\uv, c;\vv) =0, \quad \forall \vv \in \bpazocal{V}_0, \hspace{1.0cm} \text{and} \hspace{1.0cm}
      \nabla_{c} \Psi(\uv, c; w) = 0,\quad \forall w \in H^1(\Omega),
    \end{aligned}
    \label{eq:weak_form}
  \end{equation}
\end{linenomath*}
where
\begin{linenomath*}
  \begin{equation}
    \begin{aligned}
      \nabla_{\uv} \Psi(\uv, c;\vv) & = \big((1-c)^2 \boldsymbol{\sigma}^+(\uv),\boldsymbol{\varepsilon}(\vv)\big)_{\Omega} + \big(\boldsymbol{\sigma}^-(\uv),\boldsymbol{\varepsilon}(\vv)\big)_{\Omega},                  \\
      \nabla_{c} \Psi(\uv, c; w)    & = \big(2(c-1)\psi_e^+(\uv), w\big)_\Omega + \frac{\pazocal{G}_c}{l_s}[(c, w)_\Omega + l_s^2 \big(\nabla c, \nabla w\big)_\Omega ] + {\gamma}(\langle c-c^{t-1} \rangle_- , w)_\Omega.
    \end{aligned}
    \label{eq:weak_form_full}
  \end{equation}
\end{linenomath*}
Here, the symbol $(\cdot,\cdot)_\Omega$ denotes $L^2$ inner product over the domain $\Omega$.

We discretize \eqref{eq:weak_form} using a Galerkin finite element (FE) method.
To this aim, we consider a shape regular, quasi-uniform, conforming quadrilateral mesh $\pazocal{Q}_h$ of the domain $\Omega$.
Moreover, we define the continuous FE spaces
$\bpazocal{V}^t_h = \{\vv_h \in [H^1(\Omega)]^d: \vv_h \vert_Q \in \pazocal{Q}_h \text{ is bilinear }, \vv_h\vert_{(\partial\pazocal{Q}_h)_D} = \uv^t_D \}$,
and
$\pazocal{W}_h = \{w_h \in H^1(\Omega): w_h\vert_Q \in \pazocal{Q}_h \text{ is bilinear}\}$.
Both FE spaces, $\bpazocal{V}_h$ and $\pazocal{W}_h$, are spanned by nodal basis functions $\{\phi_i\}_{i=1}^n$, which we use to approximate the displacement and the phase-field as
\begin{linenomath*}
  \begin{equation}
    \label{eq:FE_approx}
    \uv \approx \uv_h(\xv) = \sum_{i=1}^{n} \sum_{k=1}^d \ev_k \phi_i(\xv) U_i, \hspace{1.0cm} \text{and} \hspace{1.0cm} c \approx c_h(\xv) = \sum_{i=1}^n \phi_i(\xv) C_i,
  \end{equation}
\end{linenomath*}
where $\{\ev_k\}_{k=1,\ldots,d}$ denotes the Euclidean basis of $\R^d$.
The vectors $\matU\in \R^{dn}$ and $\matC \in \R^n$ contain the nodal coefficients for the displacement and the phase-field, respectively.

Using the finite element approximations $\uv_h$ and $c_h$, the discretized Euler-Lagrange equations~\eqref{eq:weak_form} can be expressed by means of residuals $F_u:\R^{dn}\times\R^n \to \R^{dn}$ and $F_c:\R^{dn}\times \R^n \to \R^n$, components of which are given as follows:
\begin{linenomath*}
  \begin{equation}
    \begin{aligned}
      (F_u(\matU,\matC))_i & = \big((1-c_h)^2 \boldsymbol{\sigma}^+(\uv_h),\boldsymbol{\varepsilon}(\ev_k \phi_i)\big)_{\pazocal{Q}_h} + \big(\boldsymbol{\sigma}^-(\uv_h),\boldsymbol{\varepsilon}(\ev_k \phi_i)\big)_{\pazocal{Q}_h}, \qquad \quad \quad \ \ \text{ for } k = 1,\ldots,d, \\
      (F_c(\matU,\matC))_i & = \big(2(c_h-1)\psi_e^+(\uv_h), \phi_i\big)_{\pazocal{Q}_h} + \frac{\pazocal{G}_c}{l_s}[(c_h, \phi_i)_{\pazocal{Q}_h} + l_s^2 \big(\nabla c_h, \nabla \phi_i\big)_{\pazocal{Q}_h} ] + {\gamma}(\langle c_h-c_h^{t-1} \rangle_- , \phi_i)_{\pazocal{Q}_h}.
    \end{aligned}
    \label{eq:residual_components}
  \end{equation}
\end{linenomath*}
Here, the subscript $i$ stands for the nodal index.
Please note, that the component $(F_u(\matU,\matC))_i $ is a vector of dimension $d$.
Using~$F_u$ and $F_c$, the residual of the coupled problem $F:\R^{dn}\times\R^n \to \R^{(d+1)n}$ is defined as~$F(\matU,\matC) := \begin{bmatrix} F_u(\matU,\matC) \\ F_c(\matU,\matC) \end{bmatrix}$.

In this work, we solve the arising discrete coupled problems using nonlinearly preconditioned inexact Newton's method, which requires knowledge about the Jacobian ${F' : \R^{dn}\times \R^n \to \R^{(d+1)n \times (d+1)n}}$.
Similarly to the residual $F$, the Jacobian $\matJ = F'(\matU,\matC)$ has also block structure, i.e.,~$\matJ = \begin{bmatrix} \matJ_{uu} & \matJ_{uc} \\ \matJ_{cu} & \matJ_{cc}\end{bmatrix}$, where its components are given as
\begin{linenomath*}
  \begin{equation}
    \begin{aligned}
      (\matJ_{uu})_{i,j} & = \frac{\partial (F_u)_i}{\partial U_j}= \big((1-c_h)^2 \boldsymbol{\sigma}^+(\ev_l \phi_j),\boldsymbol{\varepsilon}(\ev_k \phi_i)\big)_{\pazocal{Q}_h} + \big(\boldsymbol{\sigma}^-(\ev_l \phi_j),\boldsymbol{\varepsilon}(\ev_k \phi_i)\big)_{\pazocal{Q}_h}, \qquad \text{for } k,l = 1,\ldots,d, \\
      (\matJ_{uc})_{i,j} & = \frac{\partial (F_u)_i}{\partial C_j}= \big( 2(c_h-1) \phi_j \boldsymbol{\sigma}^+(\uv_h),\boldsymbol{\varepsilon}(\ev_k \phi_i)\big)_{\pazocal{Q}_h}, \hspace{4.4cm} \text{for } k = 1,\ldots,d,                                                                                                  \\
      (\matJ_{cu})_{i,j} & = \frac{\partial (F_c)_i}{\partial U_j}= \big(2(c_h-1) \boldsymbol{\sigma}^+(\uv_h):\boldsymbol{\varepsilon}(\ev_k \phi_j), \phi_i \big)_{\pazocal{Q}_h}, \hspace{4.15cm} \text{for } k = 1,\ldots,d,                                                                                                \\
      (\matJ_{cc})_{i,j} & = \frac{\partial (F_c)_i}{\partial C_j}= \big(2 \phi_j \psi_e^+(\uv_h), \phi_i \big)_{\pazocal{Q}_h}
      + \frac{\pazocal{G}_c}{ l_s }[ (\phi_j, \phi_i\big)_{\pazocal{Q}_h} + l_s^2 \big(\nabla \phi_j, \nabla \phi_i \big)_{\pazocal{Q}_h}]
      + \gamma H^{-}(c_h-c_h^{t-1}) \big(\phi_j,\phi_i\big)_{\pazocal{Q}_h}.
    \end{aligned}
  \end{equation}
\end{linenomath*}
Here, the symbol $H^-(\cdot)$ denotes the Heaviside function, defined as
\begin{linenomath*}
  \begin{equation}
    H^{-}(x) = \left \{
    \begin{aligned}
       & 1, &  & \text{if}\ x < 0, \\
       & 0, &  & \text{otherwise}.
    \end{aligned} \right.
  \end{equation}
\end{linenomath*}

   \section{Standard solution strategies}
\label{sec:am}
The solution of the energy minimization problem~\eqref{eq:energy_reg} can be found by solving the nonlinear system of equations  associated with the following first-order necessary optimality conditions: Find $\matU^\ast \in \R^{dn}$ and $\matC^\ast\in \R^n$ such that
\begin{linenomath*}
  \begin{equation}
    F(\matU,\matC) = \begin{bmatrix}
      F_u(\matU,\matC) \\
      F_c(\matU,\matC)
    \end{bmatrix} = \boldsymbol{0}.
    \label{eq:algebraic_weak_form}
  \end{equation}
\end{linenomath*}
The nonlinear system~\eqref{eq:algebraic_weak_form} can be solved using a monolithic solution scheme or by employing the alternate minimization (AM) method.
The monolithic solution scheme requires a solution of the coupled problem, which is computationally challenging due to the non-convexity of energy functional induced by the coupling between the displacement and the phase-field.
In contrast, the AM method decouples the nonlinear problem into two subproblems $F_u(\matU,\matC)=0$ and $F_c(\matU,\matC)=0$ associated with the displacement and the phase-field, respectively.
The subproblems are then solved alternatively in order to update the displacement $\matU$ and the phase-field $\matC$, see \Cref{alg:AM}.
The biggest advantage of the AM method, and also the main reason for its popularity, is the fact that both subproblems are convex~\cite{bourdin2000numerical}.
Hence, the standard solution strategies can be readily applied.
Here, we point out that subproblems might not necessarily be convex if a different phase-field fracture model is employed.
For instance, if the anisotropy~\cite{gerasimov2022second}, or nonlinear constitutive laws~\cite{hesch2017framework} are considered.

\begin{algorithm}[t]
  \caption{Alternate Minimization (AM)}
  \label{alg:AM}
  \DontPrintSemicolon
  \SetKwComment{Comment}{$\triangleright$\ }{}
\KwData{$F:\R^{dn}\times\R^n \to \R^{(d+1)n}$, $F_u:\R^{dn}\times\R^n \to \R^{dn}$, $F_c:\R^{dn}\times\R^n\to \R^n$, $\matU^{(0)}\in \R^{dn}$, $\matC^{(0)}\in\R^n$}
  \KwResult{$\matU^{(k)},\matC^{(k)}$}
  \BlankLine
  $k \mapsfrom 0$

  \While{$\|F(\matU^{(k)},\matC^{(k)})\| \geqslant \epsilon_{\mathrm{rel}\_\mathrm{glob}\_{\mathrm{nonl}}}  \|F(\matU^{(0)},\matC^{(0)})\| $}{
    For fixed $\matC^{(k)}$, find $\matU^{(k+1)}$ by solving $F_u(\matU^{(k+1)},\matC^{(k)}) = 0$       \Comment*[r]{use ND/NK/INK method}
    For fixed $\matU^{(k+1)}$, find $\matC^{(k+1)}$ by solving $F_c(\matU^{(k+1)},\matC^{(k+1)}) = 0$   \Comment*[r]{use ND/NK/INK method}
    $k\mapsfrom k+1$
  }
\end{algorithm}

\subsection{Newton's method}
\label{sec:newton_method}
The monolithic solution scheme requires a solution of a coupled nonlinear system of equations at each (pseudo-)time-step, while AM method requires a solution of two nonlinear systems of equations at each iteration.
A standard approach to solve such a system of nonlinear equations is to employ Newton's method, as it is simple and locally quadratically convergent.
However, the convergence of the method relies on a good initial guess.
Global convergence, convergence to a stationary point irrespectively of the choice of the initial guess, can be ensured by employing some globalization strategy, e.g.,~line-search or trust-region method~\cite{nocedal2006numerical}.
In this work, we employ a cubic backtracking line-search algorithm~\cite[Algorithm A6.3.1, pages 325-327]{dennis1996numerical} with strong Wolfe conditions, see~\ref{sec:strong_wolfe_conditions} for more details.

Let us consider a system of nonlinear equations of the following form: $R(\xv)=0$, where ${R:\R^n \mapsto \R^n}$ and ${\xv \in \R^n}$.
An iteration of Newton's method within a line-search framework is given as
\begin{linenomath*}
  \begin{equation}
    \xv^{(k+1)} = \xv^{(k)} + \alpha^{(k)} \pv^{(k)},
  \end{equation}
\end{linenomath*}
where $\pv^{(k)}$ denotes a search-direction and $\alpha^{(k)} \in \R^{+}$ is step size obtained by the backtracking algorithm.
At each iteration, the search-direction $\pv^{(k)}$ is acquired by solving the linear system of equations
\begin{linenomath*}
  \begin{equation}
    R'(\xv^{(k)})\pv^{(k)} = -R(\xv^{(k)}).
    \label{eq:newton_eq}
  \end{equation}
\end{linenomath*}
A solution of this linear system can be found by employing a direct or an iterative method, giving rise to different variants of Newton's method.
In this work, we consider three particular variants of Newton's method, see also \Cref{alg:Newtons}.

\subsubsection{Newton-direct method (ND)}
\label{sec:ND}
In the context of phase-field fracture simulations, direct methods, such as LU or Cholesky factorization, are typically used to solve the linear systems arising at each Newton's iteration~\cite{lampron2021efficient, bourdin2000numerical}.
The popularity of direct solvers results from their remarkable robustness as well as the open-source availability of highly efficient implementations.
However, the computational cost and memory requirements of direct methods grow with increasing problem size.
In particular, modern sparse direct solvers minimize the amount of required floating-point operations as well as storage by exploiting the sparsity pattern of the matrices.
In this way, the factorization of the system with $n$ unknowns can be performed using $\pazocal{O}(n^{3/2})$ flops in two spatial dimensions and $\pazocal{O}(n^{2})$ flops in three spatial dimensions.
Thus, when $n$ becomes too large ($> 500,000$ dofs), the computational cost and memory demands induced by the fill-in become computationally intractable.
Moreover, direct methods can not be easily parallelized, which in turn limits their scalability.

\subsubsection{Newton-Krylov method (NK)}
\label{sec:NK}
An alternative to the ND method is the Newton-Krylov~(NK) method, which employs a Krylov subspace method to solve the linear system of equations arising at each Newton's iteration.
The Krylov subspace methods are the iterative methods, which find the solution of the linear system of equations by projecting a sequence of iterates onto the Krylov subspaces.
Krylov methods have smaller memory requirements compared to direct methods and are well-suited for parallel processing environments, as their main building blocks e.g.,~ matrix-vector multiplication, or inner product, allow for hybrid parallelization.

Based on the problem at hand, a specific Krylov method might be suitable to solve the problem efficiently.
In the context of phase-field fracture problems considered in this work, the arising coupled linear systems can be solved efficiently using the MINRES method.
This is due to the fact that the Jacobian matrices of the coupled nonlinear system are symmetric but not necessarily positive definite due to the non-convex nature of the problem.
For the AM method, the conjugate gradient (CG) method represents an ideal choice, as the Jacobians that arise while solving displacement and the phase-field subproblems are symmetric positive definite, due to the convexity of $F_u$ and $F_c$.

A major drawback of Krylov methods is that their convergence speed deteriorates with increasing problem size and increasing condition number.
As a consequence, preconditioning strategies have to be employed in order to accelerate convergence.
In this work, we precondition the coupled linear systems using a nonlinear additive and multiplicative field-split preconditioner, described in \cref{sec:spin}.
In order to precondition the linear systems, which arise while solving the displacement and the phase-field subproblems, we employ an algebraic multigrid (AMG) method. 

\begin{algorithm}[t]
  \caption{Newton/Inexact Newton (ND/NK/INK)}
  \label{alg:Newtons}
  \DontPrintSemicolon
  \SetKwComment{Comment}{$\triangleright$\ }{}
  \SetKwInOut{Input}{Input}
  \SetKwInOut{Output}{Output}
  \KwData{$R: \R^n \to \R^n, \xv^{(0)} \in \R^n$}
  \KwResult{$\xv^{(k)}$}
  \BlankLine
  $k \mapsfrom 0$ 

  \While {$\|R(\xv^{(k)})\| > \epsilon_{{\mathrm{rel}}\_{\mathrm{sub}}\_{\mathrm{nonl}}} \| R(\xv^{(0)})\| $}{
    Find $\pv^{(k)}$ by solving $R'(\xv^{(k)}) \pv^{(k)} = - R(\xv^{(k)})$ \Comment*[r]{solve linear system}
    \If{$\mathrm{ND}$}{
      solve the system using a direct solver \;
    }
    \If{$\mathrm{NK}$}{
      solve the system using a Krylov method, such that \\
      \hspace{3.15cm} $\| R'(\xv^{(k)}) \pv^{(k)} + R(\xv^{(k)})\| \leqslant \max(\epsilon_{\text{rel}\_\text{lin}} \|R(\xv^{(k)})\|, \epsilon_{\text{abs}\_\text{lin}}) $ \\
}
    \If{$\mathrm{INK}$}{
      solve the system using a Krylov method, such that \\
      \hspace{3.15cm} $\|R'(\xv^{(k)}) \pv^{(k)} + R(\xv^{(k)})\| \leqslant \max(\eta\|R(\xv^{(k)})\|, \epsilon_{\text{abs}\_\text{lin}})$\;
    }
    Find $\alpha^{(k)}$ using a backtracking algorithm   \Comment*[r]{use line-search}
    $\xv^{(k+1)} \mapsfrom \xv^{(k)} + \alpha^{(k)} \pv^{(k)}$ \Comment*[r]{update iterate}
    $k\mapsfrom k+1$
  }
\end{algorithm}

\subsubsection{Inexact Newton-Krylov method (INK)}
\label{sec:INK}
Although NK methods are an improvement upon the ND methods, solving the arising linear systems accurately is computationally expensive.
Moreover, if the current iterate $\xv^{(k)}$ is far from a solution, the local linear model $R(\xv^{(k)}) + R'(\xv^{(k)}) \pv^{(k)}$ might not approximate the functional $R$ well enough.
Since ${R(\xv^{(k)}) + R'(\xv^{(k)}) \pv^{(k)}}$ is also the residual of the linear system~\eqref{eq:newton_eq}, obtaining $\pv^{(k)}$ by solving ${R'(\xv^{(k)}) \pv^{(k)}=-R(\xv^{(k)})}$ exactly might cause over-solving, which in turn leads to little or no progress towards a solution.
As a result, far from the solution, it might be cheaper as well as more effective to obtain less accurate Newton's steps.

Inexact Newton's (IN) methods~\cite{eisenstat1996choosing} incorporate the aforementioned observations  and
determine how accurately the linear systems should be solved at each Newton's iteration.
More precisely, the search-direction $\pv^{(k)}$ is obtained by solving~\eqref{eq:newton_eq}, such that
\begin{linenomath*}
  \begin{equation}
    \|R'(\xv^{(k)}) \pv^{(k)} + R(\xv^{(k)})\| \leqslant \eta\|R(\xv^{(k)})\|.
    \label{eq:inex_Newton}
  \end{equation}
\end{linenomath*}
Following~\cite{cai1994newton}, we employ a simple, but effective, choice of the parameter $\eta$, i.e.,~$\eta= 10^{-4}$.
In this way, the induced tolerance $ \eta\|R(\xv^{(k)})\|$, with which~\eqref{eq:newton_eq} is solved, is updated at each iteration by taking into account the information about the current residual $R(\xv^{(k)})$.
This tolerance adapts proportionally to $\|R(\xv^{(k)})\|$, which avoids over-solving in the initial phases of the solution process and ensures the quadratic convergence of IN method close to a solution~\cite{eisenstat1996choosing}.
In this work, we solve the arising linear systems inexactly, using preconditioned Krylov methods, which gives rise to the inexact Newton-Krylov (INK) method.

   \section{Field-split Schwarz Preconditioned Inexact Newton (SPIN) method}
\label{sec:spin}
Standard Newton's method converges quadratically to a solution if a good initial guess is provided.
However, if a good initial guess is not known, Newton's methods might exhibit very slow convergence until a local neighborhood of a solution is approached.
Slow convergence is typically associated with unbalanced and highly localized nonlinearities.
In the context of the phase-field fracture problems considered in this work, the unbalanced nonlinearities occur due to coupling between the displacement and the phase-field, and due to locally varying material stiffness and steep gradients of the phase-field function.

Nonlinear preconditioning strategies, such as SPIN~\cite{cai2002nonlinearly, liu2015field}, RASPEN~\cite{dolean2016nonlinear}, or nonlinear elimination~\cite{lanzkron1996analysis} can enhance the convergence of Newton's method by rebalancing the nonlinearities, or by transforming the basis of the solution space.
These nonlinear preconditioning strategies utilize the decomposition of the solution space into multiple subspaces, for instance, related to different parts of the computational domain or the different fields.

In this work, we adopt the SPIN methodology and instead of solving the original coupled nonlinear system~\eqref{eq:algebraic_weak_form}, we solve an equivalent nonlinearly preconditioned system.
This nonlinearly preconditioned system is constructed such that it has the same root as the original nonlinear system, but its nonlinearities are more balanced.
Motivated by the robustness of the AM method, we construct the nonlinearly preconditioned system of equations using the field-split approach, thus by decomposing the original problem into two subproblems that are associated with the displacement and the phase-field.
From the algorithmic point of view, the SPIN strategy can be performed in an additive or multiplicative manner, which gives rise to two distinct variants of the SPIN method, named ASPIN and MSPIN, respectively.

The additive approach is more suitable for parallel computing environments, as nonlinear subproblems associated with different subspaces can be solved simultaneously.
In contrast, the multiplicative approach is inherently sequential, but it usually exhibits faster convergence than the additive approach.

\subsection{Field-split ASPIN and MSPIN methods}
Let $G: \R^{dn}\times\R^n \to \R^{(d+1)n}$ be a nonlinear preconditioner.
Using the field-split approach, the preconditioner $G$ is defined as follows
\begin{linenomath*}
  \begin{equation}
    G(\matU, \matC) = \begin{bmatrix}
      {G}_u (\matU, \matC) \\
      {G}_c (\matU, \matC)
    \end{bmatrix},
  \end{equation}
\end{linenomath*}
where $ G_u: \R^{dn}\times\R^n \to \R^{dn}$ and $ G_c: \R^{dn}\times\R^n \to \R^{n}$ are solution operators associated with displacement and phase-field, respectively.
Ideally, the functions $G_u$ and $G_c$ approximate the inverse of nonlinear operators $F_u$ and $F_c$, i.e.,~$G_u \approx F_u^{-1}$, and $G_c \approx F_c^{-1}$.
We, therefore, expect that at the fixed point, i.e., at the solution of $F(\matU, \matC)=0$, the following relation holds:
\begin{linenomath*}
  \begin{equation}
  [ \matU,\matC ]^{\top} = G(\matU,\matC).
  \end{equation}
\end{linenomath*}
Using this fixed point iteration, we can construct nonlinearly preconditioned system of equations ${\pazocal{F}: \R^{dn}\times\R^n \to \R^{(d+1)n}}$ as follows
\begin{linenomath*}
  \begin{equation}
    \pazocal{F}(\matU,\matC) := [ \matU,\matC ]^{\top} - G(\matU,\matC) = 0.
    \label{eq:prec_algebraiic_weak_form_new}
  \end{equation}
\end{linenomath*}
The nonlinearly preconditioned system of equations~\eqref{eq:prec_algebraiic_weak_form_new} has the same solution as the original nonlinear system associated with the coupled phase-field fracture problem~\eqref{eq:algebraic_weak_form}.
However, the preconditioned system~\eqref{eq:prec_algebraiic_weak_form_new} should be easier to solve, as the nonlinearities arising due to the coupling between displacement and phase-field should be more balanced.

An exact form of the preconditioner $G$ is not known apriori.
Therefore, we construct the preconditioned residual~$\pazocal{F}$ by solving the nonlinear subproblems associated with the displacement and the phase-field.
In this case, $G_u(\matU, \matC)$ can be obtained as $\matU_{\text{prec}}$, which solves $F_u(\matU_{\text{prec}}, \matC)=0$, for a fixed $\matC$.
Similarly, $G_c(\matU, \matC)$ can be obtained as $\matC_{\text{prec}}$, which solves $F_c(\matU, \matC_{\text{prec}})=0$, for a fixed $\matU$.
In practice, the nonlinear preconditioner $\pazocal{F}$ can be constructed in an additive or a multiplicative manner, which gives rise to two different variants of the SPIN method, namely ASPIN and MSPIN.
In the additive case, the preconditioned nonlinear system~$\pazocal{F}^{\mathrm{add}}$ is defined as follows
\begin{linenomath*}
  \begin{equation}
    \pazocal{F}^\mathrm{add}(\matU,\matC) := \begin{bmatrix}
      \matU - G_u(\matU,\matC) \\
      \matC - G_c(\matU,\matC)
    \end{bmatrix} = \boldsymbol{0}.
    \label{eq:additive}
  \end{equation}
\end{linenomath*}
Here, both subproblems can be solved simultaneously, as their solution depends only on the values of displacement $\matU$ and phase-field $\matC$.
In the multiplicative case, the preconditioned residual~$\pazocal{F}^{\mathrm{mult}}$ has the following form:
\begin{linenomath*}
  \begin{equation}
    \pazocal{F}^\mathrm{mult}(\matU,\matC) := \begin{bmatrix}
      \matU - G_u(\matU,\matC) \\
      \matC - G_c(\matU - G_u(\matU,\matC), \matC)
    \end{bmatrix} = \boldsymbol{0}.
    \label{eq:multiplicative}
  \end{equation}
\end{linenomath*}
Here, the subproblems have to be solved sequentially, as the solution of the phase-field subproblem relies on the solution of the displacement subproblem, namely $G_u(\matU,\matC)$.

Once, the nonlinearly preconditioned system $\pazocal{F}$ is formed, we can find its root using the inexact-Newton method.
Thus, on each $k$-th iteration, the search direction $\pv^{(k)}$ is found by solving the following linear system of equations:
\begin{linenomath*}
  \begin{equation}
    \pazocal{F}'(\matU^{(k)},\matC^{(k)}) \pv^{(k)} = \pazocal{F}(\matU^{(k)},\matC^{(k)}),
    \label{eq:global_linear}
  \end{equation}
\end{linenomath*}
where $\pazocal{F}'$ denotes the Jacobian of $\pazocal{F}$.
However, it is not straightforward to construct the Jacobian $\pazocal{F}'$ explicitly, or to evaluate its action on a vector, as the preconditioned residual $\pazocal{F}$ is defined implicitly.
Following~\cite{liu2015field}, we can approximate the Jacobian $\pazocal{F}'$ utilizing the Jacobian~$F'$ of the original residual $F$.
In the additive case, the approximation of the Jacobian $(\pazocal{F}^{\mathrm{add}})'$ evaluated at iterate $k$ is given as
\begin{linenomath*}
  \begin{equation}
    (\pazocal{F}^{\mathrm{add}})'(\matU^{(k)},\matC^{(k)}) \approx \underbrace{\begin{bmatrix}
      \matJ_{uu}^{(k)} &                  \\
                       & \matJ_{cc}^{(k)}
    \end{bmatrix}^{-1}}_{=: \matP^{(k)}_{\mathrm{add}}}
    \underbrace{\begin{bmatrix}
        \matJ_{uu}^{(k)} & \matJ_{uc}^{(k)} \\
        \matJ_{cu}^{(k)} & \matJ_{cc}^{(k)}
      \end{bmatrix}}_{= \matJ^{(k)}},
    \label{eq:Jac_add}
  \end{equation}
\end{linenomath*}
where $\matP_{\mathrm{add}}^{(k)}$ denotes the additive preconditioner and $\matJ^{(k)} = F'(\matU^{(k)},\matC^{(k)})$.
In the multiplicative case, the approximation of the Jacobian $(\pazocal{F}^{\mathrm{mult}})'$ evaluated at iterate $k$ is given as
\begin{linenomath*}
  \begin{equation}
    (\pazocal{F}^{\mathrm{mult}})'(\matU^{(k)},\matC^{(k)}) \approx  \underbrace{\begin{bmatrix}
      \matJ_{uu}^{(k)} &                  \\
      \matJ_{cu}^{(k)} & \matJ_{cc}^{(k)}
    \end{bmatrix}^{-1}}_{=: \matP^{(k)}_{\mathrm{mult}}}
    \underbrace{\begin{bmatrix}
        \matJ_{uu}^{(k)} & \matJ_{uc}^{(k)} \\
        \matJ_{cu}^{(k)} & \matJ_{cc}^{(k)}
      \end{bmatrix}}_{= \matJ^{(k)}},
    \label{eq:Jac_mult}
  \end{equation}
\end{linenomath*}
where $\matP_{\mathrm{mult}}^{(k)}$ denotes the multiplicative preconditioner.

Employing~\eqref{eq:Jac_add} or \eqref{eq:Jac_mult},  the linear system of equation~\eqref{eq:global_linear} can be reformulated as follows:
\begin{linenomath*}
  \begin{equation}
    \matP^{(k)} \matJ^{(k)} \pv^{(k)} = \pazocal{F}(\matU^{(k)},\matC^{(k)}),
    \label{eq:global_linear_new}
  \end{equation}
\end{linenomath*}
where $\matP^{(k)}$ denotes additive or multiplicative preconditioner matrix.
In view of using the inexact-Newton method, the linear system~\eqref{eq:global_linear_new} is solved inexactly as discussed in \cref{sec:am}.

\begin{algorithm}[t]
  \caption{ASPIN / MSPIN}
  \label{alg:spin}
  \DontPrintSemicolon
  \SetKwComment{Comment}{$\triangleright$\ }{}
  \SetKwInOut{Input}{Input}
  \SetKwInOut{Output}{Output}

  \KwData{$F:\R^{dn}\times\R^n \to \R^{(d+1)n}$, $F_u:\R^{dn}\times\R^n \to \R^{dn}$, $F_c:\R^{dn}\times\R^n\to \R^n$, $\matU^{(0)}\in \R^{dn}$, $\matC^{(0)}\in\R^n$}
  \KwResult{$\matU^{(k)}, \matC^{(k)}$}

  \BlankLine
  $k \mapsfrom 0$

  \While {$\|F(\matU^{(k)},\matC^{(k)})\| \geqslant \epsilon_{{\mathrm{rel}}\_{\mathrm{glob}}\_{\mathrm{nonl}}}  \|F(\matU^{(0)},\matC^{(0)})\| $}   {
\tikzmk{A}
    \If{additive SPIN}{
      For fixed $\matC^{(k)}$ find $\matU_{\mathrm{prec}}$ by solving $F_u(\matU_{\mathrm{prec}},\matC^{(k)}) = 0$ \Comment*[r]{use INK method}
      For fixed $\matU^{(k)}$ find $\matC_{\mathrm{prec}}$ by solving $F_c(\matU^{(k)},\matC_{\mathrm{prec}}) = 0$ \Comment*[r]{use INK method}
      Form $\matP_{\mathrm{add}}^{(k)}$ {by means of} \eqref{eq:Padd} \Comment*[r]{construct ASPIN preconditioner}
    }
\If{multiplicative SPIN}{
      For fixed $\matC^{(k)}$ find $\matU_{\mathrm{prec}}$ by solving $F_u(\matU_{\mathrm{prec}},\matC^{(k)}) = 0$ \Comment*[r]{use INK method}
      For fixed $\matU_{\mathrm{prec}}$ find $\matC_{\mathrm{prec}}$ by solving $F_c(\matU_{\mathrm{prec}},\matC_{\mathrm{prec}}) = 0$ \Comment*[r]{use INK method}
      Form $\matP_{\mathrm{mult}}^{(k)}$ {by means of} \eqref{eq:Pmult} \Comment*[r]{construct MSPIN preconditioner}
    }
    $\sv^{(k)} \mapsfrom [ \matU_{\mathrm{prec}}-\matU^{(k)}, \matC_{\mathrm{prec}}-\matC^{(k)}]^\top$  \Comment*[r]{evaluate nonlinearly preconditioned residual}
    \tikzmk{B}\boxit{brownmy!70!yellow}

\tikzmk{A}
    $\matJ^{(k)} \mapsfrom F'(\matU^{(k)},\matC^{(k)})$ \Comment*[r]{assemble monolithic Jacobian}
    Find $\pv^{(k)}$ by approximately solving $ \matP_{\mathrm{add/mult}}^{(k)}  \matJ^{(k)} \pv^{(k)} = \sv^{(k)} $  \Comment*[r]{use a Krylov method}
    \pushline \nonl
    \hspace{4.4cm} such that  $\left\lVert \matP_{\mathrm{add/mult}}^{(k)}  \matJ^{(k)} \pv^{(k)} - \sv^{(k)} \right\rVert \leqslant \eta \left\lVert \sv^{(k)} \right\rVert$
\popline
    \DontPrintSemicolon

    Find $\alpha^{(k)}$ using a backtracking algorithm \Comment*[r]{use line-search}
    $[\matU^{(k+1)},\matC^{(k+1)}]^\top \mapsfrom [ \matU^{(k)},\matC^{(k)}  ]^\top+ \alpha^{(k)} [\pv_u^{(k)},\pv_c^{(k)} ]^\top$ \Comment*[r]{update the iterates}
    \tikzmk{B}\boxit{gray!70!}
    $k\mapsfrom k+1$
  }
\end{algorithm}

\subsection{Algorithmic description of SPIN method}
\label{sec:alg_descr_SPIN}
In this section, we briefly discuss the algorithmic components and the practical implementation of the additive and the multiplicative SPIN methods.
In general, the SPIN methods, summarized in \Cref{alg:spin}, consist of two steps.
The first step (highlighted by yellow color) is associated with the construction of the additively/multiplicatively preconditioned nonlinear residual,  while the second step (highlighted by gray color) is associated with the solution of a global nonlinearly preconditioned linear system of equations.

\subsubsection{Subproblem step - construction of the preconditioned residual}
On $k$-th iteration of the SPIN method, the nonlinearly preconditioned residual $\pazocal{F}$ is constructed by solving the subproblems associated with the displacement and the phase-field using the inexact Newton's method (\Cref{alg:Newtons}).
Thus, we seek the solution $\matU_{\mathrm{prec}}$ of the displacement subproblems $F_u=0$ and the solution $\matC_{\mathrm{prec}}$ of the phase-field subproblem $F_c=0$.
In case of ASPIN method, the subproblem $F_u = 0$ is solved with a fixed value of the phase-field $\matC^{(k)}$ and the subproblem $F_c = 0$ is solved with a fixed value of displacement $\matU^{(k)}$.
The most appealing feature of the ASPIN method is that both subproblems can be solved simultaneously.
In the case of the MSPIN method, the nonlinear subproblem $F_u = 0$ is solved first, while the phase-field variable $\matC^{(k)}$ is held fixed.
Afterward, the subproblem $F_c = 0$ is solved with the fixed value of the displacement $\matU_{\mathrm{prec}}$, where $\matU_{\mathrm{prec}}$ is
the solution of the displacement subproblem $F_u = 0$, obtained at the previous step.
Here, we highlight the fact that the construction of multiplicatively preconditioned residual coincides with a single iteration of the AM method (\Cref{alg:AM}).

After both subproblems are solved, we can use the obtained solutions $\matU_{\mathrm{prec}}$ and $\matC_{\mathrm{prec}}$ to construct the nonlinearly preconditioned residual.
This is achieved by evaluating the corrections computed during the preconditioning step as
\begin{linenomath*}
  \begin{equation}
    \sv_u = \matU_{\mathrm{prec}} - \matU^{(k)},\qquad \qquad \sv_c = \matC_{\mathrm{prec}} - \matC^{(k)}.
    \label{eq:corrections}
  \end{equation}
\end{linenomath*}
Here, $\sv_u$ and $\sv_c$ denote the algebraic representation of the preconditioned residuals $\pazocal{F}_u(\matU^{(k)}, \matC^{(k)})$ and $\pazocal{F}_c(\matU^{(k)},\matC^{(k)})$, respectively.

\subsubsection{Global step - construction and solution of the preconditioned monolithic linear system}
Once the preconditioned residual is formed, we can obtain a search direction~$\pv^{(k)}$ by solving the linear system of equations $(\matP^{(k)}\matJ^{(k)}) \pv^{(k)} = \sv^{(k)}$.
To form this linear system, we have to assemble the Jacobian $\matJ^{(k)}$ associated with the coupled nonlinear residual $F$ and the preconditioning matrix~$\matP^{(k)}$.
Unfortunately, the explicit assembly of the matrix $\matP^{(k)}$ is computationally and memory exhaustive, as it requires an inverse of the block-diagonal, or lower-triangular part of the Jacobian $\matJ^{(k)}$, recall~\ref{eq:Jac_add} and~\ref{eq:Jac_mult}.
To overcome the difficulty, we solve the arising linear systems using Krylov subspace methods, which do not require an explicit representation of the preconditioning matrix $\matP^{(k)}$.
Instead, it is only necessary to evaluate an application of a matrix $\matP^{(k)}\matJ^{(k)}$ to a vector $\vv$, i.e., $\yv = \matP^{(k)} \matJ^{(k)} \vv$.
This evaluation can be carried out sequentially in two steps.
Firstly, we multiply the Jacobian matrix $\matJ^{(k)}$ with a vector $\vv$, i.e., $\wv = \matJ^{(k)}\vv$.
Secondly, we multiply the preconditioning matrix $\matP^{(k)}$ with the vector $\wv$, i.e., $\yv = \matP^{(k)} \wv$.
In order to perform this second multiplication efficiently, matrices $\matP^{(k)}_{\mathrm{add}}$ and $\matP^{(k)}_{\mathrm{mult}}$ defined in \eqref{eq:Jac_add} and \eqref{eq:Jac_mult} have to be recast to computationally more suitable form.
In particular, we reformulate the additive preconditioner $\matP^{(k)}_{\mathrm{add}}$ as
\begin{linenomath*}
  \begin{equation}
    \matP^{(k)}_{\mathrm{add}} : = \begin{bmatrix}
      \matJ_{uu}^{(k)} &                  \\
                       & \matJ_{cc}^{(k)}
    \end{bmatrix}^{-1} = \begin{bmatrix}
      (\matJ_{uu}^{(k)})^{-1} &                         \\
                              & (\matJ_{cc}^{(k)})^{-1}
    \end{bmatrix},
    \label{eq:Padd}
  \end{equation}
\end{linenomath*}
while the multiplicative preconditioner $\matP^{(k)}_{\mathrm{mult}}$ is reformulated as
\begin{linenomath*}
  \begin{equation}
    \matP^{(k)}_{\mathrm{mult}} : = \begin{bmatrix}
      \matJ_{uu}^{(k)} &                  \\
      \matJ_{cu}^{(k)} & \matJ_{cc}^{(k)}
    \end{bmatrix}^{-1} :=
    \begin{bmatrix}
      \matI_{uu} &                         \\
                 & (\matJ_{cc}^{(k)})^{-1}
    \end{bmatrix}
    \begin{bmatrix}
      \matI_{uu}         &            \\
      - \matJ^{(k)}_{cu} & \matI_{cc}
    \end{bmatrix}
    \begin{bmatrix}
      (\matJ_{uu}^{(k)})^{-1} &            \\
                              & \matI_{cc}
    \end{bmatrix},
    \label{eq:Pmult}
  \end{equation}
\end{linenomath*}
where $\matI_{uu} \in \R^{dn \times dn}$ and $\matI_{cc} \in \R^{n \times n}$ denote identity matrices.

By examining the structure of the matrices $\matP_{\mathrm{add}}^{(k)}$ and $\matP_{\mathrm{mult}}^{(k)}$, it is now clear that the multiplication of the matrix $\matP^{(k)}$ with a vector requires the solution of two linear systems, associated with blocks $\matJ_{uu}^{(k)}$ and $\matJ_{cc}^{(k)}$.
Here, we solve these linear systems using preconditioned Krylov methods.
Since the solution of these linear systems is only used to approximate the matrix-vector product, therefore, it is sufficient to solve these linear systems inexactly.
The process of applying the nonlinearly preconditioned operator $\matP^{(k)} \matJ^{(k)}$ to a vector $\vv$ is summarized in \Cref{alg:action_ASPIN} and \Cref{alg:action_MSPIN}.

From the structure of operators $\matP^{(k)}_{\mathrm{add}}\matJ^{(k)}$ and $\matP^{(k)}_{\mathrm{mult}} \matJ^{(k)}$, we can also see that the additively preconditioned operator $\matP^{(k)}_{\mathrm{add}}\matJ^{(k)}$ is symmetric, while multiplicatively preconditioned operator $\matP^{(k)}_{\mathrm{mult}} \matJ^{(k)}$ is non-symmetric.
As a consequence, a different Krylov method is suitable to solve the additively and multiplicatively preconditioned linear system.
For instance, the generalized minimal residual (GMRES) method can be employed to efficiently solve the non-symmetric linear systems arising at each MSPIN iteration, while the minimal residual (MINRES) method can be used at each ASPIN iteration.
We also point out that the preconditioned Jacobian $\matP^{(k)}\matJ^{(k)}$ has a significantly smaller condition number than the Jacobian of the original system $\matJ^{(k)}$.
This is due to the fact that the SPIN preconditioner implicitly provides robust linear preconditioning by means of the preconditioning matrix $\matP^{(k)}$.
As a consequence, solving~\eqref{eq:global_linear_new} typically requires only very few iterations of a Krylov method.

\begin{algorithm}[t]
  \caption{Evaluation of $\matP^{(k)}_{\mathrm{add}} \matJ^{(k)} \vv$}
  \label{alg:action_ASPIN}
  \DontPrintSemicolon
  \SetKwComment{Comment}{$\triangleright$\ }{}
  \SetKwInOut{Input}{Input}
  \SetKwInOut{Output}{Output}

  \KwData{$\matJ^{(k)} \in \R^{(d+1)n \times (d+1)n}, \vv_u\in\R^{dn},  \vv_c \in \R^n$}
  \KwResult{$\yv_u, \yv_c$}

  \BlankLine
  $ [\wv_u, \wv_c]^\top \mapsfrom \matJ^{(k)}  [\vv_u, \vv_c]^\top$ \Comment*[r]{apply Jacobian to vector}
  Obtain $\matJ^{(k)}_{uu}, \matJ^{(k)}_{cc}$ by decomposing $\matJ^{(k)}$\Comment*[r]{get block components of Jacobian}
  Find $\yv_u$ by approximately solving $\matJ^{(k)}_{uu} \yv_u = \wv_u$  \Comment*[r]{use a Krylov method}
  \pushline \nonl
  \hspace{1.75cm}  such that  $\|\matJ^{(k)}_{uu} \yv_u - \wv_u\| \leqslant \epsilon_{\text{app}\_{\text{lin}}} \| \wv_u \|$
  \popline
  \DontPrintSemicolon

  Find $\yv_c$ by approximately solving $\matJ^{(k)}_{cc} \yv_c = \wv_c$ \Comment*[r]{use a Krylov method}
  \pushline \nonl
  \hspace{1.2cm} such that  $\|\matJ^{(k)}_{cc} \yv_c - \wv_c\| \leqslant \epsilon_{\text{app}\_{\text{lin}}} \| \wv_c \|$
\end{algorithm}

\begin{algorithm}[t]
  \caption{Evaluation of $\matP^{(k)}_{\mathrm{mult}} \matJ^{(k)} \vv$}
  \label{alg:action_MSPIN}
  \DontPrintSemicolon
  \SetKwComment{Comment}{$\triangleright$\ }{}
  \SetKwInOut{Input}{Input}
  \SetKwInOut{Output}{Output}

  \KwData{$\matJ^{(k)} \in \R^{(d+1)n \times (d+1)n}, \vv_u\in\R^{dn},  \vv_c \in \R^n$}
  \KwResult{$\yv_u, \yv_c$}

  \BlankLine
  $ [\wv_u, \wv_c]^\top \mapsfrom \matJ^{(k)}  [\vv_u, \vv_c]^\top$ \Comment*[r]{apply Jacobian to vector}
  Obtain $\matJ^{(k)}_{uu}, \matJ^{(k)}_{cc}, \matJ^{(k)}_{cu}$ by decomposing $\matJ^{(k)}$\Comment*[r]{get block components of Jacobian}
  Find $\yv_u$ by approximately solving $\matJ^{(k)}_{uu} \yv_u = \wv_u$  \Comment*[r]{use a Krylov method}
  \pushline \nonl
  \hspace{1.75cm} such that  $\|\matJ^{(k)}_{uu} \yv_u - \wv_u\| \leqslant \epsilon_{\text{app}\_{\text{lin}}} \| \wv_u \|$
\popline
  \DontPrintSemicolon

  ${\zv_c} \mapsfrom \wv_c - \matJ^{(k)}_{cu} \yv_u$ \Comment*[r]{apply off-diagonal block}
  Find $\yv_c$ by approximately solving $\matJ^{(k)}_{cc} \yv_c = \zv_c$ \Comment*[r]{use a Krylov method}
  \pushline \nonl
  \hspace{1.2cm} such that  $\|\matJ^{(k)}_{cc} \yv_c - \zv_c\| \leqslant \epsilon_{\text{app}\_{\text{lin}}} \| \zv_c \|$
\end{algorithm}

After the search direction~$\pv^{(k)}$ is obtained by solving the nonlinearly preconditioned system of equations~\eqref{eq:global_linear_new}, the line-search method is used to determine the appropriate step size $\alpha^{(k)}$.
This is followed by the update of variables $\matU^{(k)}$ and $\matC^{(k)}$, performed as $\matU^{(k+1)}=\matU^{(k)}+\alpha^{(k)} \pv^{(k)}_u$ and $\matC^{(k+1)}=\matC^{(k)}+\alpha^{(k)} \pv^{(k)}_c$, where $\pv^{(k)}_u$ and $\pv^{(k)}_c$ denote the components of $\pv^{(k)}$ associated with the displacement and the phase-field, respectively.

Compared to the AM method, one iteration of the SPIN algorithm is computationally more expensive.
The additional computational cost is associated with the global step of the SPIN algorithm.
More precisely, we have to assemble monolithic Jacobian $\matJ^{(k)}$, which can be performed efficiently by reusing the block $\matJ_{uu}$ constructed during the first iteration of Newton's method, while solving the subproblem $F_u=0$.
In the case of the ASPIN method, block $\matJ_{cc}$ can be also reused.
In contrast, for the MSPIN method, the block~$\matJ_{cc}$ has to be re-evaluated at $\matU^{(k)}$ and $\matC^{(k)}$.
Moreover, the off-diagonal coupling blocks $\matJ_{uc}$ and $\matJ_{cu}$ have to be assembled for both SPIN methods.
Besides the assembly cost, the preconditioned coupled system of equations has to be solved using the Krylov method, which requires the solution of two linear systems, associated with $\matJ_{uu}$ and $\matJ_{cc}$, at each iteration.
Even though one iteration of the SPIN method is more expensive than an iteration of the AM method, its superior convergence properties outweigh the additional computational cost.
Hence, the overall simulation time of the phase-field fracture problems is greatly reduced, as shown in \cref{sec:conv_study}.

   \section{Benchmark Problems}
\label{sec:benchmark_problems}
In this section, we present five benchmark problems, which we use to demonstrate the convergence properties and the performance of the proposed SPIN methods.
For all test problems, we prepare a finite element mesh by locally refining the expected crack path.
The length scale parameter $l_s$ is then prescribed as $l_s=2h$, where $h$ denotes the size of the largest element from the locally pre-refined zone.
\Cref{tab:material_params} contains a summary of the material parameters, the tolerance used for enforcing the irreversibility condition, and the number of degrees of freedom (dofs) for all numerical examples.

All considered benchmark problems are solved using both, additive and multiplicative, SPIN methods as well as the AM method.
During all numerical experiments, the SPIN methods solve the phase-field and displacement subproblems using the INK method with the biconjugate gradient stabilized (BCGSTAB) method preconditioned with the AMG method, as described in \cref{sec:am}.
The solution of the global preconditioned linear system is obtained using the GMRES method for both SPIN methods.
As discussed in \cref{sec:alg_descr_SPIN}, applying the SPIN preconditioner requires a solution of two linear systems.
Since these linear systems have the same structure as the linear systems arising while solving the subproblems, we also solve them using the BCGSTAB method, preconditioned with the AMG method.

In this work, we consider four variants of AM method, which employ different subproblem solvers and stopping criteria.
More precisely, the AM-ND and AM-ST methods make use of Newton's method with a direct linear solver.
The AM-NK employs Newton's method with the BCGSTAB linear solver, which is preconditioned with the AMG method.
We also consider the AM-INK method, which utilizes the inexact Newton's method with the AMG preconditioned BCGSTAB method.

\Cref{tab:tolerances} provides a summary of all tolerances used for the termination of all employed solution strategies.
We point out that the AM-NK, AM-INK, and AM-ND methods use the same stopping criteria as the SPIN methods, in order to ensure a fair comparison.
In contrast, AM-ST method employs a termination criterium commonly used in the phase-field fracture literature~\cite{bourdin2000numerical}, i.e.,~the method terminates as soon as the change in the phase-field drops below a tolerance $\epsilon_{\text{c}_\text{diff}}$, i.e., $\|c^{(k+1)} - c^{(k)}\|_{\infty} \leqslant \epsilon_{\text{c}_\text{diff}}$.

\begin{remark}
  The stopping criteria used in this work employ standard Euclidean norms and therefore they are mesh-dependent.
  Alternatively, the mesh-independent stopping criteria could be employed, for example by using the correction in the energy norm~\cite{deuflhard2011newton} given as $\|[\matU^{(k+1)}, \matC^{(k+1)}]^{\top}-[\matU^{(k+1)}, \matC^{(k+1)}]^{\top}\|_{F'} \leqslant \epsilon$, or the relative correction in the $L^2$-norm~\cite{gupta2020auto} given as $\max \Big( \frac{\|\matU^{(k+1)}-\matU^{(k)}\|_{L^2}}{ \|\matU^{(k+1)}\|_{L^2} }, \frac{\|\matC^{(k+1)}-\matC^{(k)}\|_{L^2}}{ \|\matC^{(k+1)}\|_{L^2} }\Big) \leqslant \epsilon$.
\end{remark}

The convergence study and exact configuration of all solution strategies are discussed in detail in \cref{sec:conv_study}.
Here, we only report simulation results, in terms of obtained crack patterns.
If all solution strategies converge to the identical solution, only a single snapshot of the result is reported.
Moreover, the evolution of elastic and fracture energy is monitored for all simulations in order to assess the quality of all employed solution strategies.
To increase the readability of our results, we do not report the energy evolution for AM-NK/INK methods, as their global convergence behavior and stopping criterium is identical to those of the AM-ND method.

\begin{table}[t]
  \centering
  \caption{Values of material parameters and the number of dofs for the presented numerical experiments.}
  \label{tab:material_params}
  \begin{tabular}{|c|c|c|c|c|c|}
    \hline
    \multirow{2}{*}{Parameters} & \multicolumn{5}{c|}{Test problems}                                                                                                                                                                                   \\ \cline{2-6}
                                & Tension                            & Shear               & \begin{tabular}{@{}c@{}}Three-point \\ bending \end{tabular} & L-shaped            & \begin{tabular}{@{}c@{}}Asymmetrically \\ notched beam \end{tabular} \\ \hline \hline
    $\pazocal{G}_c$ (kN/mm)     & $2.7 \cdot 10^{-3}$                & $2.7 \cdot 10^{-3}$ & $5.4 \cdot 10^{-4}$                                          & $8.9 \cdot 10^{-5}$ & $10^{-3}$                                                            \\ \hline
    $\lambda$ (kN/mm$^2$)       & $121.15$                           & $121.15$            & $12.00$                                                      & $6.16$              & $12.00$                                                              \\ \hline
    $\mu$ (kN/mm$^2$)           & $80.77$                            & $80.77$             & $8.00$                                                       & $10.95$             & $8.00$                                                               \\ \hline
    $l_s$                       & $0.003$                            & $0.006$             & $0.01$                                                       & $2.0$               & $0.06$                                                               \\ \hline
    $\tau_{\text{irr}}$         & $10^{-2}$                          & $10^{-2}$           & $10^{-2}$                                                    & $10^{-2}$           & $10^{-2}$                                                            \\ \hline
    $\gamma$                    & $9\cdot10^{3}$                     & $4.5\cdot10^{3}$    & $540$                                                        & $0.45$              & $166.67$                                                             \\ \hline
    $\#$ dofs                   & $180,117$                          & $92,004$            & $99,135$                                                     & $141,849$           & $139,020$                                                            \\ \hline
  \end{tabular}
\end{table}

\subsection{Single-edge notched plate under tension}
\label{sec:tension_test}
\begin{figure}[t]
  \begin{minipage}{0.31\linewidth}
    \includegraphics{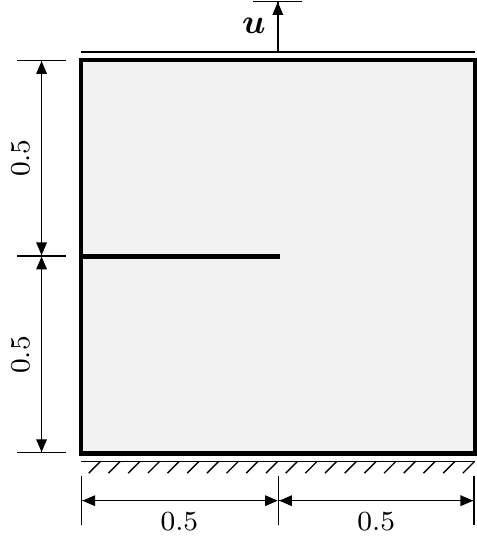}

\end{minipage}
  \hfill
  \begin{minipage}{0.31\linewidth}
    \includegraphics[scale=0.0425]{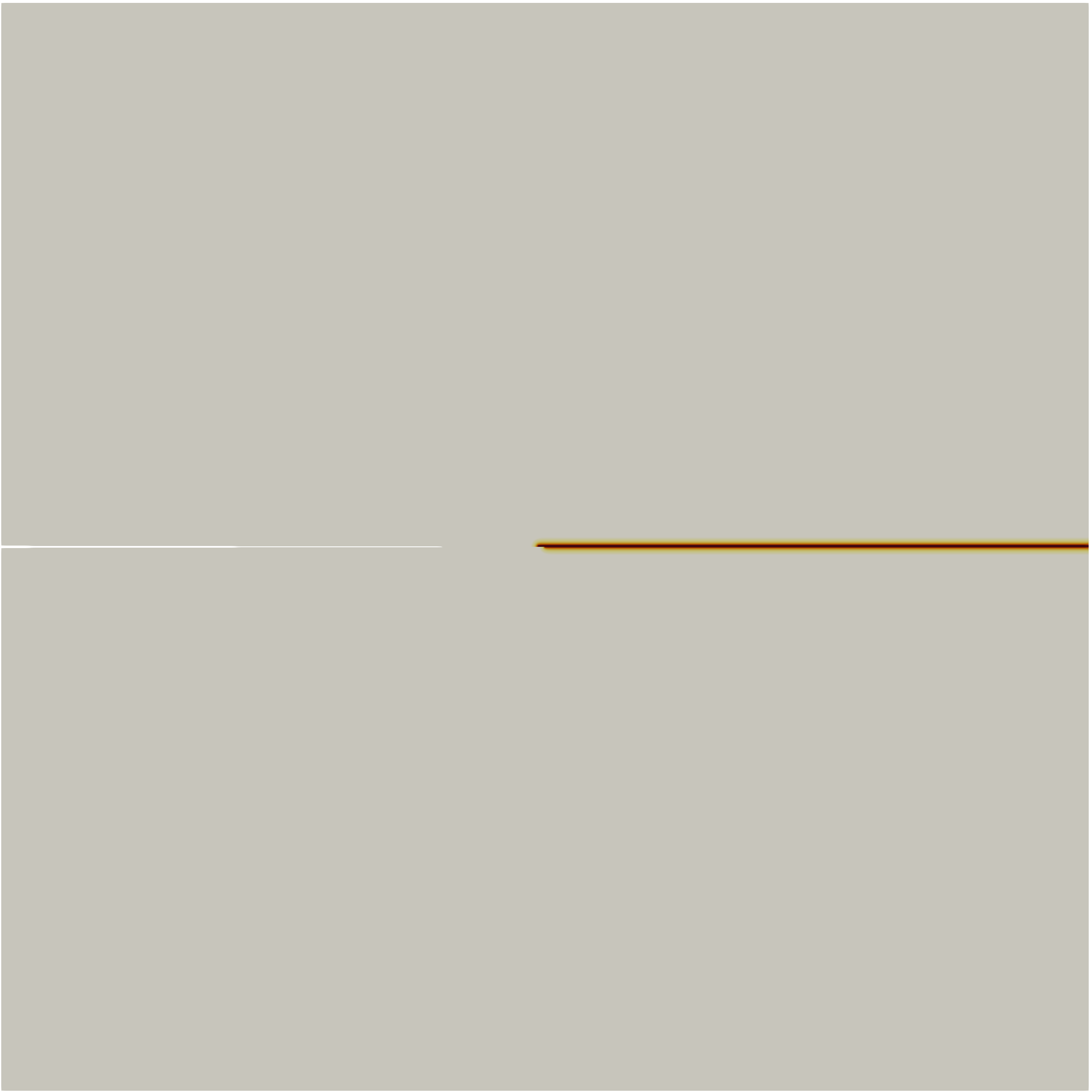}
    \includegraphics[scale=0.0425]{c_scale}
  \end{minipage}
  \hfill
  \begin{minipage}{0.36\linewidth}
    \includegraphics{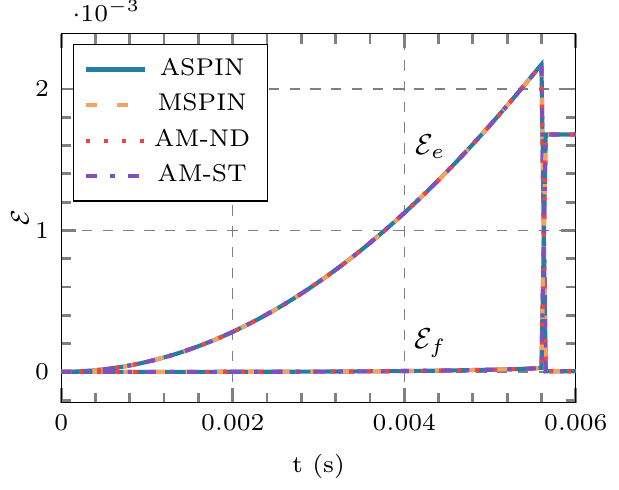}

\end{minipage}
  \caption{Tension test. Left: Geometry (units are in mm) and boundary value problem setup.
    Middle: Simulation result.
    Right: Evolution of elastic $(\pazocal{E}_e)$ and fracture $(\pazocal{E}_f)$ energy over (pseudo-)time for ASPIN, MSPIN, AM-ND and AM-ST solution strategies.}
  \label{fig:tension_setup_sim}
\end{figure}

As the first benchmark problem, we investigate the single-edge notched plate under uniaxial tension.
We consider plate of size $1 \times 1$\,mm with $0.5$\,mm initial notch, subjected to boundary conditions as depicted in \Cref{fig:tension_setup_sim} on the left.
In particular, we impose Dirichlet boundary conditions ${(\uv_D^t)_y := t \bar{u}}$, where $\bar{u}=1\,\text{mm/s}$ on the top of the boundary.
During this experiment, the pseudo time-step is chosen as ${\delta t = 5 \cdot 10^{-5}}$.
For this particular boundary value problem setup, the initial crack is known to extend across the plate within a single loading step.
Therefore, the experiment is widely employed in the literature~\cite{heister2015primal, GERASIMOV2016276, lampron2021efficient} to investigate the robustness of the solution strategies with respect to the ``brutal" crack propagation.

\Cref{fig:tension_setup_sim} in the middle depicts the simulation result, while \Cref{fig:tension_setup_sim} on the right illustrates the evolution of elastic and fracture energy as a function of time.
As we can see, the fracture energy remains zero and elastic energy is increasing until the critical time $t_c= 5.65 \cdot 10^{-3}$\,second is reached.
Once $t_c$ is reached, the ``brutal" crack propagation takes place.
The formation of the crack causes an increase in fracture energy and the degradation in elastic energy.
As we can also observe from \Cref{fig:tension_setup_sim}, the evolution of the elastic and fracture energy is comparable for all four solution strategies.

\subsection{Single-edge notched plate under shear}
\begin{figure}[t]
  \begin{minipage}{0.31\linewidth}
    \includegraphics{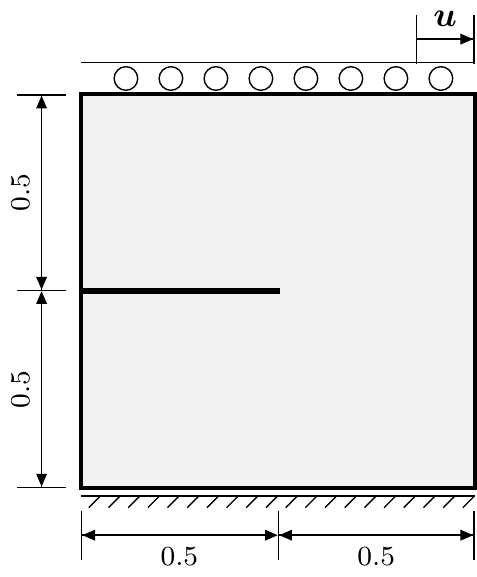}

\end{minipage}
  \hfill
  \begin{minipage}{0.31\linewidth}
    \includegraphics[scale=0.0425]{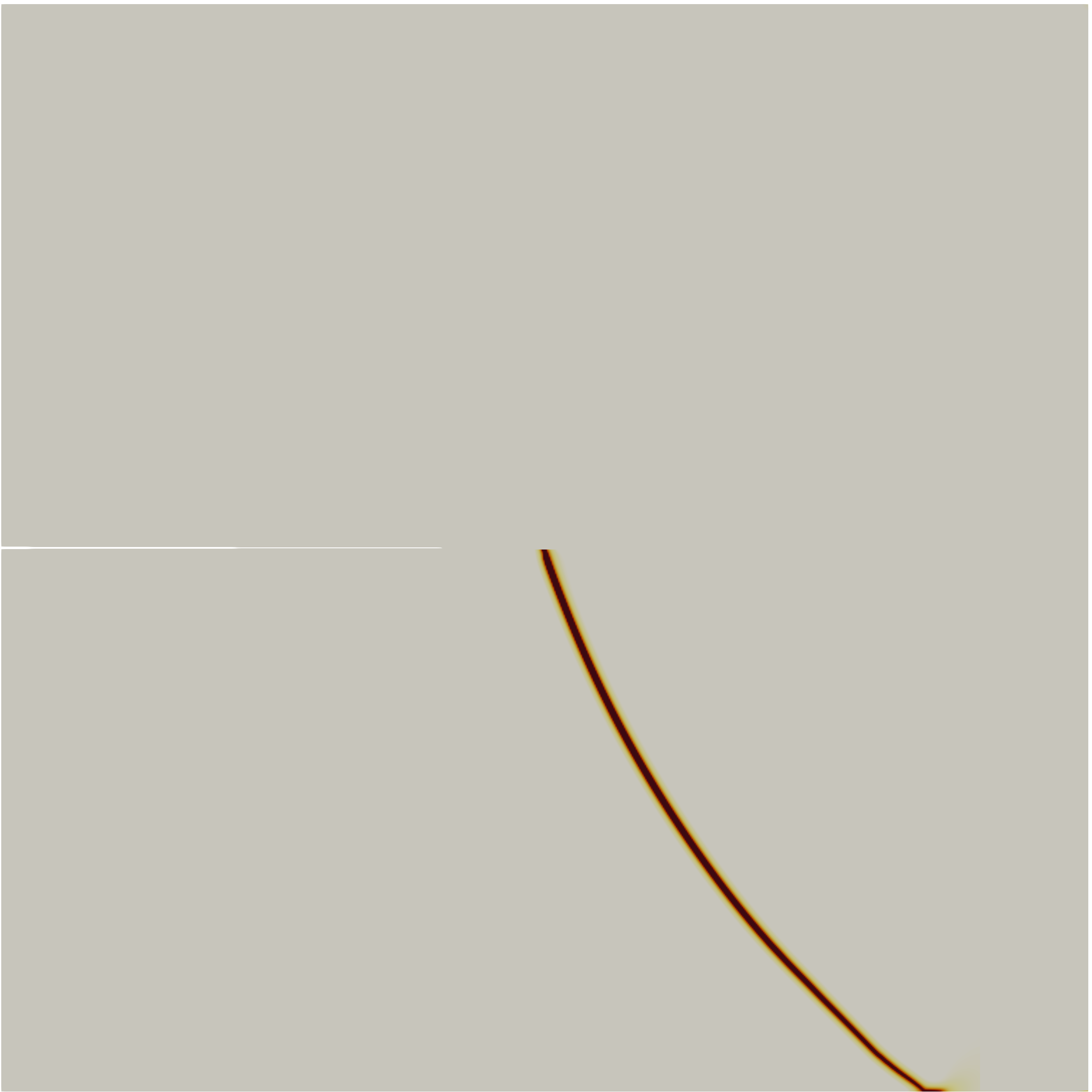}
    \includegraphics[scale=0.0425]{c_scale}
  \end{minipage}
  \hfill
  \begin{minipage}{0.36\linewidth}
    \includegraphics{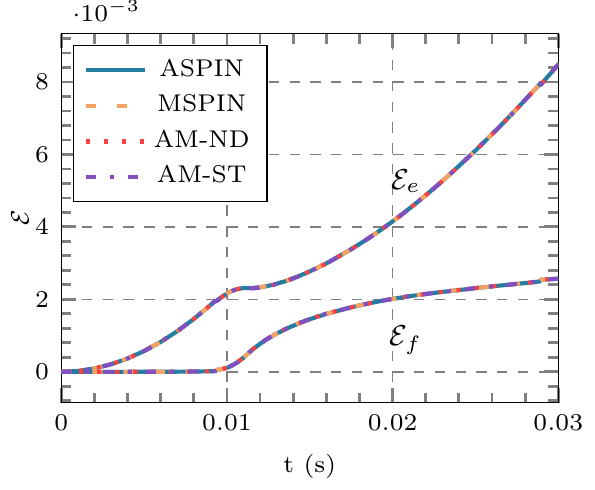}

\end{minipage}
  \caption{Shear test. Left: Geometry (units are in mm) and boundary value problem setup.
    Middle: Simulation result.
    Right: Evolution of elastic $(\pazocal{E}_e)$ and fracture $(\pazocal{E}_f)$ energy over (pseudo-)time for ASPIN, MSPIN, AM-ND and AM-ST solution strategies.}
  \label{fig:shear_setup_sim}
\end{figure}

As our second test, we investigate the single-edge notched plate under shear loading.
Here, we consider the same domain as for the tension test, but different boundary conditions are imposed.
In particular, the bottom of the specimen is held fixed, while the incremental horizontal displacement is prescribed, see also \Cref{fig:shear_setup_sim} on the left.
During this experiment, we prescribe the Dirichlet boundary conditions $(\uv_D^t)_x := t \bar{u}$, where $\bar{u} = 1\,\text{mm/s}$, on the top of the boundary.
Moreover, we employ $\delta t=10^{-3}$ for the first $8$ steps and $\delta t = 7.5 \cdot 10^{-5}$ subsequently until the complete failure occurs.

This example is often considered a canonical one in the phase-field fracture literature~\cite{miehe2010thermodynamically, bourdin2000numerical, heister2015primal, amor2009regularized}, as it enables to test the capability of the model to ensure the asymmetric nucleation in tension and compression.
As described in \cref{sec:pf}, we employ the phase-field model with the spectral decomposition of the strain tensor.
In this case, the initial crack is extended only in the bottom-right part of the specimen and forms a curved crack, see \Cref{fig:shear_setup_sim} in the middle.
In contrast to the tension test, the crack propagation is stable and occurs consistently over multiple loading steps.
This enables us to assess the convergence properties and the efficiency of the solution strategies over time.

\Cref{fig:shear_setup_sim} on the right depicts the evolution of the elastic and fracture energies as a function of time.
As we can see, the behavior of the elastic and fracture energies coincides for all three solution strategies.
We can also observe that complete failure is never attained.
After the cracked specimen reaches a stage, where no further evolution of the phase-field in the lower-right corner is possible, the behavior corresponds to the linearly elastic response of the cracked specimen clamped at the undamaged bottom-right portion of the boundary.
Comparable numerical results have also been reported in~\cite{miehe2010thermodynamically, ambati2015review}.

\subsection{Three-point bending test}
\label{sec:3pt_bending_description}
\begin{figure}[t]
  \begin{minipage}{0.49\linewidth}
    \includegraphics{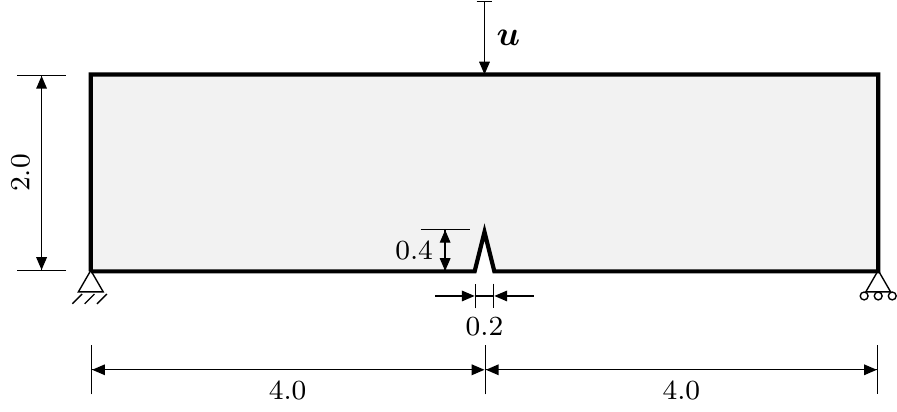}

    \includegraphics[scale=0.022]{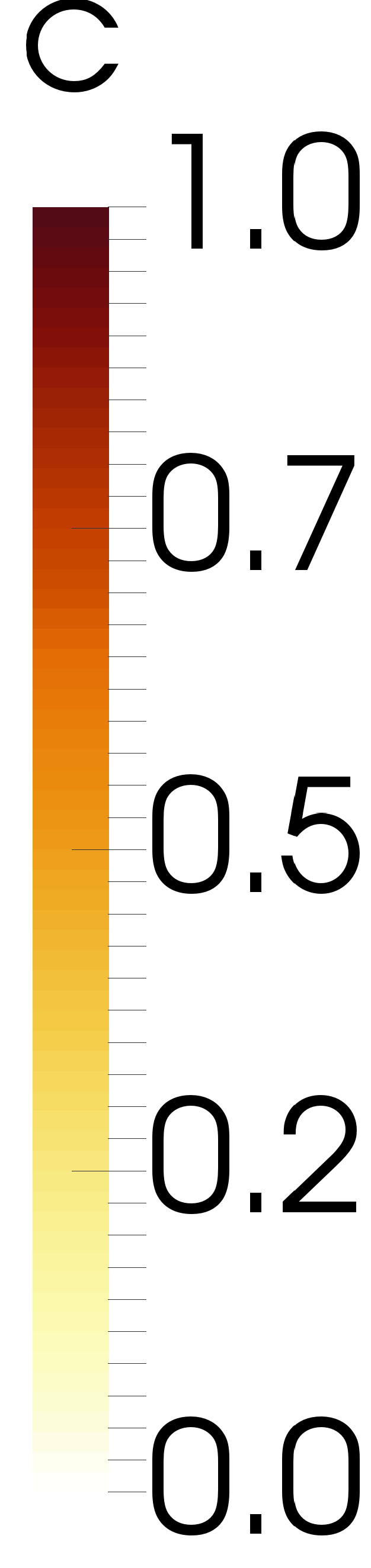}
    \hspace{0.19cm}
    \includegraphics[scale=0.09]{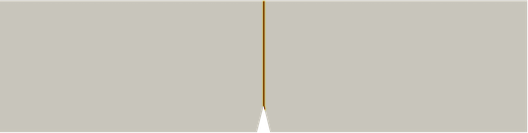}
\end{minipage}
  \hfill
  \begin{minipage}{0.42\linewidth}
    \vspace{0.6cm}
    \includegraphics{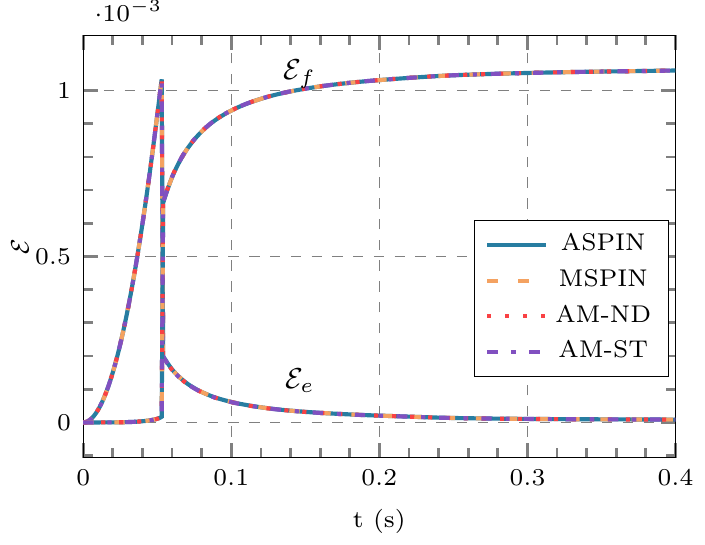}

\end{minipage}

  \caption{
    Three-point bending test.
    Left top: Geometry (units are in mm) and boundary value problem setup.
    Left bottom: Simulation result.
    Right: Evolution of elastic $(\pazocal{E}_e)$ and fracture $(\pazocal{E}_f)$ energy over (pseudo-)time for ASPIN, MSPIN, AM-ND and AM-ST solution strategies.}
  \label{fig:3pt_bending_setup}
\end{figure}

In this experiment, we model the three-point bending of a simply supported beam, as presented in~\cite{ambati2015review, miehe2010phase}.
\Cref{fig:3pt_bending_setup} on the left illustrates the geometry and the boundary value problem setup.
The time-dependent boundary conditions are prescribed on the top as $(\uv_D^t)_y := t \bar{u}$, where $\bar{u}=-1\,\text{mm/s}$.
During this experiment, we use constant pseudo time-step $\delta t=10^{-3}$.
For this particular loading scenario, the crack propagates from the notch in a straight line towards the upper boundary, see also \Cref{fig:3pt_bending_setup} on the left.
The crack propagation exhibits two phases, which can be observed from the curves illustrating the energy evolution over time (\Cref{fig:3pt_bending_setup} on the right).
The first phase consists of the crack propagation in the tension mode, thus we can observe ``brutal" crack propagation at time $t=0.058$\,seconds.
This is followed by the second phase, where the crack propagates gradually.
As for the previous examples, the same energy evolution is obtained with all three solution strategies under consideration.

Since this experiment contains both ``brutal" and gradual crack propagation, we also decided to use it in order to assess our implementation of the phase-field fracture model.
To this aim, we study the converge behavior of the phase-field fracture model under refinement.
We consider a hierarchy of three meshes, obtained by uniformly refining the adaptively refined finite element mesh considered in the experiments above.
With each refinement step, the mesh size decreases by a factor of two, which allows us to also reduce the associated value of the length-scale parameter $l_s$.
We can observe from the \Cref{fig:3pt_bending_conv_ref} on the left, that the mesh resolution and the value of length-scale parameter influence the time of crack propagation, as well as the values of the elastic and fracture energy.
As expected from the $\Gamma$-convergence theory~\cite{ambrosio1990approximation}, as the value of the length-scale parameter~$l_s$ approaches zero, the fracture energy converges towards $\pazocal{E}_f = \pazocal{G}_c W$, where $W$ symbolizes the area of the fracture surface.

\subsection{L-shaped panel test}
\begin{figure}[t]
  \includegraphics{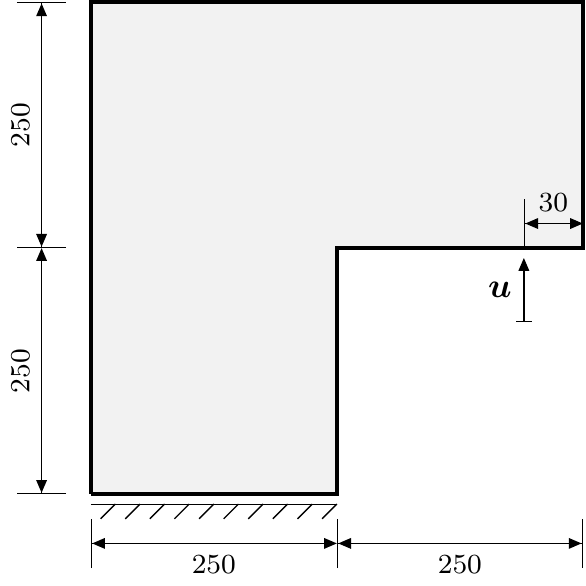} \hfill  \includegraphics{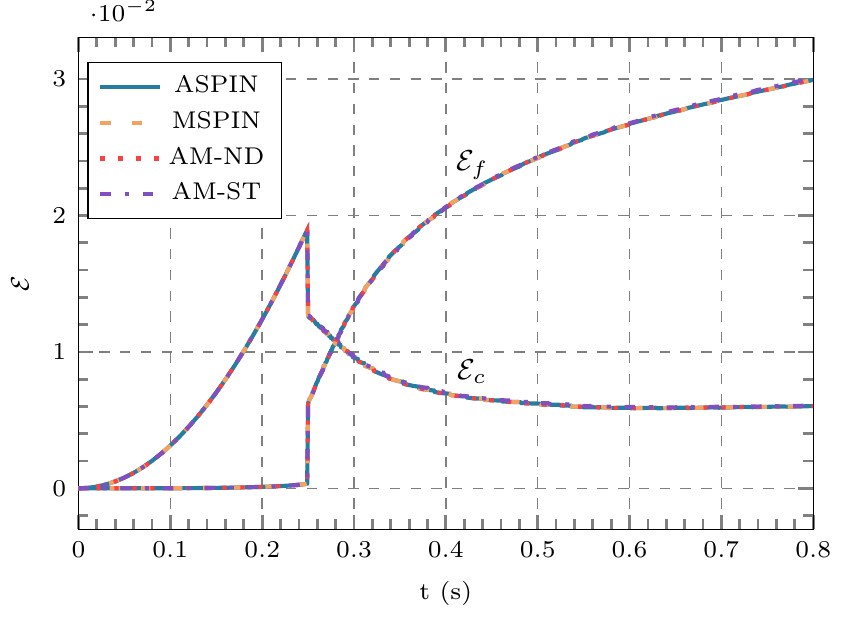}

\caption{The L-shaped panel test. Left: Geometry (units are in mm) and boundary value problem setup.
    Right: Evolution of elastic $(\pazocal{E}_e)$ and fracture $(\pazocal{E}_f)$ energy over (pseudo-)time for ASPIN, MSPIN, AM-ND and AM-ST solution strategies.}
  \label{fig:l_shape_setup}
\end{figure}
\begin{figure}[htb!]
  \includegraphics{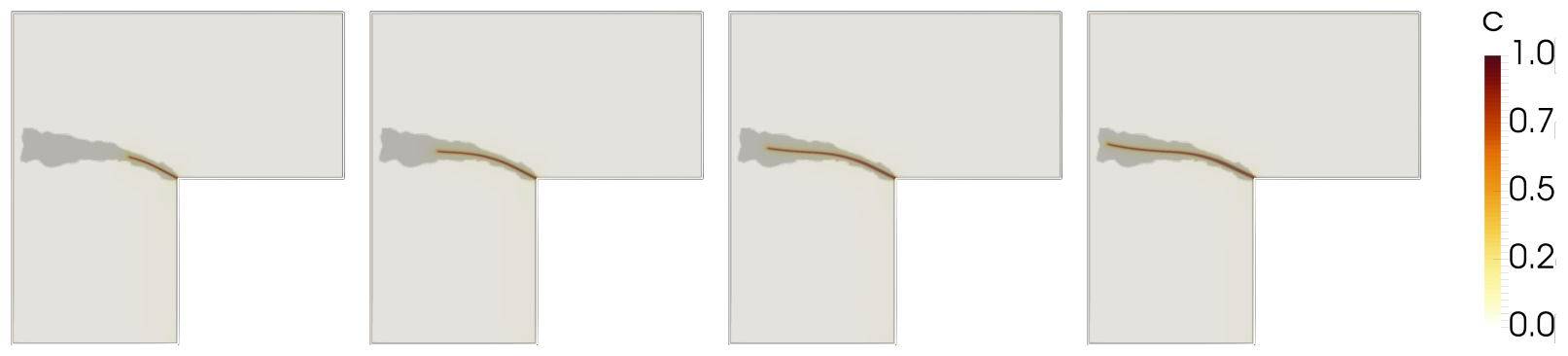}

\caption{Simulation result for L-shaped panel test.
    The simulation result (the value of the phase-field) overlays the experimentally obtained results reported in~\cite{winkler2001lastlastuntersuchungen}.
    The crack evolution depicted at time $t \in \{0.27, 0.365, 0.515, 0.792 \} (\text{s})$, from left to the right.}
  \label{fig:l_shape_sol}
\end{figure}

Our next experiment considers crack propagation in an L-shaped panel.
Here, we adopt the geometry and material parameters from~\cite{ambati2015review}, which mimics the experimental setup presented in~\cite{winkler2001lastlastuntersuchungen}.
The details about geometry and the problem setup can be found in \Cref{fig:l_shape_setup} on the left.
The simulation is performed by imposing the Dirichlet boundary conditions of the following form: $(\uv_D^t)_y := t \bar{u}$, where $\bar{u}=1\,\text{mm/s}$.
The pseudo time-step $\delta t$ is set to $10^{-2}$ for the first $20$ steps and to $10^{-3}$ afterwords.
\Cref{fig:l_shape_sol} presents the obtained results, i.e., the simulated crack path at different loading steps.
As we can see, the simulated crack path is in agreement with the experimental results reported in~\cite{winkler2001lastlastuntersuchungen}.
Furthermore, \Cref{fig:l_shape_setup} on the right demonstrates that the identical evolution of the energies is obtained irrespectively of the choice of the solution strategy at the hand.

\subsection{Asymmetrically notched beam test}
\begin{figure}[t]
  \begin{minipage}{0.48\linewidth}
    \includegraphics{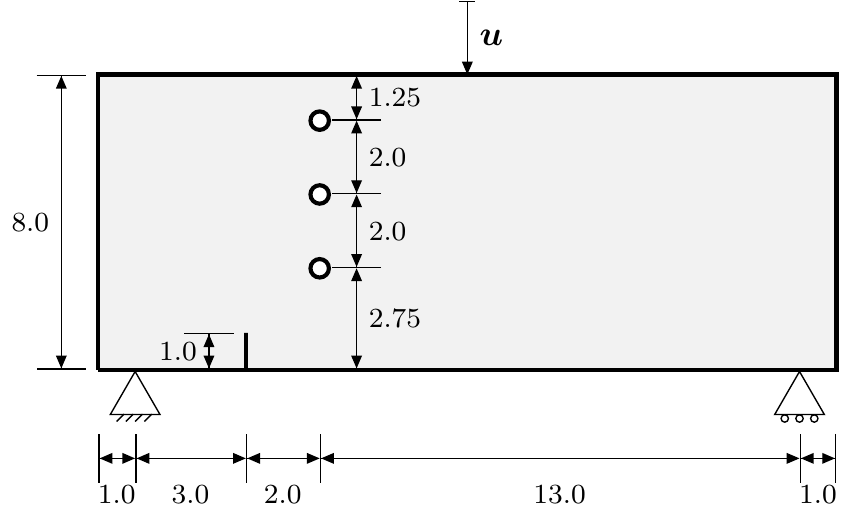}

\end{minipage}
  \hfill
  \begin{minipage}{0.48\linewidth}
    \includegraphics{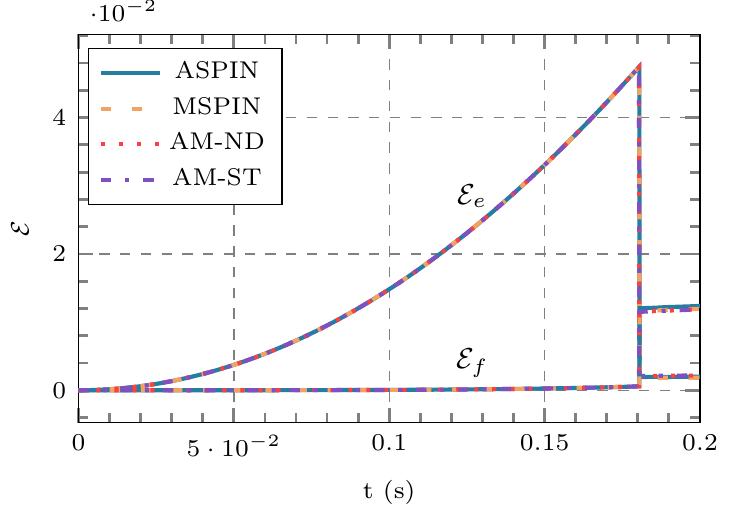}

\end{minipage}
  \caption{Asymmetrically notched beam test.
    Left: Geometry (units are in mm) and boundary value problem setup.
    Right: Evolution of elastic $(\pazocal{E}_e)$ and fracture $(\pazocal{E}_f)$ energy over (pseudo-)time for ASPIN, MSPIN, AM-ND and AM-ST solution strategies.}
  \label{fig:assym_bending_setup}
\end{figure}
\begin{figure}[t]
  \includegraphics[scale=0.085]{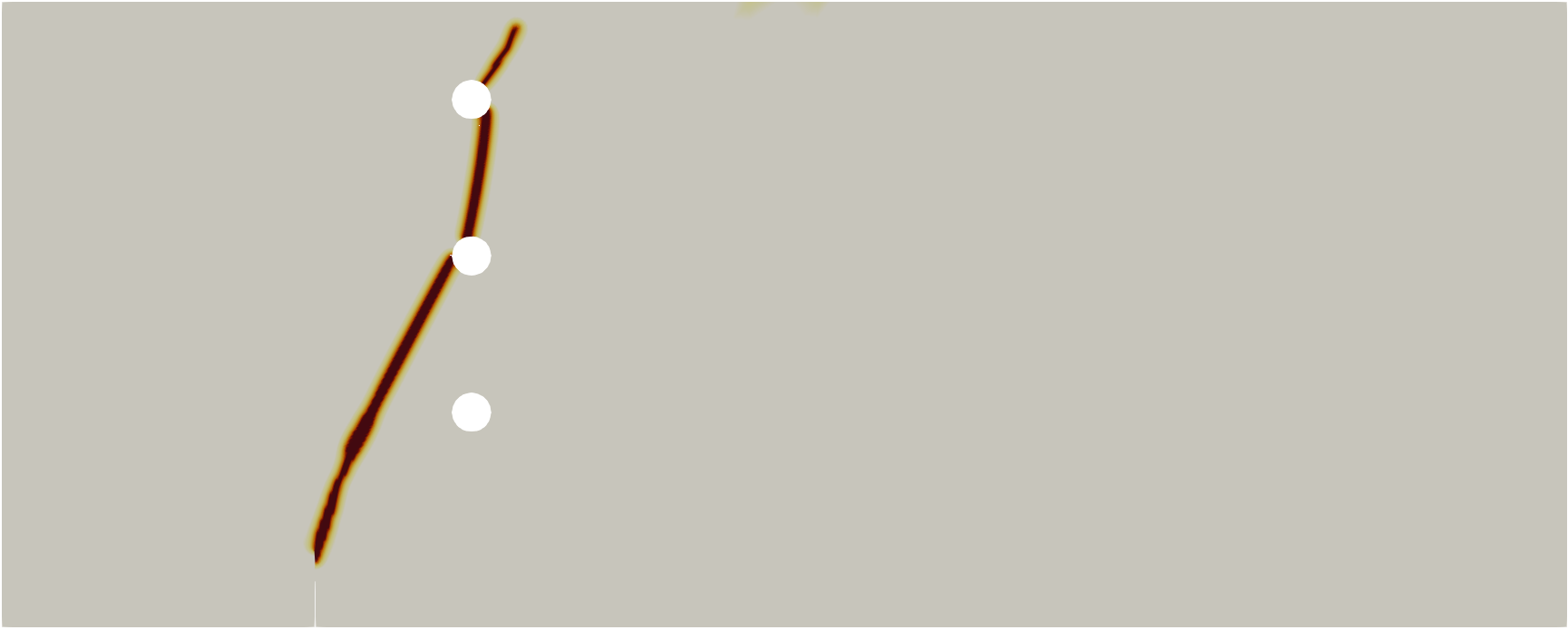}
  \includegraphics[scale=0.085]{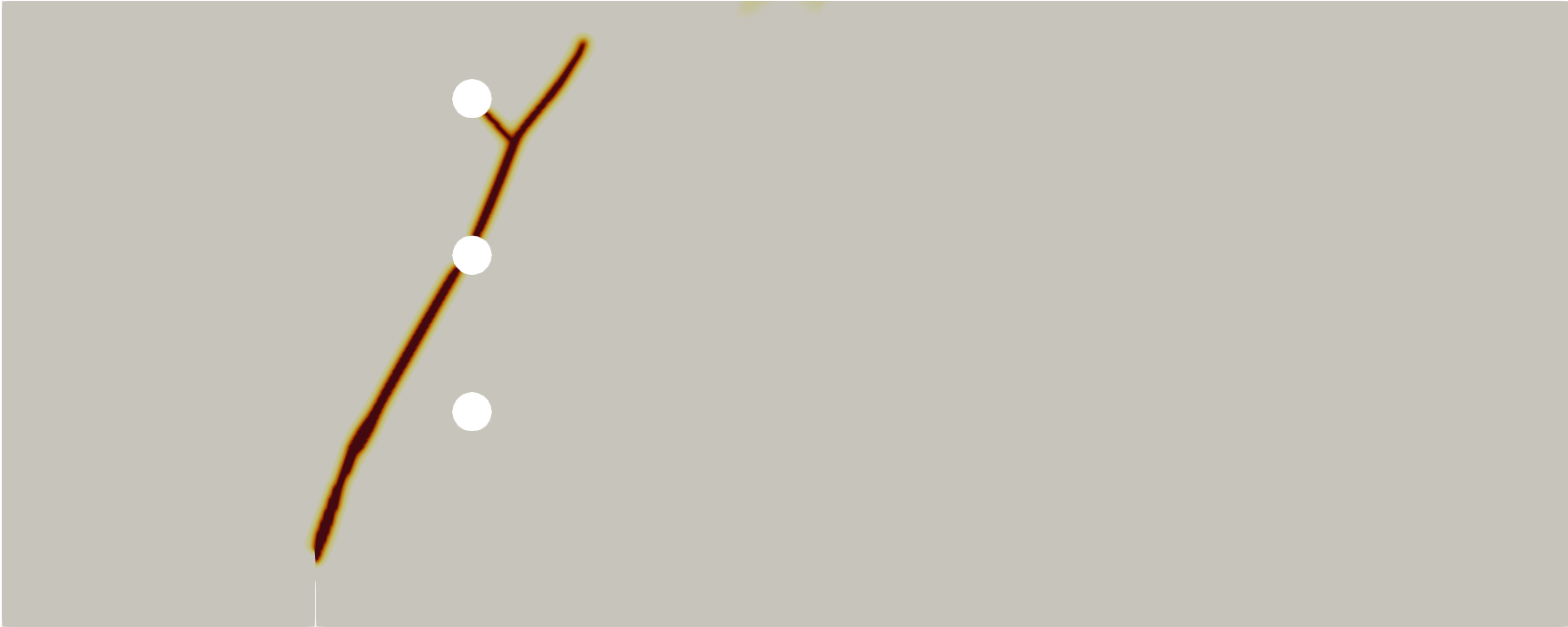}
  \includegraphics[scale=0.085]{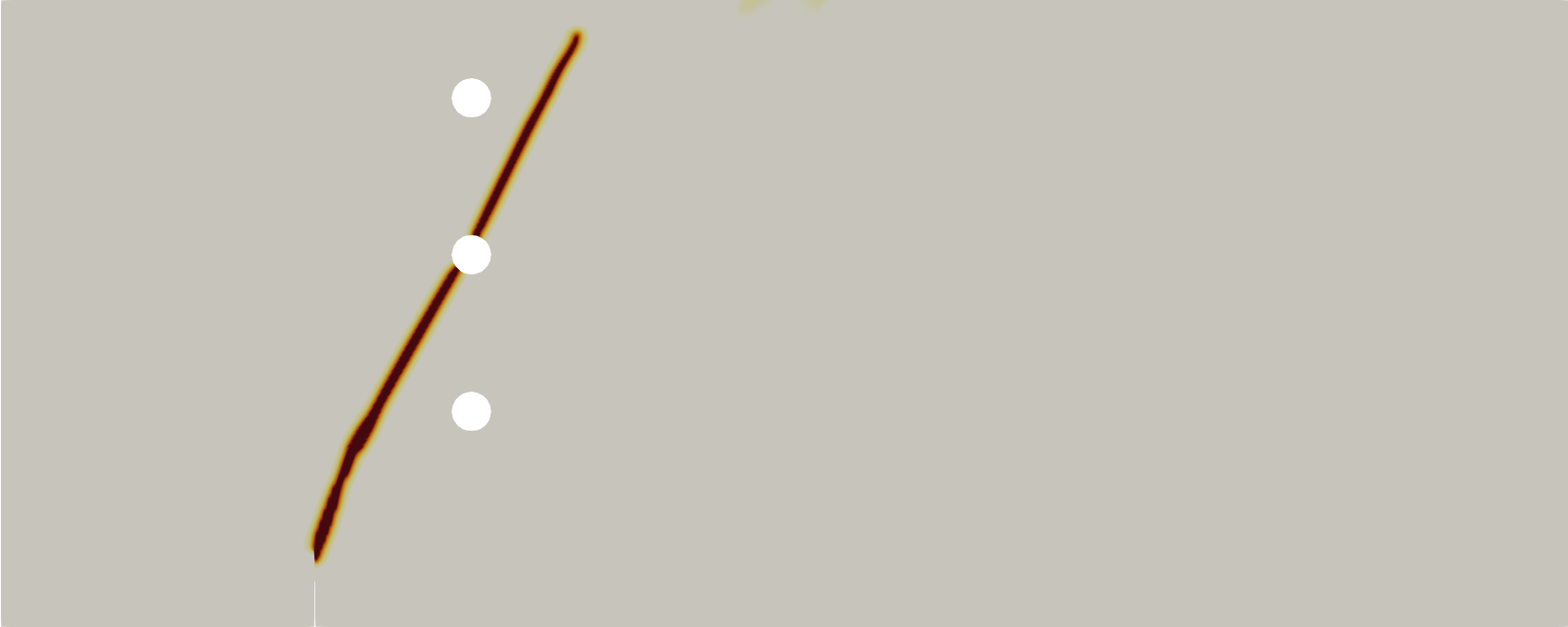}
  \includegraphics[scale=0.023]{c_scale}
  \caption{Simulation results for asymmetrically notched beam test.
    Left: The solution obtained using the AM methods.
    Middle: The solution obtained using the ASPIN method.
    Right: The solution obtained using the MSPIN method.}
  \label{fig:asymetric_sol}
\end{figure}
In the last benchmark problem, we consider an asymmetrically notched beam, with the notch and three holes in the left part, as illustrated in \Cref{fig:assym_bending_setup}.
The beam is subjected to displacement loading in the center of the top boundary while being supported at two points at the bottom boundary.
Thus, we prescribe the Dirichlet boundary conditions as $(\uv_D^t)_y := t \bar{u}$, where $\bar{u}=-1\,\text{mm/s}$.
The time-step is chosen as $\delta t = 10^{-3}$ for the first $160$ steps and as $\delta t = 10^{-4}$ for the remaining part of the simulation.
The material parameters are selected as in~\cite{hesch2014thermodynamically, brun2020iterative, storvik2021accelerated} in order to mimic the experimental setup from~\cite{bittencourt1996quasi}.

For a given setup, the experimental results have shown that crack should propagate towards the second hole, which we can also observe in our numerical experiments, as depicted in \Cref{fig:asymetric_sol}.
However, we obtain three distinct numerical solutions.
This is not surprising, as the energy functional is non-convex and therefore the minimization admits multiple solutions.
Interestingly, all variants of the AM method converge to the same solution, while the ASPIN and MSPIN method find two different crack patterns.
Close examination of the energies suggests that the MSPIN method converges to the solution with the smallest value of the total potential energy $\pazocal{E}$.
In particular, the value of the total potential energy is $\pazocal{E} \approx 1.379\cdot 10^{-2}$ for the MSPIN method, $\pazocal{E} \approx 1.402\cdot 10^{-2}$ for the AM methods and, $\pazocal{E} \approx 1.435\cdot 10^{-2}$ for the ASPIN method.
However, we also point out that additional solutions might exist.
The recent work of Gerasimov et al.~\cite{gerasimov2020stochastic} provides a methodology for uncovering different solutions by modifying loading increments and/or by perturbing the finite element mesh.
To the best of our knowledge, an approach that would allow us to discover all possible solutions for a given finite element mesh and loading scenario has not been developed so far in the phase-field fracture literature.

   \section{Convergence and performance study}
\label{sec:conv_study}
In this section, we investigate the convergence properties of the proposed ASPIN and MSPIN methods, described in \Cref{alg:spin}.
The proposed methods are compared to the alternate minimization (AM) method, summarized in \Cref{alg:AM}.
The assessment of the solution strategies is performed in terms of the number of iterations, execution times, and memory requirements.
Moreover, we also demonstrate the robustness and the performance of the SPIN methods with respect to the refinement level.

\subsection{Implementation}
Our implementation of the phase-field fracture model is based on the finite-element library MOOSE~\cite{gaston2009moose,chakraborty2016multi}.
All solution strategies and their components were implemented using the open-source library Utopia~\cite{utopia,zulian2021large}.
Utopia is an embedded domain-specific language, which allows us to express complex numerical procedures with few lines of code, whereas the complexity of the parallelization/hardware-specific optimizations is hidden in the different back-ends.
For this work, we have utilized linear algebra provided by PETSc~\cite{balay2014petsc} back-end.
The direct solver is provided by MUMPS library~\cite{amestoy2000mumps}.
Moreover, we employ the algebraic multigrid (AMG) method provided by the BoomerAMG package from Hypre library~\cite{yang2002boomeramg}.
The employed AMG method is configured using default settings.
The developed code, including implementation of the model, solution strategies, and benchmark problems, is publicly available at~\url{https://bitbucket.org/alena_kopanicakova/pf_frac_spin}.

The presented numerical experiments were performed on a machine equipped with two {Intel Xeon~E5-2650~v3} processors (10 cores per processor) and memory of size~$64$\,GB.
Thus, the numerical experiments were performed using $20$ cores.

\subsection{Sensitivity of the SPIN methods to choice of $\epsilon_{\text{app}\_{\text{lin}}}$}
\begin{figure}[t]
  \centering
  \includegraphics{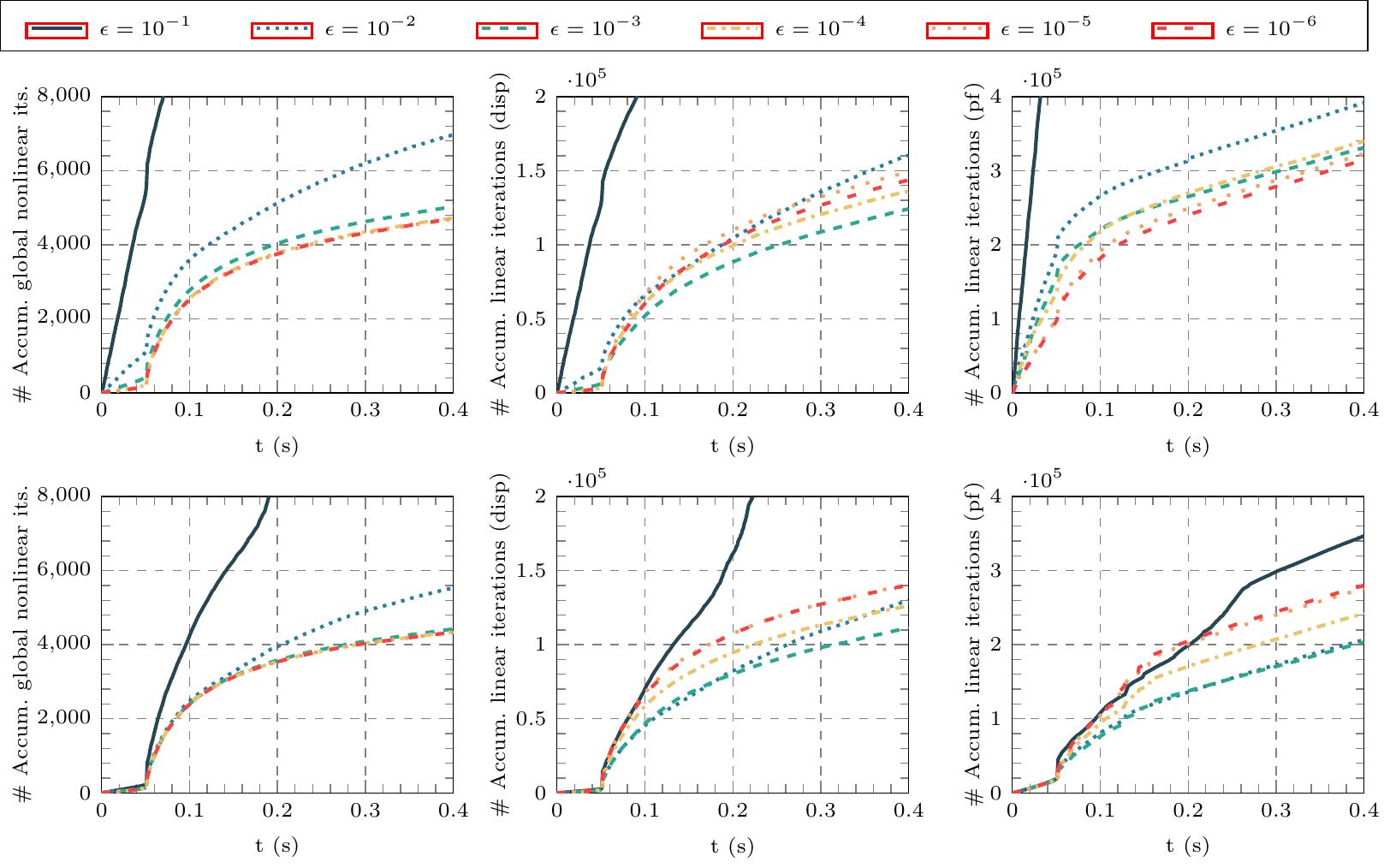}

\caption{Three-point bending test. The convergence behavior of the ASPIN (top row) and MSPIN (bottom row) methods with respect to different values of $\epsilon_{\text{app}\_{\text{lin}}}$.
    Left: Accumulated number of global nonlinear iterations.
    Middle: The accumulated number of linear iterations required for solving the linear systems associated with the displacement field.
    Right: The accumulated number of linear iterations required for solving the linear systems associated with the phase field.}
  \label{fig:3pt_bending_action_rtol_impact_mspin}
\end{figure}

As discussed in \cref{sec:alg_descr_SPIN}, applying the SPIN preconditioners requires a solution of two linear systems.
However, it is not necessary to solve these linear systems exactly, but only with some relative tolerance $\epsilon_{\text{app}\_{\text{lin}}}$.
\Cref{fig:3pt_bending_action_rtol_impact_mspin} demonstrates the convergence behavior of SPIN methods with respect to different values of $\epsilon_{\text{app}\_{\text{lin}}}$.
For this study, the experiments are performed using the three-point bending test with $6,417$ dofs, where the mesh size is approximately $4$ times larger than the mesh size reported in \Cref{tab:tolerances}.
We also remark that the length-scale parameter is adjusted accordingly for the given mesh size.
As we can observe from the obtained results, it is necessary to choose $\epsilon_{\text{app}\_{\text{lin}}}\leqslant 10^{-2}$, as larger values of $\epsilon_{\text{app}\_{\text{lin}}}$ can cause the stagnation or even divergence of SPIN method.
We can also see that for $\epsilon_{\text{app}\_{\text{lin}}} \leqslant 10^{-3}$, the number of required global ASPIN/MSPIN iterations remains constant.
However, as the tolerance $\epsilon_{\text{app}\_{\text{lin}}}$ decreases, the number of required linear iterations varies, which causes an increase in the computational cost of the overall method.
Identical behavior can be also observed in terms of execution time, reported in \Cref{tab:computational_time_epsilon}.
As we can see, the MSPIN method configured with $\epsilon_{\text{app}\_{\text{lin}}} = 10^{-4}$ is the fastest.
While for the ASPIN method the fastest performance can be observed with $\epsilon_{\text{app}\_{\text{lin}}} = 10^{-5}$.
As a compromise, all numerical results presented in this work employ $\epsilon_{\text{app}\_{\text{lin}}} = 10^{-4}$.

\begin{table}[t]
  \caption{Summary of the numerical tolerances used for the termination of all solution strategies.}
  \label{tab:tolerances}
  \centering
  \begin{tabular}{ |l|l|l|}
    \hline
    \multicolumn{1}{|c|}{Symbol}                        & \multicolumn{1}{c|}{Description of the termination criteria}                     & \multicolumn{1}{c|}{Value} \\ \hline
    $\epsilon_{\text{abs}\_\text{glob}\_{\text{nonl}}}$ & Absolute residual norm of global algorithms (AM/ASPIN/MSPIN).                    & $10^{-7}$                  \\ \hline
    $\epsilon_{\text{rel}\_\text{glob}\_{\text{nonl}}}$ & Relative residual norm of global algorithms (AM/ASPIN/MSPIN).                    & $10^{-6}$                  \\ \hline
    $\epsilon_{\text{abs}\_\text{sub}\_{\text{nonl}}}$  & Absolute residual norm of ND/NK/INK methods, used for solving subproblems.       & $10^{-7}$                  \\ \hline
    $\epsilon_{\text{rel}\_\text{sub}\_{\text{nonl}}}$  & Relative residual norm of ND/NK/INK methods, used for solving subproblems.       & $10^{-6}$                  \\ \hline
    $\epsilon_{\text{abs}\_\text{lin}}$                 & Absolute residual norm of Krylov methods, used inside of NK method.              & $10^{-9}$                  \\ \hline
    $\epsilon_{\text{rel}\_\text{lin}}$                 & Relative residual norm of Krylov methods, used inside of NK method.              & $10^{-9}$                  \\ \hline
    $\epsilon_{\text{app}\_{\text{lin}}}$               & Relative residual norm of Krylov methods, used for applying SPIN preconditioner. & $10^{-4}$                  \\ \hline
    $\epsilon_{\text{c}\_{\text{diff}}}$                & Absolute difference in change of phase-field measured in infinity norm.          & $10^{-4}$                  \\ \hline
  \end{tabular}
\end{table}

\begin{table}[t]
  \caption{The execution time (mins) of the ASPIN and MSPIN method required for the simulation of the three-point bending test.
    The time is reported for different values of $\epsilon_{\text{app}\_{\text{lin}}}$.}
  \label{tab:computational_time_epsilon}
  \centering
  \begin{tabular}{|c|c|c|c|c|c|c|}
    \hline
    $\epsilon_{\text{app}\_{\text{lin}}}$ & $10^{-1}$ & $10^{-2}$ & $10^{-3}$ & $10^{-4}$ & $10^{-5}$ & $10^{-6}$ \\ \hline \hline
    ASPIN                                 & $206.05$  & $24.26$   & $16.13$   & $15.02$   & $14.74$   & $15.80$   \\ \hline
    MSPIN                                 & $73.79$   & $18.65$   & $15.60$   & $14.26$   & $14.98$   & $14.78$   \\ \hline
  \end{tabular}
\end{table}

\subsection{Convergence properties of solution strategies}
\label{sec:solver_perf_summary_plot_and_table}
\begin{figure}
  \centering
  \includegraphics{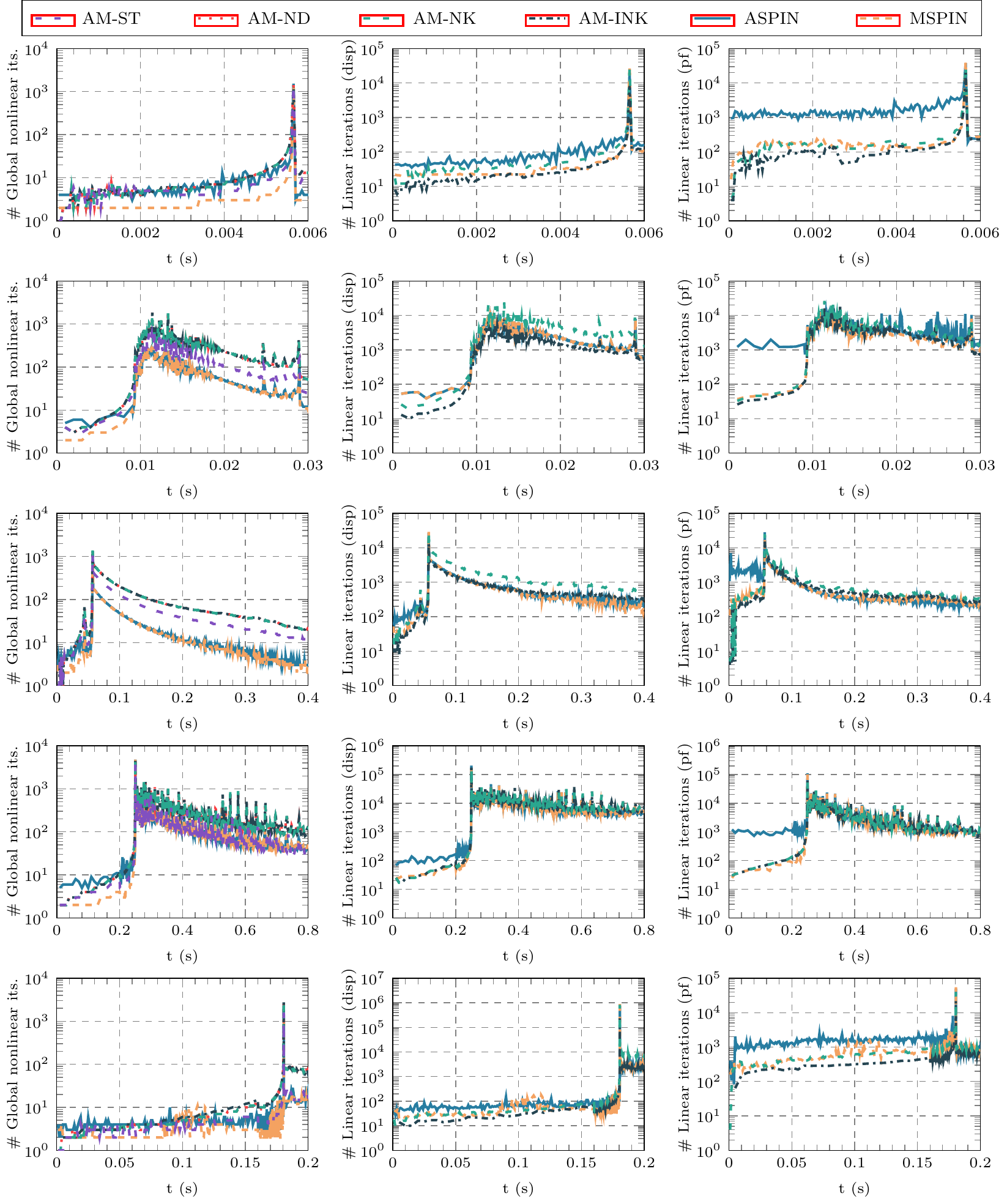}

  \caption{Convergence history of ASPIN, MSPIN, AM-ST, AM-ND, AM-NK, and AM-INK solution strategies.
    The experiments are performed for tension, shear, three-point bending, L-shaped panel, and asymmetrically notched beam tests (from top to bottom).
    Left: The number of global nonlinear iterations over time.
    Middle: The number of linear iterations performed for solving linear systems related to the displacement subproblem.
    Right: The number of linear iterations performed for solving linear systems related to the phase-field subproblem.  }
  \label{fig:all_solvers_perf}
\end{figure}

\begin{table}[ht!]
  \caption{Summary of computational cost for ASPIN, MSPIN, AM-ST, AM-ND, AM-NK, and AM-INK methods.
  Summary information is expressed in terms of the accumulated number of linear and nonlinear iterations required to solve the coupled problem as well as the internal displacement and phase-field subproblems.
  The symbols $*,\dagger,+$ denote convergence to a different solution.}
  \label{tab:overall_results}
  \centering
  \begin{tabular}{|c|l|rr|rr|rr|}
    \hline
    \multirow{2}{*}{Example}                             & \multicolumn{1}{c|}{ \multirow{2}{*}{Solver}} & \multicolumn{2}{c|}{Global} & \multicolumn{2}{c|}{Disp. subproblem} & \multicolumn{2}{c|}{PF subproblem}                                         \\
                                                         &                                               & Nonlinear                   & Linear                                & Nonlinear                          & Linear      & Nonlinear & Linear      \\ \hline \hline
    \multirow{6}{*}{Tension}                             & AM-ND                                         & $2,507$                     & --                                    & $5,421$                            & --          & $ 4,422$  & --          \\ \cline{2-8}
                                                         & AM-NK                                         & $2,507$                     & --                                    & $5,421$                            & $34,670$    & $ 4,422$  & $ 57,516$   \\ \cline{2-8}
                                                         & AM-INK                                        & $2,526$                     & --                                    & $6,373$                            & $12,336$    & $ 4,422$  & $ 34,330$   \\ \cline{2-8}
                                                         & AM-ST                                         & $2,135$                     & --                                    & $5,052$                            & --          & $ 4,033$  & --          \\ \cline{2-8}
                                                         & ASPIN                                         & $2,375$                     & $4,925$                               & $5,490$                            & $31,580$    & $11,322$  & $222,009$   \\ \cline{2-8}
                                                         & MSPIN                                         & $1,823$                     & $5,336$                               & $5,125$                            & $37,034$    & $ 6,328$  & $ 63,237$   \\ \hline \hline

    \multirow{6}{*}{Shear}
                                                         & AM-ND                                         & $81,940$                    & --                                    & $91,066$                           & --          & $86,908$  & --          \\ \cline{2-8}
                                                         & AM-NK                                         & $81,940$                    & --                                    & $91,066$                           & $1,437,712$ & $86,909$  & $1,430,913$ \\ \cline{2-8}
                                                         & AM-INK                                        & $82,106$                    & --                                    & $115,785$                          & $  459,834$ & $87,088$  & $  980,977$ \\ \cline{2-8}
                                                         & AM-ST                                         & $40,626$                    & --                                    & $50,023$                           & --          & $45,710$  & --          \\ \cline{2-8}
                                                         & ASPIN                                         & $20,250$                    & $59,584$                              & $26,237$                           & $711,607$   & $45,312$  & $1,212,635$ \\ \cline{2-8}
                                                         & MSPIN                                         & $19,824$                    & $58,766$                              & $25,318$                           & $697,591$   & $32,732$  & $849,583$   \\ \hline \hline

    \multirow{6}{2.5cm}{\centering Three-point  bending} & AM-ND                                         & $36,037$                    & --                                    & $42,249$                           & --          & $39,232$  & --          \\ \cline{2-8}
                                                         & AM-NK                                         & $36,037$                    & --                                    & $42,250$                           & $600,923$   & $39,232$  & $380,221$   \\ \cline{2-8}
                                                         & AM-INK                                        & $36,045$                    & --                                    & $42,622$                           & $290,456$   & $39,243$  & $291,308$   \\ \cline{2-8}
                                                         & AM-ST                                         & $22,407$                    & --                                    & $28,716$                           & --          & $25,639$  & --          \\ \cline{2-8}
                                                         & ASPIN                                         & $ 8,875$                    & $26,021$                              & $12,968$                           & $304,249$   & $22,096$  & $362,201$   \\ \cline{2-8}
                                                         & MSPIN                                         & $ 8,409$                    & $25,452$                              & $12,168$                           & $287,991$   & $17,299$  & $244,672$   \\ \hline \hline

    \multirow{6}{*}{L-shaped panel}
                                                         & AM-ND                                         & $152,121$                   & --                                    & $168,475$                          & --          & $166,997$ & --          \\ \cline{2-8}
                                                         & AM-NK                                         & $152,121$                   & --                                    & $168,475$                          & $5,549,803$ & $166,998$ & $1,774,117$ \\ \cline{2-8}
                                                         & AM-INK                                        & $158,205$                         & --                                    & $175,214 $                                & $1,989,176 $         & $167,005 $       & $948,725$         \\ \cline{2-8}
                                                         & AM-ST                                         & $64,833$                    & --                                    & $81,656$                           & --          & $ 79,357$ & --          \\ \cline{2-8}
                                                         & ASPIN                                         & $67,653$                    & $201,748$                             & $71,651$                           & $4,433,556$ & $111,782$ & $1,981,303$ \\ \cline{2-8}
                                                         & MSPIN                                         & $66,326$                    & $198,405$                             & $70,099$                           & $4,370,314$ & $ 97,412$ & $1,719,563$ \\ \hline \hline

    \multirow{6}{2.5cm}{\centering Asymmetrically notched beam}
                                                         & AM-ND$^\ast$                                  & $21,961$                    & --                                    & $29,959$                           & --          & $29,989$  & --          \\ \cline{2-8}
                                                         & AM-NK$^\ast$                                  & $21,961$                    & --                                    & $29,515$                           & $2,262,266$ & $29,958$  & $559,184$   \\ \cline{2-8}
                                                         & AM-INK$^\ast$                                 & $21,376$                    & --                                    & $29,543$                           & $  707,202$ & $28,725$  & $287,935$   \\ \cline{2-8}
                                                         & AM-ST$^\ast$                                  & $ 7,200$                    & --                                    & $13,957$                           & --          & $15,391$  & --          \\ \cline{2-8}
                                                         & ASPIN$^\dagger$                               & $ 6,868$                    & $21,918$                              & $37,882$                           & $1,343,071$ & $36,408$  & $784,603$   \\ \cline{2-8}
                                                         & MSPIN$^+$                                     & $ 5,978$                    & $22,367$                              & $48,558$                           & $1,710,703$ & $27,193$  & $520,496$   \\ \hline
  \end{tabular}
\end{table}

In this section, we study the convergence and the performance of the AM and SPIN methods for all benchmark problems.
\Cref{fig:all_solvers_perf} depicts the number of iterations required for reaching the desirable stopping criterion at each loading step.
In particular, we plot a number of global nonlinear iterations, the total number of linear iterations required for solving linear systems associated with displacement and phase-field subproblems, at each loading step.
Furthermore, \Cref{tab:overall_results} reports an accumulated number of nonlinear and linear iterations required for solving arising coupled problems, as well as the displacement and phase-field subproblems.
Here, we point out that the AM method requires a solution of linear systems only during the local (subproblem) phase.
In contrast, the SPIN methods require a solution of the preconditioned coupled linear system at each global nonlinear iteration.
This is achieved using Krylov methods, which rely on an efficient evaluation of the matrix-vector product $(\matP \matJ) \vv$, which in turn necessitates approximation of $\matJ_{uu}^{-1}$ and $\matJ_{cc}^{-1}$, recall \Cref{alg:action_ASPIN} and \Cref{alg:action_MSPIN}.
As specified in \cref{sec:alg_descr_SPIN}, we approximate both inverses by solving associated linear systems using the BCGSTAB method preconditioned with the AMG method.
Therefore, the reported number of linear iterations for displacement and phase-field subproblems also includes the iterations required while solving the preconditioned coupled system.

As we can see from \Cref{fig:all_solvers_perf}, the number of linear and nonlinear iterations increases rapidly as soon as the crack propagation starts.
For the benchmark problems, which exhibit only the ``brutal" crack-propagation, i.e.,~tension and asymmetrically notched beam, a rapid increase in iteration count typically occurs only once per simulation, at the critical time $t_c$.
By taking a closer look at the number of required global nonlinear iterations, we can see that all variants of the AM method except the AM-ST method exhibit almost identical convergence behavior.
This is due to the fact that the AM-ST method satisfies its stopping criterion earlier than the other three AM methods.
Thus, the choice of linear solvers does not impact the overall convergence speed of the AM method, only its performance in terms of the execution time, as we see in \cref{sec:exec_time}.
We also point out that the SPIN methods require a lower number of nonlinear global iterations than all variants of AM method for almost all benchmark problems.
The difference is more prevalent for the problems with the gradual crack propagation (Shear, Three-point bending, L-shaped panel), as the obtained savings accumulate over several loading steps.
For instance, for the shear test, an accumulated number of global nonlinear iterations are around $40,000$ and $82,000$ for AM-ST and AM-INK/NK/ND methods, respectively.
In contrast, the SPIN methods require only around~$20,000$ global nonlinear iterations, see \Cref{tab:overall_results}.

The displacement and phase-field subproblems are solved once per each global nonlinear iteration.
As a consequence, the accumulated number of nonlinear iterations required to solve both subproblems is significantly higher for AM methods than for SPIN solution strategies.
For instance, for the L-shaped panel benchmark problem, the number of accumulated nonlinear iterations required for solving all arising displacement subproblems is approximately $170,000$ for AM-ND/NK/INK methods, $80,000$ for the AM-ST methods while SPIN methods require approximately $70,000$ iterations.
If we examine the iteration count more closely, we also notice that the benchmark problems with gradual crack propagation require fewer nonlinear subproblem iterations per global step than benchmark problems that exhibit brutal crack propagation.
For example, the AM-ND method requires on average approximately $1.11$ nonlinear iterations for displacement subproblem per global step for the shear test, while approximately $2.16$ iterations are needed for the tension test.

SPIN methods require a solution of preconditioned coupled linear system~\eqref{eq:global_linear_new} at each global nonlinear iteration.
Results reported in \Cref{tab:overall_results} suggest that solving these preconditioned systems is fairly cheap as the Krylov methods require only a few iterations to converge.
In particular, for the tension test, the GMRES method requires on average approximately $2.07$ and $2.92$ iterations to solve the preconditioned coupled linear system, when invoked by the ASPIN method and the MSPIN method, respectively.

In the end, we note that the number of linear iterations for both, displacement and phase-field, subproblems mimics the trend observed for the global nonlinear iterations.
Thus, AM-NK and AM-INK methods require a substantially higher amount of linear iterations than SPIN methods.
For instance, for the three-point bending test, AM-NK requires approximately $1.55$ times more linear phase-field subproblem iterations and $2.09$ times more linear displacement subproblem iterations than the MSPIN method.
Here, we also point out that the AM-INK method requires fewer iterations than the AM-NK method, due to the use of the adaptive stopping criterion induced by the inexact Newton's method, recall \cref{sec:INK}

\subsubsection{Analysis of convergence history}
\label{sec:res_iters}
\begin{figure}[t]
  \begin{subfigure}[b]{1\textwidth}
    \begin{minipage}{0.3\linewidth}\centering
      \includegraphics{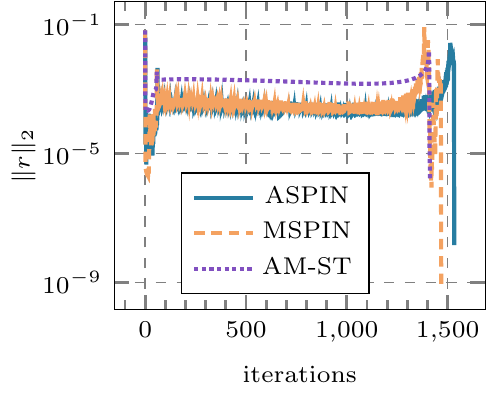}

\end{minipage}\hfill \begin{minipage}{0.3\linewidth}\centering
      \includegraphics{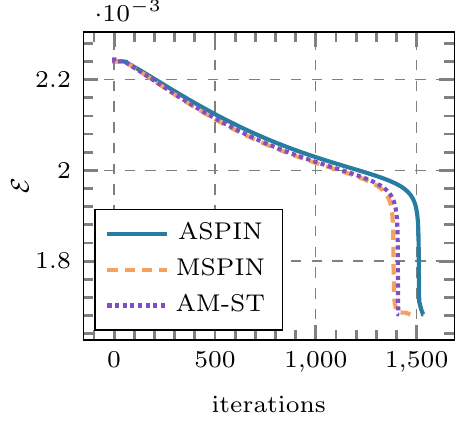}
\end{minipage}\hfill \begin{minipage}{0.33\linewidth}\centering
      \includegraphics[scale=0.0565]{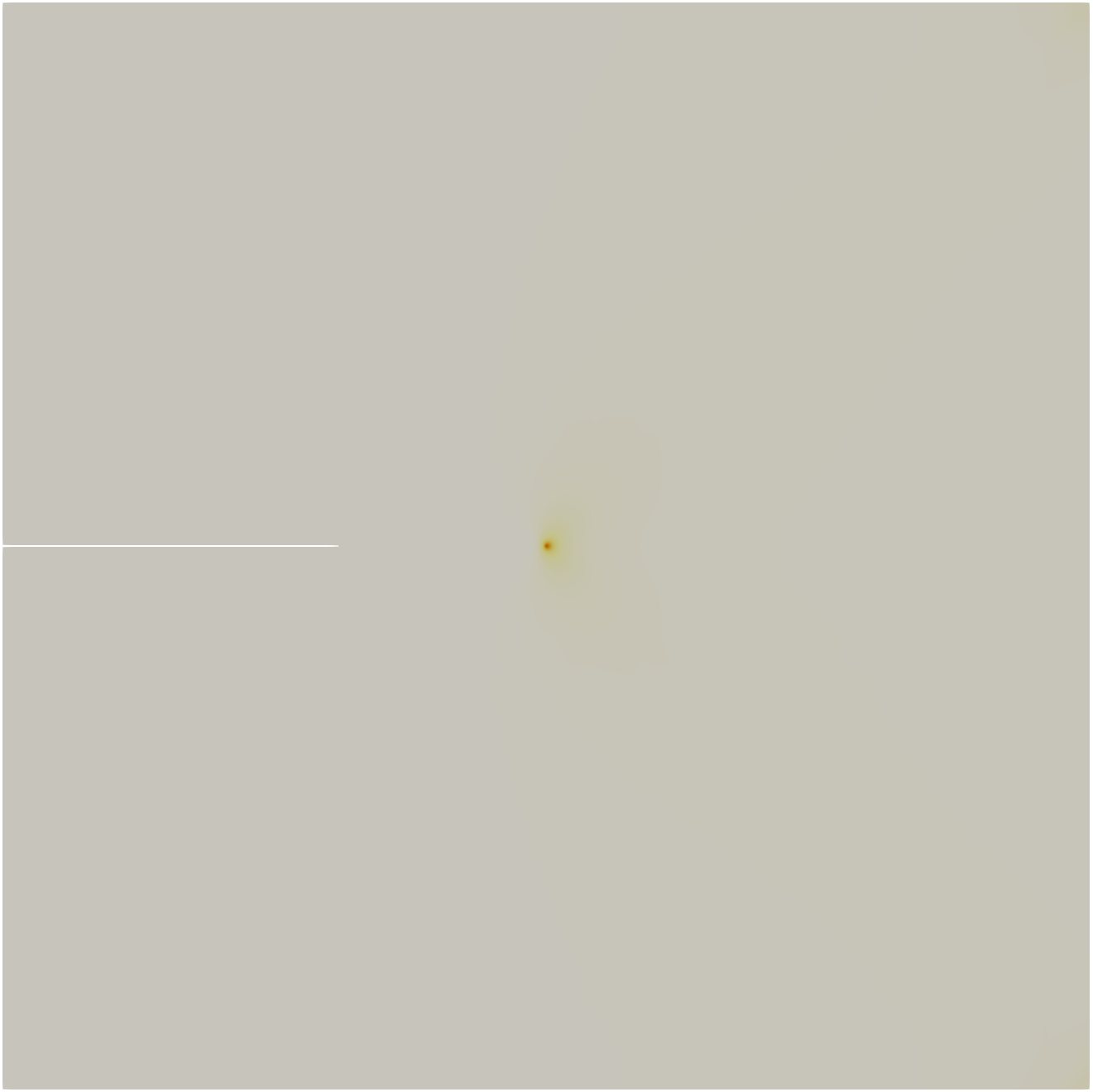}\hspace{3pt}\includegraphics[scale=0.0565]{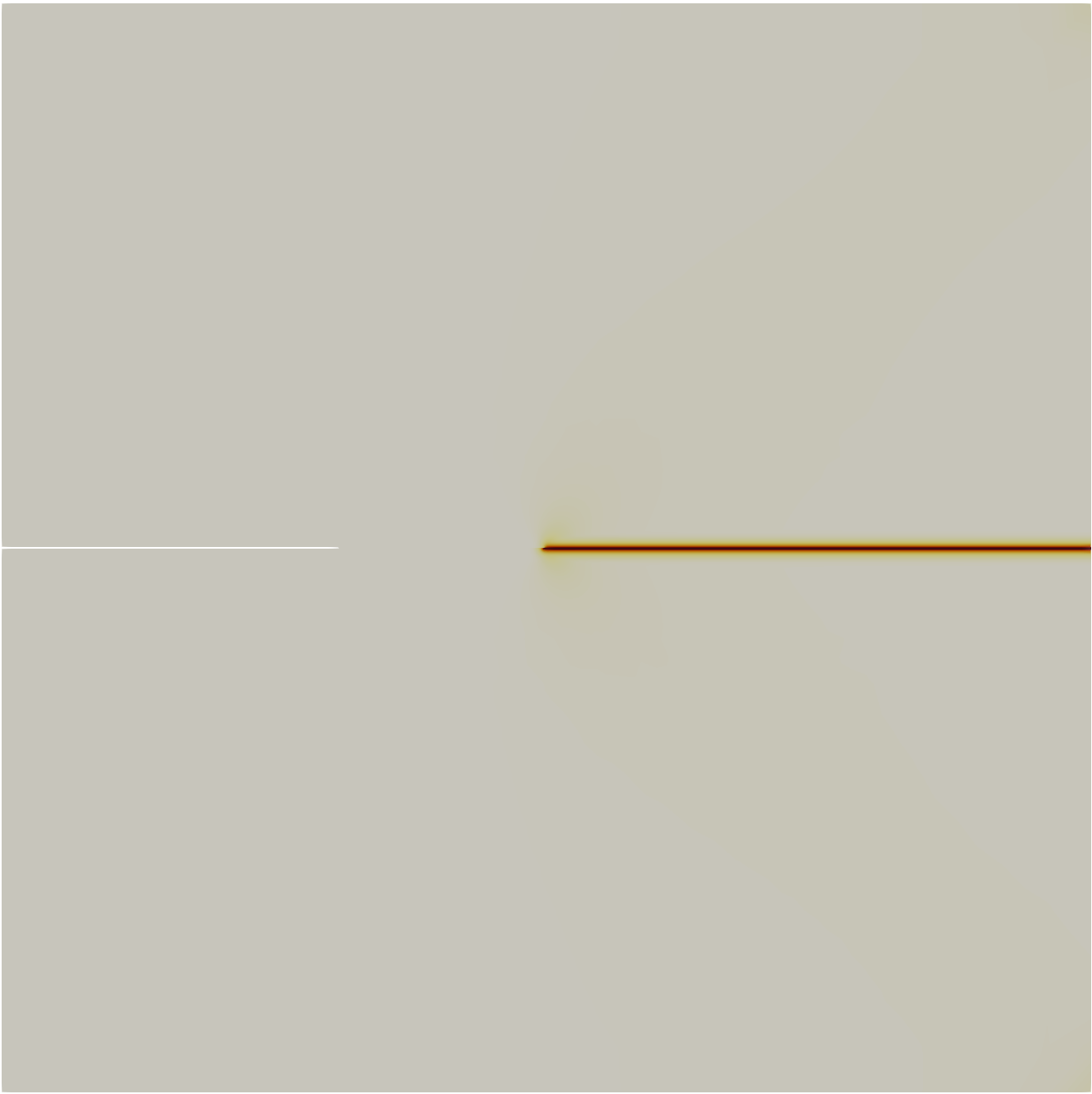}
       \vspace{1pt}
    \end{minipage}\hfill \begin{minipage}{0.05\linewidth}\centering
      \includegraphics[scale=0.022]{c_scale.png}
    \end{minipage}
    \caption{Convergence history in terms of residual norm and value of energy for tension test at loading step 113.}
    \label{fig:residual_energy_01}
  \end{subfigure}
  \begin{subfigure}[b]{1\textwidth}
    \begin{minipage}{0.3\linewidth}\centering
      \includegraphics{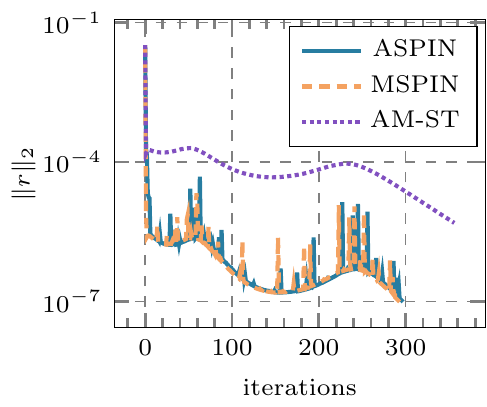}

\end{minipage}\hfill \begin{minipage}{0.3\linewidth}\centering
      \includegraphics{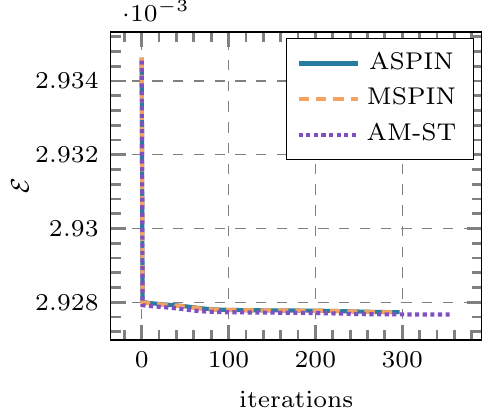}
\end{minipage}\hfill \begin{minipage}{0.33\linewidth}\centering
    \hspace{5pt}
      \includegraphics[scale=0.065]{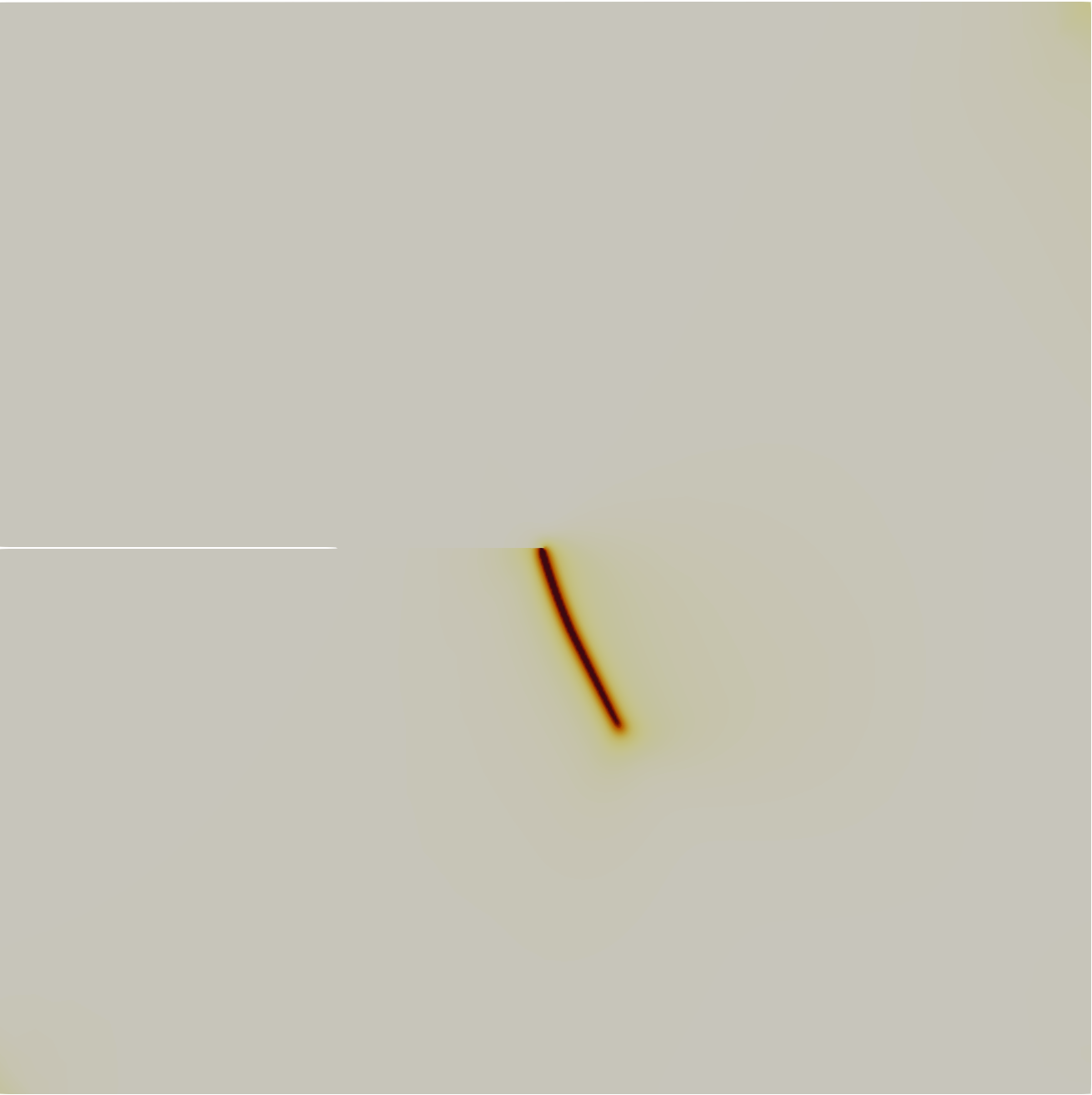}\vspace{1pt}
    \end{minipage}\hfill \begin{minipage}{0.05\linewidth}\centering
      \includegraphics[scale=0.022]{c_scale.png}
    \end{minipage}
    \caption{Convergence history in terms of residual norm and value of energy for shear test at loading step 43.}
    \label{fig:residual_energy_02}
  \end{subfigure}
  \begin{subfigure}[b]{1\textwidth}
    \begin{minipage}{0.3\linewidth}\centering
      \includegraphics{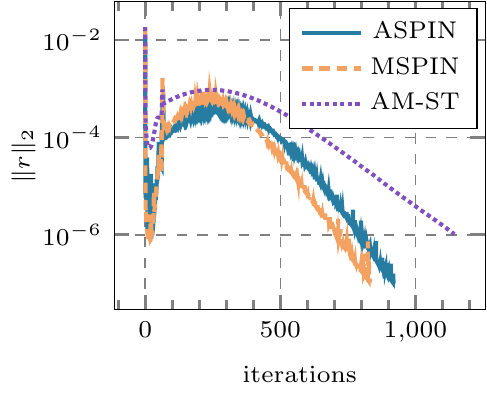}

\end{minipage}
    \hfill
    \begin{minipage}{0.3\linewidth}\centering
      \includegraphics{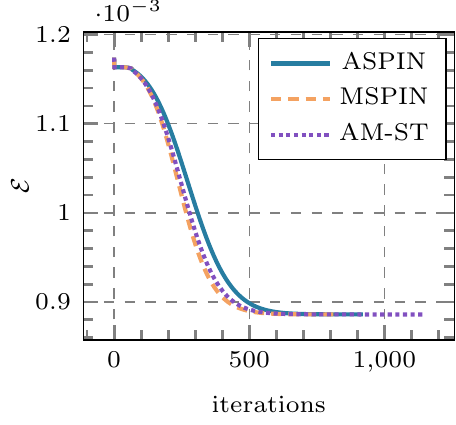}

\end{minipage}\hfill
    \begin{minipage}{0.33\linewidth}\centering
\includegraphics[scale=0.085]{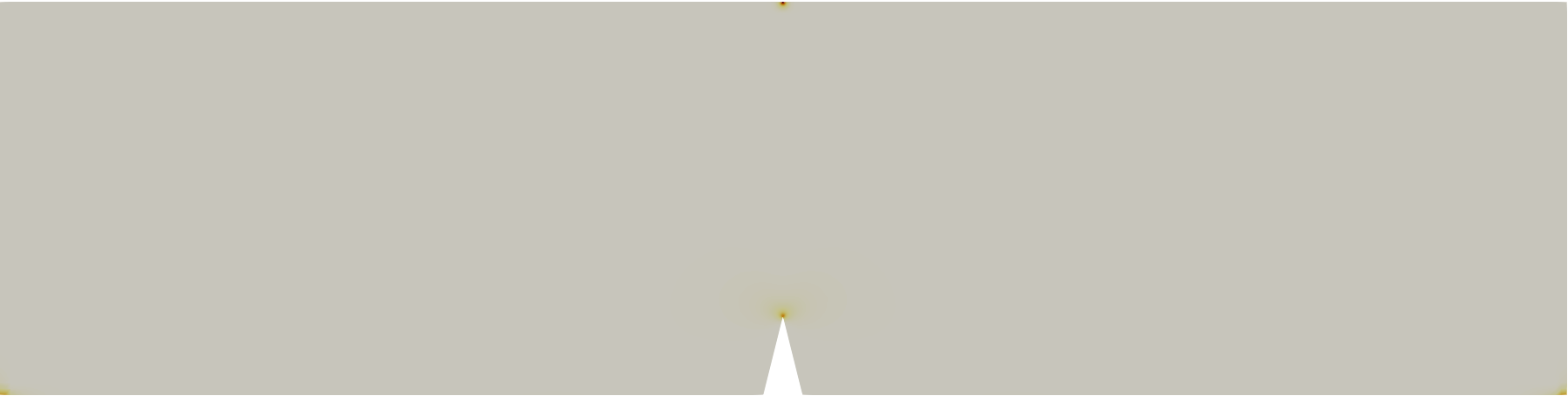}
      \vspace{5pt}
      \includegraphics[scale=0.085]{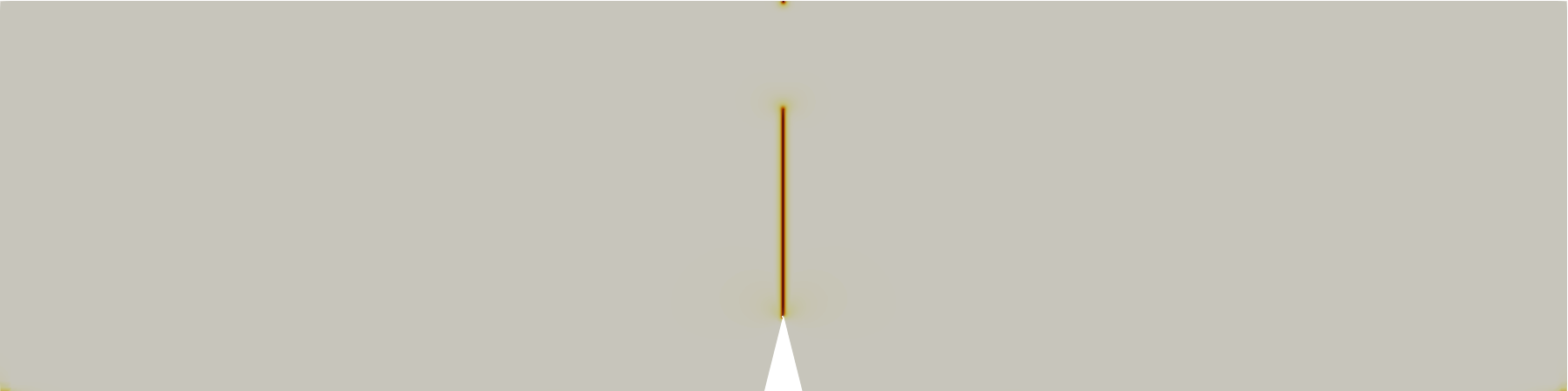}
       \vspace{1pt}
    \end{minipage}
    \hfill
    \begin{minipage}{0.05\linewidth}
      \includegraphics[scale=0.022]{c_scale.png}
    \end{minipage}
    \caption{Convergence history in terms of residual norm and value of energy for three point bending test loading step 57.}
    \label{fig:residual_energy_03}
  \end{subfigure}
  \caption{Convergence history and snapshots of crack evolution for  a given the loading step.}
  \label{fig:res_iter_1}
\end{figure}

In this section, we analyze the convergence behavior of the AM-ST and SPIN methods by examining the norm of residual and the value of energy during the solution process.
For the tension test, three-point bending test, and asymmetrically notched beam test, we focus our attention on the loading steps where the ``brutal" crack propagation occurs.
For the shear test, the crack propagation takes place in a gradual manner, so we pick a random loading step.
In the case of the L-shaped panel test, a crack is first propagated brutally, which is followed by gradual crack propagation.
Therefore, we consider the loading step where the crack propagation initiates.

\begin{figure}
  \begin{subfigure}[b]{1\textwidth}
    \begin{minipage}{0.3\linewidth}\centering
      \includegraphics{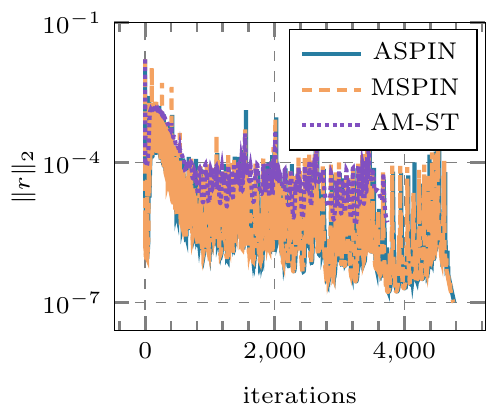}

\end{minipage}\hfill \begin{minipage}{0.3\linewidth}\centering
      \includegraphics{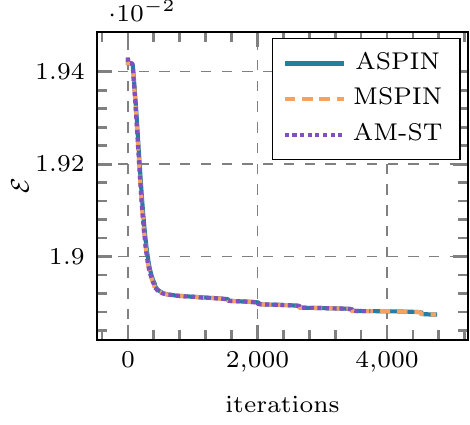}
\end{minipage}\hfill \begin{minipage}{0.33\linewidth}\centering
      \includegraphics[scale=0.0565]{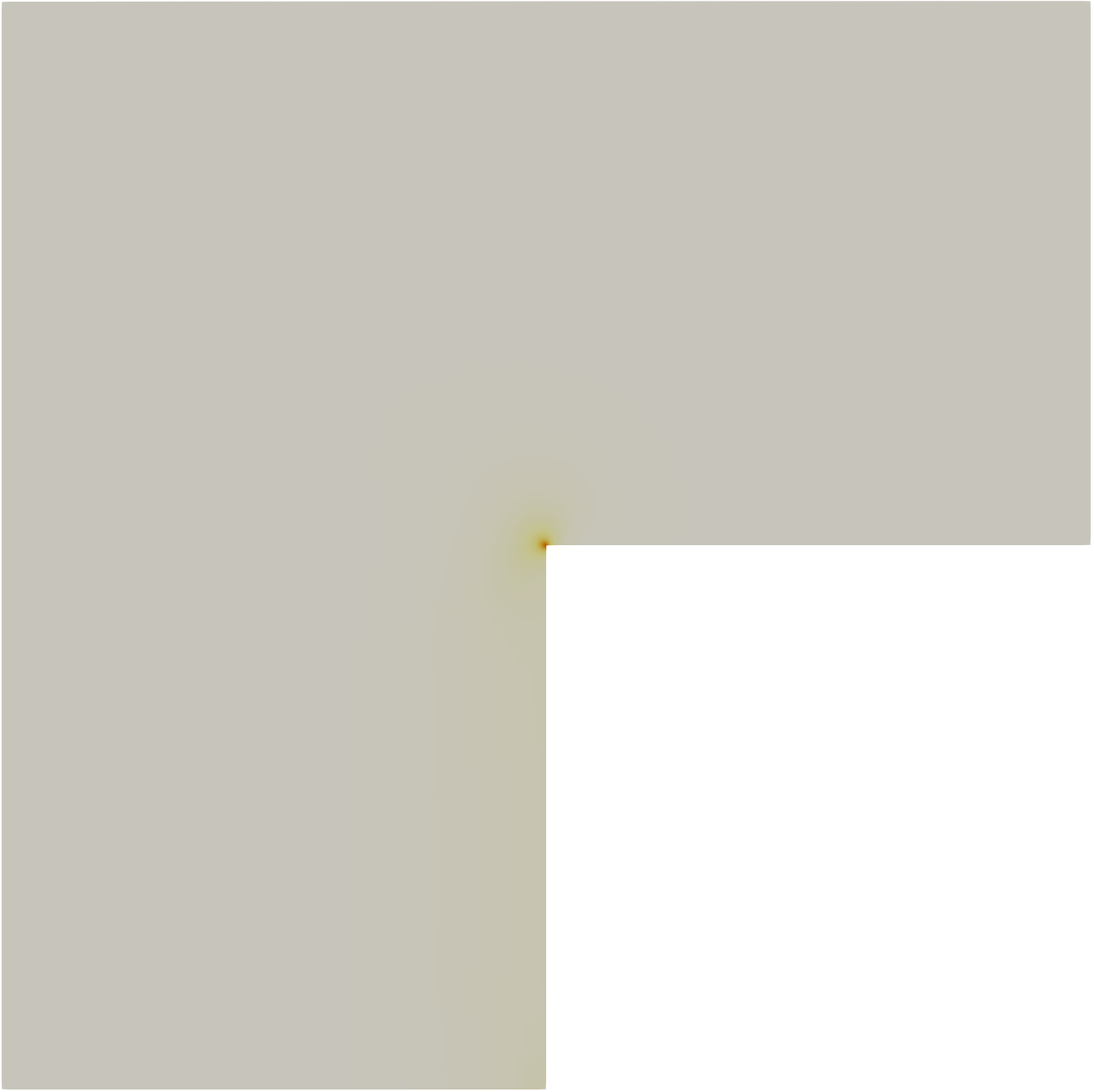}\hspace{3pt}\includegraphics[scale=0.0565]{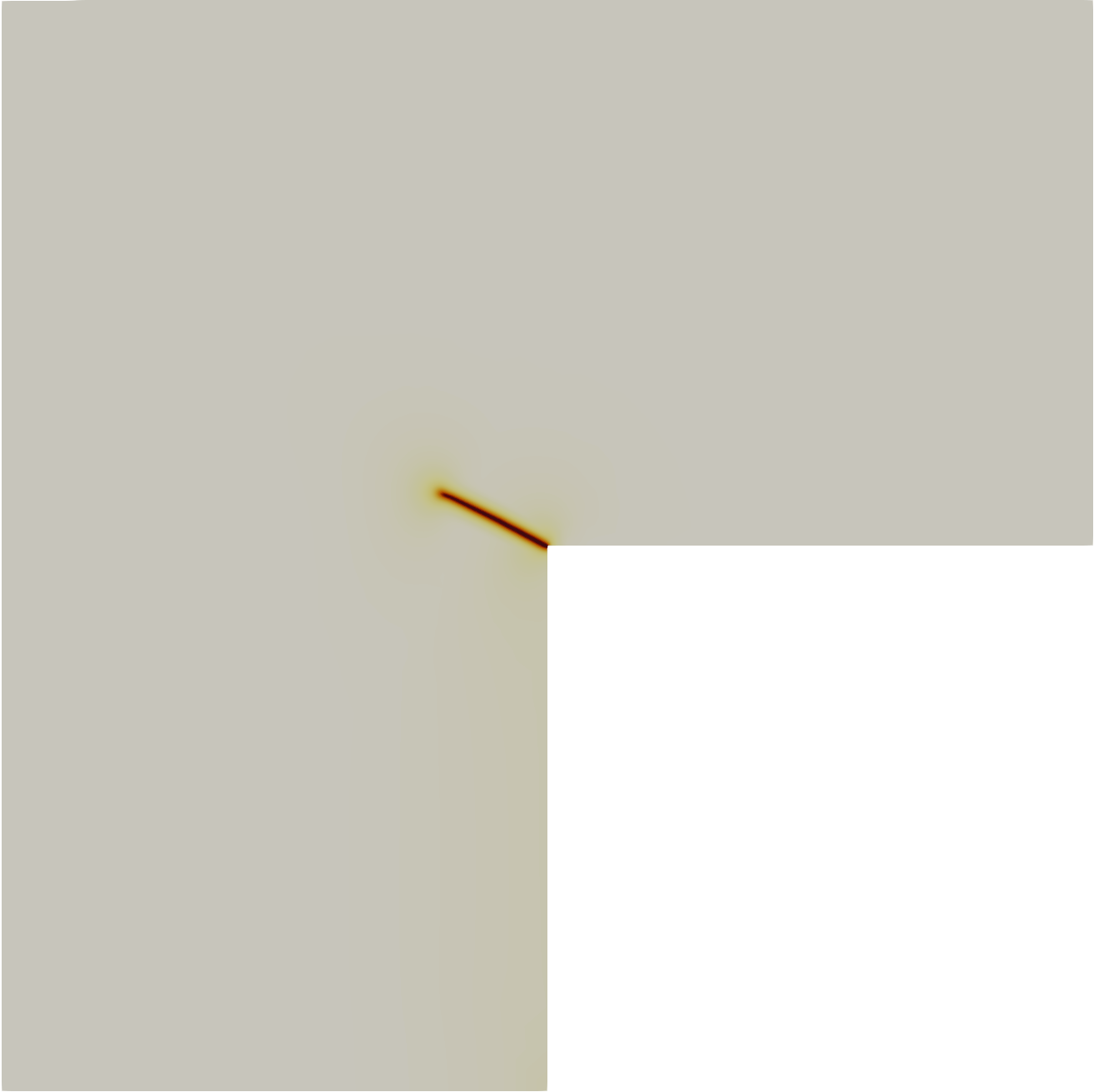}
      \vspace{5pt}
    \end{minipage}\hfill \begin{minipage}{0.05\linewidth}\centering
      \includegraphics[scale=0.022]{c_scale.png}
    \end{minipage}
\caption{Convergence history in terms of residual norm and value of energy for L-shaped panel test at loading step $70$.}
    \label{fig:residual_energy_04}
  \end{subfigure}
  \begin{subfigure}[b]{1\textwidth}
    \begin{minipage}{0.3\linewidth}\centering
      \includegraphics{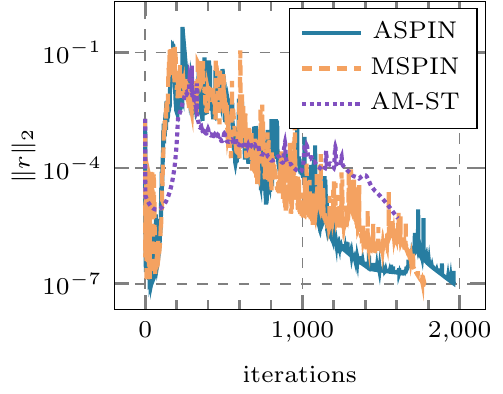}

\end{minipage}\hfill \begin{minipage}{0.3\linewidth}\centering
      \includegraphics{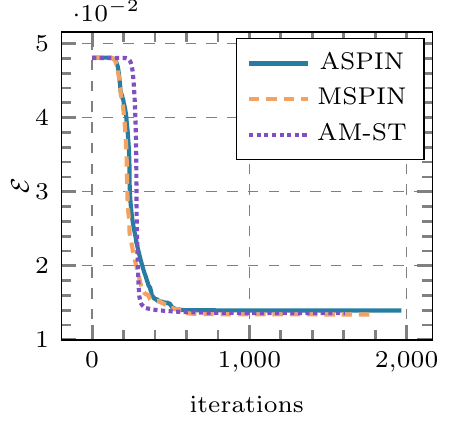}
\end{minipage}\hfill
    \begin{minipage}{0.33\linewidth}\centering
\includegraphics[scale=0.045]{05_AMST.png}
      \vspace{1pt}
      \includegraphics[scale=0.045]{05_ASPIN.png}
      \vspace{1pt}
      \includegraphics[scale=0.045]{05_MSPIN.png}
    \end{minipage}
    \hfill
    \begin{minipage}{0.05\linewidth}\centering
      \includegraphics[scale=0.022]{c_scale.png}
    \end{minipage}
    \caption{Convergence history in terms of residual norm and value of energy for asymmetrically notched beam test. We observe three distinct crack patterns for AM-ST, ASPIN and MSPIN methods (from top to bottom).}
    \label{fig:residual_energy_05}
  \end{subfigure}
  \caption{Convergence history and snapshots of crack evolution for  a given the loading step.}
  \label{fig:res_iter_22}
\end{figure}

\Cref{fig:res_iter_1} illustrates the norm of residual and energy over iterations as well as the crack pattern before and after the specific time step for the tension test, the shear test, and the three-point bending tests.
\Cref{fig:res_iter_22} depicts the same quantities for the L-shaped panel test and the asymmetrically notched beam tests.
From Figures \ref{fig:res_iter_1} and \ref{fig:res_iter_22}, we observe that the AM-ST method terminates sooner than SPIN methods.
We can notice that the energy is monotonically decreasing as we have employed the backtracking line search method.
In contrast, the norm of residual is allowed to oscillate.
For the considered test cases, we see that the residual is oscillating more for the SPIN methods than for the AM-ST method.

For the tension test (\Cref{fig:residual_energy_01}), we see that the norm of residual reduces rapidly in the first few iterations, but then it increases and oscillates before actually achieving the quadratic convergence after more than $1,000$ iterations.
In the last few iterations, we also observe that the energy reduces dramatically.
This is quite contrary to the shear test (\Cref{fig:residual_energy_02}), where the energy reduces rapidly only in the first iteration.
The three-point bending test also shows a different behavior (\Cref{fig:residual_energy_03}) as the energy is reduced in the first $500$ iterations and then a flat region is encountered.
Although, we notice that the norm of residual is reducing at a higher rate. 

The norm of residual is more oscillatory for the L-shaped panel test than for other benchmark problems, as can be seen from \Cref{fig:residual_energy_04}.
Here, we notice that the AM-ST method requires fewer iterations to converge than the SPIN method as it employs weaker termination criteria.
More interesting results can be observed in~\Cref{fig:residual_energy_05}, as the obtained crack pattern for the asymmetrically notched beam differs for different solution strategies.
For this test case, we can see that the MSPIN method achieves a solution with the smallest value of energy in comparison with the other methods.

\subsection{Execution time}
\label{sec:exec_time}
\begin{table}[ht!]
  \caption{The execution time required to simulate all benchmark problems using ASPIN, MSPIN, AM-ST, AM-ND, AM-NK, and AM-INK methods.
  The symbols $*,\dagger,+$ denote convergence to a different solution.}
  \label{tab:computational_time}
  \centering
  \begin{tabular}{|c|l|r|ccccc|}
    \hline
    \multirow{2}{*}{Example} & \multicolumn{1}{c|}{\multirow{2}{*}{Solver}} & \multirow{2}{*}{Time (min)} & \multicolumn{5}{c|}{Speedup with respect to}                                     \\ \cline{4-8}
                             &                                              &                             & AM-ND                                        & AM-NK  & AM-INK & AM-ST  & ASPIN  \\ \hline \hline
    \multirow{6}{*}{Tension} & AM-ND                                        & $247.72$                    & --                                           & --     & --     & --     & --     \\ \cline{2-8}
                             & AM-NK                                        & $168.92$                    & $1.47$                                       & --     & --     & --     & --     \\ \cline{2-8}
                             & AM-INK                                       & $157.02$                    & $1.58$                                       & $1.08$ & --     & --     & --     \\ \cline{2-8}
                             & AM-ST                                        & $254.10$                    & $0.97$                                       & $0.66$ & $0.62$ & --     & --     \\ \cline{2-8}
                             & ASPIN                                        & $183.07$                    & $1.35$                                       & $0.92$ & $0.86$ & $1.39$ & --     \\ \cline{2-8}
                             & MSPIN                                        & $131.55$                    & $1.88$                                       & $1.28$ & $1.19$ & $1.93$ & $1.39$ \\ \hline \hline

    \multirow{6}{*}{{Shear}} & AM-ND                                        & $2526.34$                   & --                                           & --     & --     & --     & --     \\ \cline{2-8}
                             & AM-NK                                        & $1,716.82$                   & $1.47$                                       & --     & --     & --     & --     \\ \cline{2-8}
                             & AM-INK                                       & $1,757.77$                   & $1.44$                                       & $0.98$ & --     & --     & --     \\ \cline{2-8}
                             & AM-ST                                        & $1,318.44$                   & $1.92$                                       & $1.30$ & $1.33$ & --     & --     \\ \cline{2-8}
                             & ASPIN                                        & $498.36 $                   & $5.07$                                       & $3.44$ & $3.53$ & $2.65$ & --     \\ \cline{2-8}
                             & MSPIN                                        & $453.78 $                   & $5.57$                                       & $3.75$ & $3.87$ & $2.91$ & $1.10$ \\ \hline \hline

    \multirow{6}{2.5cm}{\centering Three-point bending}
                             & AM-ND                                        & $1,378.35$                   & --                                           & --     & --     & --     & --     \\ \cline{2-8}
                             & AM-NK                                        & $988.34 $                   & $1.39$                                       & --     & --     & --     & --     \\ \cline{2-8}
                             & AM-INK                                       & $855.32 $                   & $1.61$                                       & $1.16$ & --     & --     & --     \\ \cline{2-8}
                             & AM-ST                                        & $810.00 $                   & $1.70$                                       & $1.22$ & $1.06$ & --     & --     \\ \cline{2-8}
                             & ASPIN                                        & $291.01 $                   & $4.74$                                       & $3.40$ & $2.94$ & $2.78$ & --     \\ \cline{2-8}
                             & MSPIN                                        & $260.63 $                   & $5.29$                                       & $3.79$ & $3.28$ & $3.11$ & $1.12$ \\ \hline \hline

    \multirow{6}{*}{{L-shaped panel}}
                             & AM-ND                                        & $8,662.82$                   & --                                           & --     & --     & --     & --     \\ \cline{2-8}
                             & AM-NK                                        & $5,748.19$                   & $1.51$                                       & --     & --     & --     & --     \\ \cline{2-8}
                             & AM-INK                                       & $5,862.43$                   & $1.48$                       & $0.98$   & --     & --     & --     \\ \cline{2-8}
                             & AM-ST                                        & $3,539.46$                   & $2.45$                       & $1.62$ & $1.66$    & --     & --     \\ \cline{2-8}
                             & ASPIN                                        & $2,674.32$                   & $3.24$                       & $2.15$ & $2.19$    & $1.32$ & --     \\ \cline{2-8}
                             & MSPIN                                        & $2,564.15$                   & $3.38$                      & $2.24$ & $2.29$    & $1.38$ & $1.04$ \\ \hline  \hline

    \multirow{6}{2.5cm}{\centering {Asymmetrically notched beam}}
                             & AM-ND$^\ast$                                 & $1,263.54$                   & --                                           & --     & --     & --     & --     \\ \cline{2-8}
                             & AM-NK$^\ast$                                 & $1,060.75$                   & $1.19$                                       & --     & --     & --     & --     \\ \cline{2-8}
                             & AM-INK$^\ast$                                & $ 952.49$                   & $1.33$                                       & $1.11$ & --     & --     & --     \\ \cline{2-8}
                             & AM-ST$^\ast$                                 & $ 563.77$                   & $2.24$                                       & $1.88$ & $1.69$ & --     & --     \\ \cline{2-8}
                             & ASPIN$^\dagger$                              & $ 667.23$                   & $1.89$                                       & $1.59$ & $1.43$ & $0.84$ & --     \\ \cline{2-8}
                             & MSPIN$^+$                                    & $ 710.64$                   & $1.78$                                       & $1.49$ & $1.34$ & $0.79$ & $0.94$ \\ \hline
  \end{tabular}
\end{table}

In this section, we compare the execution time required by the AM and SPIN methods for all benchmark problems.
The reported execution times stand for the total execution time necessary to run the entire simulation, which includes, call to nonlinear assembly routines, setup and execution of nonlinear and linear solvers, I/O, and post-processing.
The obtained results, summarized in \Cref{tab:computational_time}, demonstrate that the SPIN methods are faster than the AM methods for all benchmark problems.
One exception in these tests is the asymmetrically notched beam test.
Since the AM, the ASPIN, and the MSPIN methods converge to different solutions for this particular example, we can not compare the iteration numbers and the computation time, in a fair manner.
The obtained speedup varies across different examples, depending on the type of crack propagation, particular loading, and the problem size.
As expected, a larger speedup is obtained for problems with the gradual crack propagation.
As a consequence, we observe the speedup up to a factor of $1.8$ for the problems with ``brutal" crack propagation, and up to a factor of $5.57$ for the problems with gradual crack propagation.
Moreover, as shown in \cref{sec:mesh_res}, the speedup of SPIN methods increases with the problem size.

The obtained numerical results also suggest that the alternate minimization configured with Krylov methods (AM-NK, AM-INK) is always faster than the configuration with the direct solver (AM-ND).
This is not surprising, as the LU factorization of the sparse linear system with $n$ unknowns requires~$\pazocal{O}(n^{3/2})$ flops in $2D$, recall \cref{sec:newton_method}.
In contrast, the multigrid preconditioner requires approximately $\pazocal{O}(n)$ flops and it is also optimal in terms of memory usage.
As a consequence, the AM-ND and AM-ST methods can be very effective for small problems, but their use becomes prohibitive for larger problems.
We also observe that the AM-INK is slightly faster than AM-NK, in particular by a factor of $1.1$ on average.
The execution time of the ASPIN and the MSPIN methods is comparable, favoring MSPIN by a factor of $1.13$ on average.
We note that the behavior observed in terms of execution time corresponds to the results reported in terms of iteration count in \Cref{tab:overall_results}.

\subsubsection{Analysis of a nonlinear iteration cost}
In this section, we analyze the computational time of a single iteration of the AM and the SPIN methods.
We consider a three-point bending test with two different mesh resolutions, corresponding to coupled problems with $99,135$ and $1,572,147$ dofs.
We report the time required for solving the displacement subproblems (blue color), the phase-field subproblems (brown color), and the global/coupled linear problem (red color).
Furthermore, we also distinguish between the time required by the linear solvers, and the assembly routines, which include an evaluation of the energy, gradients, and Hessians.
The results presented in this section are based on serial executions.

\begin{figure}[t]
  \includegraphics{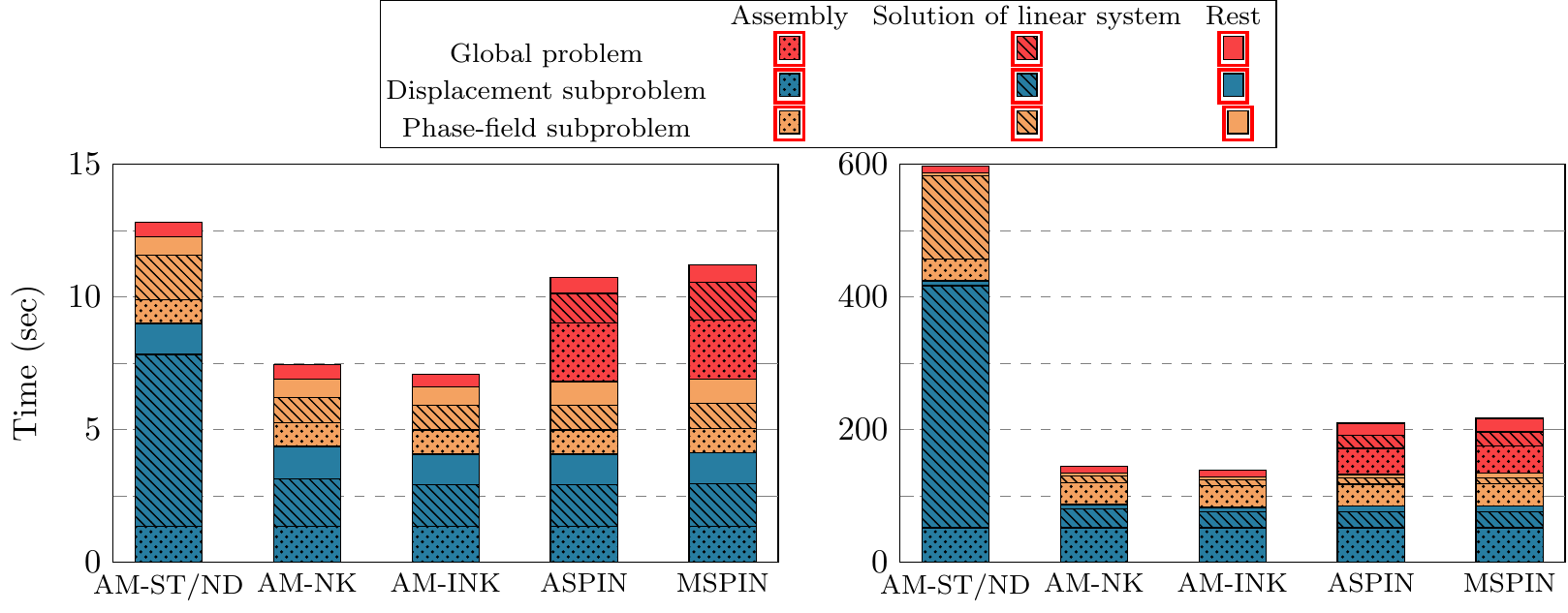}

\caption{Analysis of execution time for a nonlinear iteration of alternate minimization and SPIN methods.
    Time spent solving the displacement subproblem is denoted by blue, the phase-field subproblem by brown, and the global problem by red.
    The time spent by the assembly of energy, gradient, and Hessian is depicted by the dotted pattern, while time attributed to the solution of a linear system is depicted by the slanted line pattern. No pattern denotes the time required by all other routines.
    The experiment was performed in serial using a three-point bending test with $99,135$ dofs (left) and $1,572,147$ dofs (right).}
  \label{fig:timings}
\end{figure}

As we can observe from \Cref{fig:timings}, the execution time of a single nonlinear iteration for SPIN methods is higher than for the AM-NK/INK methods.
This is due to the fact that SPIN methods require the solution of a nonlinearly preconditioned linear system of equations, as well as the assembly of coupled gradient and Hessian.
Moreover, SPIN methods also employ a line-search backtracking strategy to ensure that the energy of the coupled problem is decreasing monotonically.
In contrast, AM methods only require the solution of displacement and phase-field subproblems, while the time attributed to the global problem is only associated with the transfer of data between different subproblems and convergence monitoring.
However, if a direct solver is used then an iteration of AM method is more time-consuming than an iteration of the SPIN method.
Moreover, this difference grows with an increasing number of dofs.
In the end, we also point out that all solution strategies which employ iterative linear solvers spent a larger amount of computational time for assembly than for solutions of linear systems.
This is due to the fact that in order to evaluate energy functional and its derivatives on each quadrature point, we have to perform the eigen-decomposition of the strain tensor.

\subsection{Memory requirements}
\label{sec:memory}
\begin{figure}[t]
  \includegraphics{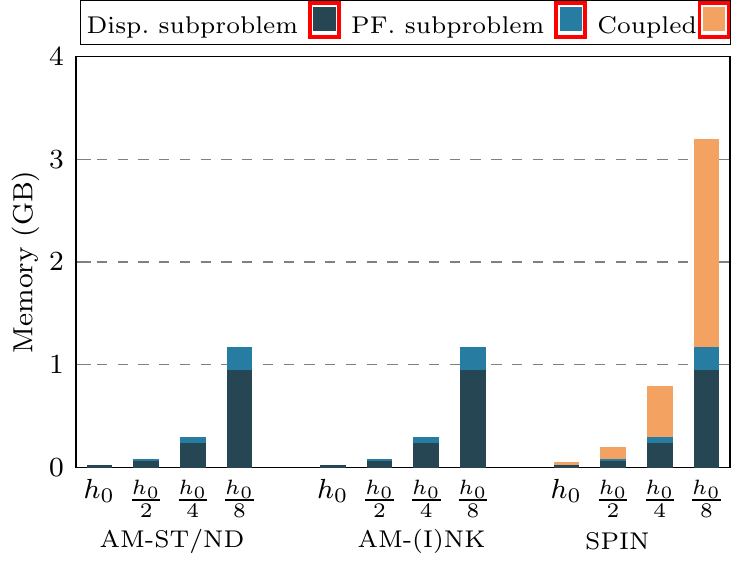}  \hfill   \includegraphics{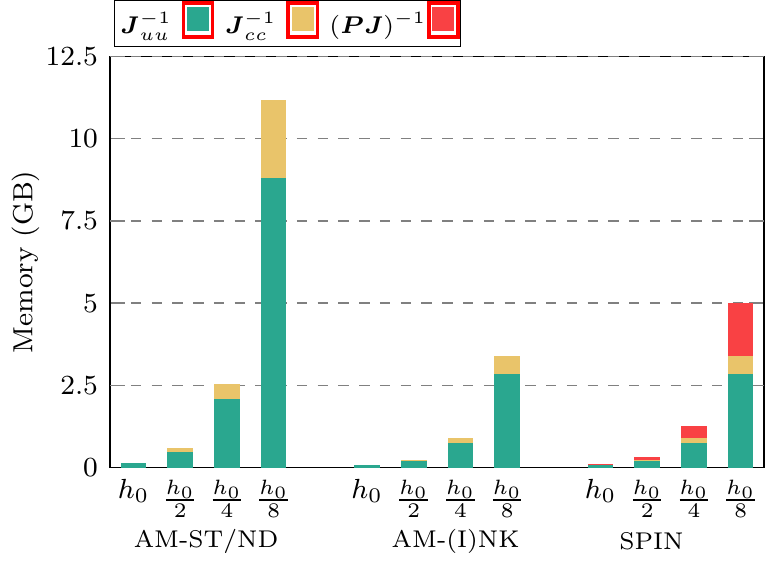}

\includegraphics{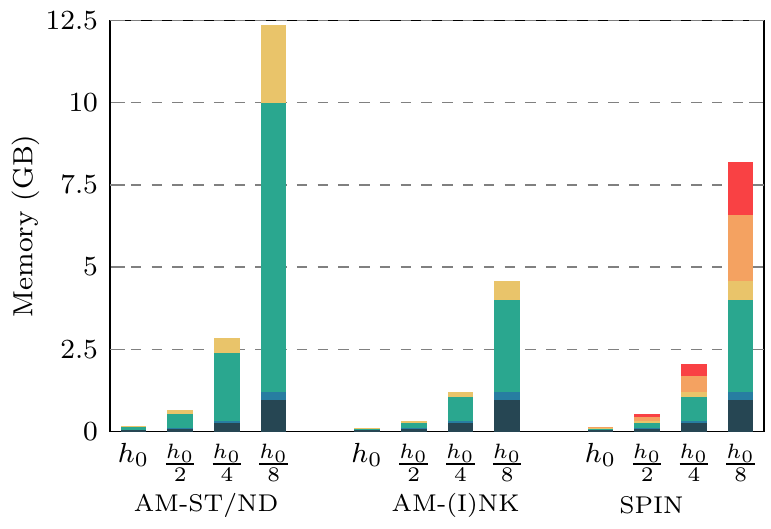}  \hfill  \includegraphics{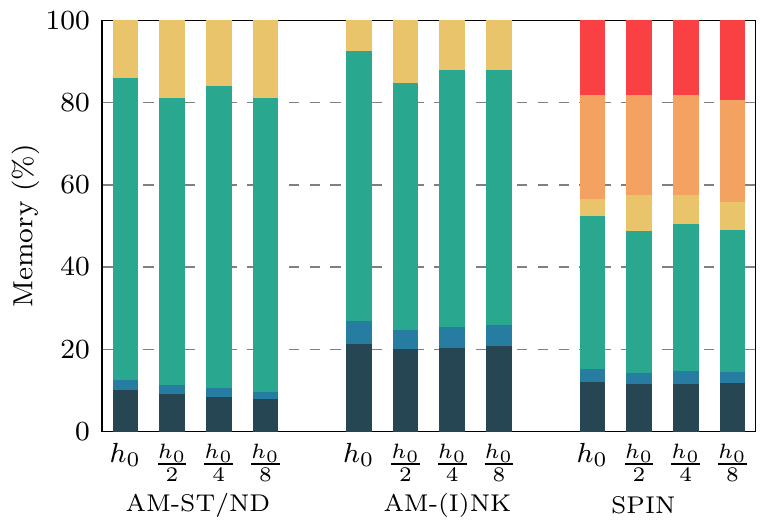}

\caption{Memory requirements of different variants of alternate minimization and SPIN methods.
    Top left: Memory requirement associated with storing Hessians for phase-field/displacement subproblems and the fully coupled problem.
    Top right: Memory requirement associated with solution process of each arising linear systems
    Bottom left: Combined memory requirement for operators and the solution process.
    Bottom right: The percentage weight of the memory requirement associated with overall system.
  }
  \label{fig:memory_requirements}
\end{figure}

In this section, we investigate the memory requirements of the employed solution strategies.
The analysis is performed using a three-point bending test with increasing mesh refinements.
During the presented study, we distinguish between the memory requirements associated with the assembly of quantities, such as gradients, and Hessians, and the memory requirements associated with the solution of arising linear systems, see Figure~\ref{fig:memory_requirements} on the top left and top right, respectively.
As we can observe, the memory requirements associated with assembly routines are higher for SPIN methods, as they require the knowledge of coupled gradient and Hessian.
In contrast, AM methods require to store only quantities associated with the phase-field and displacement subproblems.
Regarding the memory requirements associated with the solution of linear systems, it is important to distinguish if a direct or iterative linear solver is employed.
As we can see, the memory requirements of AM-ST/ND methods grow significantly with an increasing number of dofs.
This is due to the fact that the LU factorization of the Jacobians contains a large amount of non-zero elements, giving rise to a so-called memory fill-in.
We note that a number of non-zeros grow with $\pazocal{O}(n^{2/3})$ in two dimensions and with $\pazocal{O}(n^{2})$ in three dimensions.
In contrast, the memory requirements of Krylov methods preconditioned with the multigrid method grow linearly with the problem size.
As a consequence, the overall memory requirements of AM-ST/ND are significantly higher for problems with a large number of dofs, see Figure~\ref{fig:memory_requirements} on the left bottom.

\subsection{Efficiency of the solution strategies with increasing mesh resolution}
\label{sec:mesh_res}
\begin{figure}
  \includegraphics{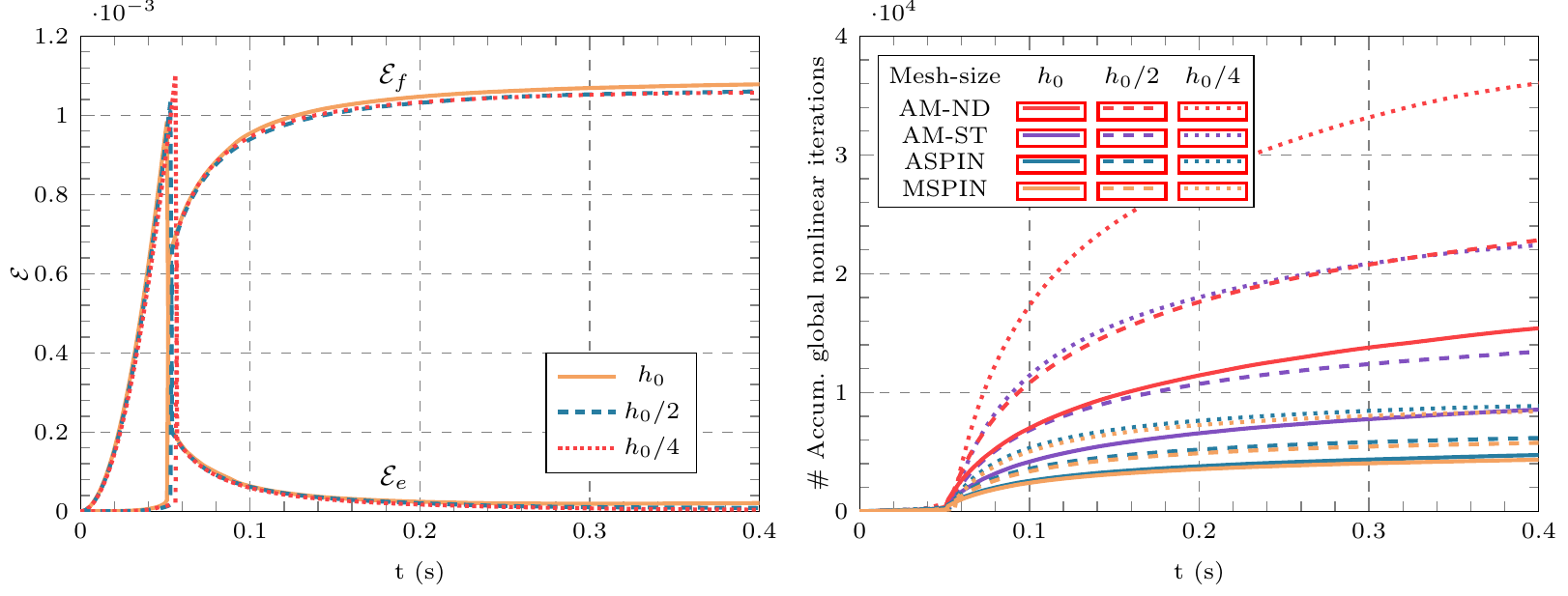}

\caption{Three-point bending test.
    Left: The evolution of elastic and fracture energy over displacement with respect to increasing refinement level.
    The experiment is performed using the MSPIN method.
    Right: A number of accumulated global nonlinear iterations as a function of refinement level for AM-ND, AM-ST, ASPIN and MSPIN methods.
    The symbol $h_0$ denotes the mesh size of the initial mesh.}
  \label{fig:3pt_bending_conv_ref}
\end{figure}
\begin{table}[ht!]
  \caption{The execution time for the three-point bending test with respect to varying mesh resolution.}
  \label{tab:computational_time_wrt_ref}
  \centering
  \begin{tabular}{|c|l|r|rrrrr|}
    \hline
    \multirow{2}{*}{Mesh size ($h$)} & \multicolumn{1}{c|}{\multirow{2}{*}{Solver}} & \multirow{2}{*}{Time (min)} & \multicolumn{5}{c|}{Speedup with respect to}                                     \\ \cline{4-8}
                                     &                                              &                             & AM-ND                                        & AM-NK  & AM-INK & AM-ST  & ASPIN  \\ \hline \hline
    \multirow{6}{2.5cm}{\centering $h_0$}
                                     & AM-ND                                        & $34.65$                     & --                                           & --     & --     & --     & --     \\ \cline{2-8}
                                     & AM-NK                                        & $34.77$                     & $1.00$                                       & --     & --     & --     & --     \\ \cline{2-8}
                                     & AM-INK                                       & $29.88$                     & $1.16$                                       & $1.16$ & --     & --     & --     \\ \cline{2-8}
                                     & AM-ST                                        & $19.15$                     & $1.81$                                       & $1.82$ & $1.56$ & --     & --     \\ \cline{2-8}
                                     & ASPIN                                        & $17.34$                     & $2.00$                                       & $2.01$ & $1.72$ & $1.10$ & --     \\ \cline{2-8}
                                     & MSPIN                                        & $16.21$                     & $2.14$                                       & $2.15$ & $1.84$ & $1.18$ & $1.07$ \\ \hline \hline
    \multirow{6}{2.5cm}{\centering $\dfrac{h_0}{2}$}
                                     & AM-ND                                        & $192.35$                    & --                                           & --     & --     & --     & --     \\ \cline{2-8}
                                     & AM-NK                                        & $167.22$                    & $1.15$                                       & --     & --     & --     & --     \\ \cline{2-8}
                                     & AM-INK                                       & $163.87$                    & $1.17$                                       & $1.02$ & --     & --     & --     \\ \cline{2-8}
                                     & AM-ST                                        & $121.61$                    & $1.58$                                       & $1.38$ & $1.35$ & --     & --     \\ \cline{2-8}
                                     & ASPIN                                        & $61.75$                     & $3.11$                                       & $2.71$ & $2.65$ & $1.97$ & --     \\ \cline{2-8}
                                     & MSPIN                                        & $56.26$                     & $3.42$                                       & $2.97$ & $2.91$ & $2.16$ & $1.10$ \\ \hline \hline
    \multirow{6}{2.5cm}{\centering $\dfrac{h_0}{4}$}
                                     & AM-ND                                        & $1378.35$                   & --                                           & --     & --     & --     & --     \\ \cline{2-8}
                                     & AM-NK                                        & $988.34 $                   & $1.39$                                       & --     & --     & --     & --     \\ \cline{2-8}
                                     & AM-INK                                       & $855.32 $                   & $1.61$                                       & $1.16$ & --     & --     & --     \\ \cline{2-8}
                                     & AM-ST                                        & $810.00 $                   & $1.70$                                       & $1.22$ & $1.06$ & --     & --     \\ \cline{2-8}
                                     & ASPIN                                        & $291.01 $                   & $4.74$                                       & $3.40$ & $2.94$ & $2.78$ & --     \\ \cline{2-8}
                                     & MSPIN                                        & $260.63 $                   & $5.29$                                       & $3.79$ & $3.28$ & $3.11$ & $1.12$ \\ \hline
  \end{tabular}
\end{table}

In this section, we investigate the convergence properties of the solution strategies with respect to increasing mesh resolution.
During this experiment, we use a three-point bending test, described in \cref{sec:3pt_bending_description}.
We consider a hierarchy of three meshes, which were obtained by uniformly refining an initial - adaptively refined finite element mesh.
The problem on the coarsest mesh ($h_0$) consists of the $6,417$ dofs and the problem on the finest level ($h_0/4$) has $99,135$ dofs.
With each refinement step, the mesh size $h$ and the value of the length-scale parameter $l_s$ decrease by a factor of two.
As a consequence, we can approximate the fracture surface more accurately.
However, this comes at a higher computational cost for two main reasons.
Firstly, the non-linearity of the underlying energy minimization problem tends to strengthen, which causes an increase in the number of nonlinear iterations.
Secondly, the computational cost per nonlinear iteration grows due to the increasing number of dofs.
Moreover, the condition number of the arising linear systems grows with each refinement level, which in turn influences the convergence speed of Krylov methods.
For instance, the number of required iterations for the BCGSTAB method grows proportionally with the condition number~\cite{saad2003iterative}.
In this work, we aim to eliminate this dependence on the condition number by employing the AMG preconditioner, which is known to be mesh independent, in the sense, that the number of iterations should stay constant for an increasing number of dofs.
But in practice, this is not easy to achieve as the penalty parameter ($\gamma$) employed to ensure the irreversibility condition also increases with each refinement.
Here, we remark that even though the SPIN method by itself is a nonlinear preconditioner for the global/coupled problem, a scalable solution strategy is required to solve linear problems arising during the minimization of the phase-field and displacement subproblems (recall \cref{sec:spin}).

\Cref{fig:3pt_bending_conv_ref} on the left demonstrates that energy can be approximated more accurately as we employ finer meshes.
Also, \Cref{fig:3pt_bending_conv_ref} on the right illustrates how the number of accumulated global nonlinear iterations increases with increasing mesh resolution.
Interestingly, an increase in the number of iterations is more prominent for the alternate minimization than for the SPIN methods.
This suggests, that using SPIN methods becomes even more beneficial for the larger problems.
Indeed, the results reported in \Cref{tab:computational_time_wrt_ref} demonstrate that the SPIN methods achieve a larger speedup compared to the AM method, as the problem size grows.
For instance, the MSPIN method achieves a speedup of $2.14$ compared to AM-ND, when the mesh with mesh size $h_0$ is used.
In contrast, the speedup of $5.29$ is observed for the problem discretized using the twice refined mesh, i.e., with the mesh size $h_0/4$.
While in comparison with the AM-ST method, the MSPIN method achieves a speedup of $1.18$ for the mesh size $h_0$, but this speedup grows almost linearly with increasing problem size, i.e. with mesh sizes $h_0/2$, $h_0/4$, we observe a speedup of $2.16$ and $3.11$, respectively.

   \section{Conclusion}
\label{sec:conclusion}
In this work, we have proposed a variant of the nonlinear Schwarz Preconditioned Inexact Newton (SPIN) method to solve the nonlinear systems arising in the phase-field fracture simulations in a monolithic manner.
Motivated by the robustness of the alternate minimization scheme, we have constructed the nonlinear preconditioner for the coupled problem utilizing the field split approach.
Thus, the SPIN preconditioner has been built by minimizing the underlying energy functional separately with respect to the displacement and the phase-field.
This minimization has been performed in an additive or multiplicative manner using inexact Newton's method.
The numerical performance and the robustness of the newly proposed SPIN methods have been demonstrated using five standard benchmark problems.
A comparison with a widely used alternate minimization scheme has been made and it has shown a significant reduction in the execution time, especially for the problems with the gradual crack propagation.
Moreover, we have also demonstrated that the speedup factor of the SPIN method grows with increasing mesh resolution.

The SPIN algorithms for solving the phase-field fracture problems can be extended in several ways.
For example, we plan to investigate the convergence behavior of the SPIN method and its efficiency using more complex phase-field fracture models, such as ones, which take into account anisotropy, or nonlinear constitutive laws.
We also aim to extend the proposed field-split SPIN method using the nesting domain decomposition approach.
In this particular case, a variant of the additive SPIN method, which would employ decomposition of the computational domain into smaller subdomains, could be employed to solve the phase-field and displacement subproblems.
This would enable an additional level of parallelism, which could in turn provide a further reduction of computational time.

  \appendix
  \section{Command-line options}
\label{sec:appendix_code_options}
The implementation of the phase-field fracture model, benchmark problems, and the solution strategies presented in this manuscript is freely available at~\url{https://bitbucket.org/alena_kopanicakova/pf_frac_spin}.
Following PETSc methodology, the developed solutions strategies can be configured using command-line options, which are summarized in \Cref{tab:command_line_options}.

\begin{table}
  \caption{The command line arguments used to configure AM/SPIN methods.}
  \label{tab:command_line_options}
  \centering
  \begin{tabular}{|r|c|l|}
    \hline
    Option                                        & Default value & Description                                            \\ \hline \hline
    \texttt{-snes\_atol}                          & $10^{-7}$     & Absolute tolerance (coupled residual).     \\ \hline
    \texttt{-snes\_rtol}                          & $10^{-6}$     & Relative tolerance (coupled residual).     \\ \hline
    \texttt{-snes\_atol}                          & $10^{-8}$     & Correction tolerance (coupled correction).        \\ \hline
    \texttt{-snes\_max\_it}                       & $50000$       & A maximum number of iterations.                        \\ \hline
\texttt{-snes\_am\_c\_diff\_tol}              & $10^{-4}$     & Absolute tolerance (change in the phase-field variable).   \\ \hline
    \texttt{-snes\_am\_disp\_diff\_tol}              & $10^{-12}$     & Absolute tolerance (change in the displacement variable).   \\ \hline
    \texttt{-snes\_am\_inexact\_solve}      	& true          & Usage of the inexact Newton's method.                  \\ \hline
    \texttt{-snes\_am\_direct\_solver}            & false         & Usage of the direct linear solver.          \\ \hline \hline
    \texttt{-snes\_spin\_additive}                & true          & Usage of additive variant of the SPIN method.          \\ \hline
    \texttt{-snes\_spin\_action\_rtol}            & $10^{-4}$     & Relative tolerance (action of SPIN operator).          \\ \hline
  \end{tabular}
\end{table}

\section{Strong-Wolfe conditions}
\label{sec:strong_wolfe_conditions}
Consider the following minimization problem
\begin{linenomath*}
  \begin{equation}
  \min_{\xv \in \R^n} f(\xv),
  \end{equation}
\end{linenomath*}
where $f:\R^n \mapsto \R$ is an objective function.
The k-th iteration of an iterative method embedded inside of a line-search framework is given as
\begin{linenomath*}
  \begin{equation}
    \xv^{(k+1)} = \xv^{(k)} + \alpha^{(k)} \pv^{(k)},
  \end{equation}
\end{linenomath*}
where $\pv^{(k)}$ denotes a descent search-direction.
In this work, the step size $\alpha^{(k)} \in \R^{+}$ is obtained by a cubic backtracking line-search method~\cite[Algorithm A6.3.1, pages 325 - 327]{dennis1996numerical}, such that it satisfies the following strong Wolfe conditions~\cite{nocedal2006numerical}:
\begin{linenomath*}
\begin{equation}
  \begin{aligned}
    f(\xv^{(k)} + \alpha^{(k)} \pv^{(k)})                                       & \leq f(x_k^{(k)}) + c_1 \alpha^{(k)}   \langle \nabla f(\xv^{(k)}), \pv^{(k)} \rangle \\
    | \langle \nabla f(\xv^{(k)} + \alpha^{(k)} \pv^{(k)}), \pv^{(k)} \rangle | & \leq c_2 |  \langle \nabla f(\xv^{(k)}), \pv^{(k)} \rangle |,
  \end{aligned}
\end{equation}
\end{linenomath*}
where $0 < c_1 < c_2 < 1$.
The numerical results reported in Section~\ref{sec:conv_study} employ $c_1=10^{-4}$ and $c_2=0.9$.

\section{Avoiding negative curvature}
The line-search methods require that the search direction $\pv^{(k)}$ is a descent search-direction~\cite{nocedal2006numerical}, i.e., $\langle \pv^{(k)}, \nabla f(\xv^{(k)}) \rangle < 0$.
As the underlying energy functional of the coupled phase-field fracture problem is non-convex, we can not guarantee that the search directions provided by the SPIN method are directions of descent.
Motivated by the trust-region strategies~\cite{Conn2000trust}, we employ a simple heuristic to ensure that the obtained search directions are always descent directions.
More precisely, if the search direction provided by the SPIN algorithm is a direction of ascent, then it is replaced with an inexact Newton's direction.
If also inexact Newton's direction is a direction of ascent, then we use a negative gradient direction instead.

Based on our empirical experience, this situation occurs rather rarely.
We encountered negative curvature only once for MSPIN and twice for ASPIN while simulating crack propagation of the asymmetrically notched beam.
For all other examples, the search directions provided by SPIN methods were directions of descent.

\end{sloppypar}

\section*{Acknowledgement}
The authors would like to thank the Swiss National Science Foundation and the Deutsche Forschungsgemeinschaft, Germany (DFG) for their support through the project SPP 1962 ``Stress-Based Methods for Variational Inequalities in Solid Mechanics:~Finite Element Discretization and Solution by Hierarchical Optimization [186407]".
Additionally, we would also like to gratefully acknowledge the support of Platform for Advanced Scientific Computing (PASC)  through projects FraNetG:~Fracture Network Growth and FASTER:~Forecasting and Assessing Seismicity and Thermal Evolution in geothermal Reservoirs.

 \bibliography{biblio_SPIN.bib}

\end{document}